\documentclass[leqno]{article}
\usepackage{latexsym,amsmath,amsfonts,mathrsfs,amssymb,amsthm}
\usepackage[usenames]{color}
\usepackage{graphicx}
\def\R{\mathbb R}
\def\N{\mathbb N}
\def\C{\mathbb C}

\def\s{\sharp}
\def\a{\alpha}
\def\e{\epsilon}
\def\d{\delta}
\def\Y{\mathbb Y}
\def\T{\mathbb T}
\def\P{\mathbb P}
\def\be{\begin{equation}}
\def\ee{\end{equation}}
\def\bs{\backslash}
\def\qed{\hfill$\Box$\bigskip}
\def\nd{\noindent Proof. }
\numberwithin{equation}{section}
\newtheorem{lem}[equation]{Lemma}
\newtheorem{pro}[equation]{Proposition}
\newtheorem{defn}[equation]{Definition}
\newtheorem{thm}[equation]{Theorem}
\newtheorem{cor}[equation]{Corollary}
\newtheorem{rem}[equation]{Remark}
\begin{document}

\begin{center} \Large\textbf{Almgren-minimality of unions of two almost orthogonal planes in $\R^4$}\end{center}

\bigskip

\centerline{\large Xiangyu Liang}

\vskip 1cm

\centerline {\large\textbf{Abstract.}}

\bigskip

In this article we prove that the union of two almost orthogonal planes in $\R^4$ is Almgren-minimal. This gives an example of a one parameter family of minimal cones, which is a phenomenon that does not exist in $\R^3$. This work is motivated by an attempt to classify the singularities of 2-dimensional Almgren-minimal sets in $\R^4$. Note that the traditional methods for proving minimality (calibrations and slicing arguments) do not apply here, we are obliged to use some more complicated arguments such as a stopping time argument, harmonic extensions, Federer-Fleming projections, etc. that are rarely used to prove minimality (they are often used to prove regularity). The regularity results for 2-dimensional Almgren minimal sets (\cite{DJT},\cite{DEpi}) are also needed here.

%. The key idea is to prove firstly that in the orthogonal case, the union of two planes is the unique minimal set verifying certain general condition, and then to use a stopping time argument to bring us to a situation where we can make some perturbation argument. The estimates of Hausdorff measures in this article use mainly the regularity results of 2-dimensional Almgren minimal sets (\cite{DJT},\cite{DEpi}), a projection argument for simple unit 2-vectors, and harmonic extensions on annuli. We use also a partial existence result for minimal sets in \cite{Fv}, to get necessary information so that we can begin the stopping time argument.

%The key point is to prove firstly that in orthogonal case, the union of two planes is the unique minimal set verifying certain general condition, and then to use a stopping time argument and estimates of Hausdorff measure by regularity results of 2-dimensional Almgren minimal sets (\cite{DJT},\cite{DEpi}), and harmonic extensions on annuli, to prove that a small perturbation of angles between the two planes preserves the minimality of their union. The proof uses also a partial existence result for minimal sets in \cite{Fv}. 

\bigskip


\bigskip


\section{Introduction and preliminaries}

\subsection{Introduction}

One of the main topics in geometric measure theory is the theory of minimal sets, currents and surfaces, which aims at important progress in understanding the regularity and existence of physical objects that have certain minimizing properties such as soap films. This is known in physics as Plateau's problem. Recall that a most simple version of Plateau's problem aims at finding a surface which minimizes area among all the surfaces having a given curve as boundary. See works of Besicovitch, Federer, Fleming, De Giorgi, Douglas, etc. for example.

Lots of notions of minimality have been introduced to modernize Plateau's problem, such as minimal surfaces, mass minimizing or size minimizing currents (see \cite{DPHa} for their definitions), varifold (c.f.\cite{Al66}). In this article, we will mainly use the notion introduced by F.Almgren \cite{Al76}, in a general setting of sets, and which gives a very good description of the behavior of soap films. Note that soap films are 2-dimensional objects, but a general definition of $d-$dimensional minimal sets in an open set $U\subset \R^n$ is not more complicated.

Intuitively, a $d-$dimensional minimal set $E$ in an open set $U\subset \R^n$ is a closed set $E$ whose $d-$dimensional Hausdorff measure could not be decreased by any local Lipschitz deformation. (See Section 1.2 for the  precise definition.)

\bigskip

The point of view here is very different from those of minimal surfaces and mass minimizing currents under certain boundary condition, which are more usually used. Comparing to the big number of results in the theory of mass minimizing currents, or classical minimal surface, in our case, very little results of regularity and existence are known. However, Plateau's problem is more like the study of size minimizing currents, for which much less results are known either (see \cite{DP09} for certain existence results). One can prove that the support of a size minimizing current is automatically an Almgren minimizer, so that all the regularity results listed below are also true for supports of size minimizing currents.

\medskip

First regularity results for minimal sets have been given by Frederick Almgren \cite{Al76} (rectifiability, Ahlfors regularity in arbitrary dimension), then generalized by Guy David and Stephen Semmes \cite{DS00} (uniform rectifiability, big pieces of Lipschitz graphs), Guy David \cite{GD03} (minimality of the limit of a sequence of minimizers). 

Since minimal sets are rectifiable and Ahlfors regular, they admit a tangent plane at almost every point. But our main interest is to study those points where there is no tangent plane, i.e. singular points.

\medskip

A first finer description of the interior regularity for minimal sets is due to Jean Taylor, who gave in \cite {Ta} an essential regularity theorem for 2-dimensional minimal sets in 3-dimensional ambient spaces: if $E$ is a minimal set of dimension 2 
in an open set of $\R^3$ , then  every point $x$ of $E$ has a neighborhood where $E$ is equivalent (modulo a negligible set) through a $C^1$ diffeomorphism to a minimal cone (that is, a minimal set which is also a cone).

In \cite{DJT}, Guy David generalized Jean Taylor's theorem to 2-dimensional minimal sets in $\R^n$, but with a local bi-H\"older equivalence, that is, every point $x$ of $E$ has 
a neighborhood where $E$ is equivalent through a bi-H\"older diffeomorphism to a minimal cone $C$ (but the minimal cone might not be unique). 

%This cone $C$ is obtained by blow-up limit of $E$ at $x$, i.e., there exists a sequence $r_n\to 0$ such that $C$ is the limit of the sequence of sets $r_n^{-1}(E-x)$. The set $E$ might admit more than one blow-up limit at a point $x\in E$, but all blow-up limit of $E$ at $x$ are bi-H\"older equivalent.

In addition, in \cite{DEpi}, David also proved that, if this minimal cone $C$ satisfies a ``full-length" condition, we will have the $C^1$ equivalence (called $C^1$ regularity). In particular, the tangent cone of $E$ at the point $x\in E$ exists and is a minimal cone, and the blow-up limit of $E$ at $x$ is unique; if the full-length condition fails, we might lose the $C^1$ regularity.

Thus, the study of singular points is transformed into the classification of singularities, i.e., into looking for a list of minimal cones. Besides, getting such a list would also help deciding locally what kind (i.e. $C^1$ or bi-H\"older) of equivalence with a minimal cone can we get.

In $\R^3$, the list of 2-dimensional minimal cones has been given by several mathematicians a century ago. (See for example \cite{La} or \cite{He}). 
%
%The idea is to find all nets on the unit sphere composed of great circles and arcs of great circles that meet by three and make angles of $120^\circ$. This is a condition satisfied by the intersection of all 2-dimensional minimal cones with the unit sphere (this is also true for 2-dimensional minimal cones in higher ambient dimension). We obtain a finite list of this kind of cones (in fact, 10 cones). Then we eliminate one by one those cones which are not minimal, by finding better competitors. In the end there are, 
They are, modulo isomorphism: a plane, a $\Y$ set (the union of 3 half planes that meet along a straight line where they make angles of 120 degrees), and a $\T$ set (the cone over the 1-skeleton of a regular tetrahedron centered at the origin). See the pictures below.

\centerline{\includegraphics[width=0.2\textwidth]{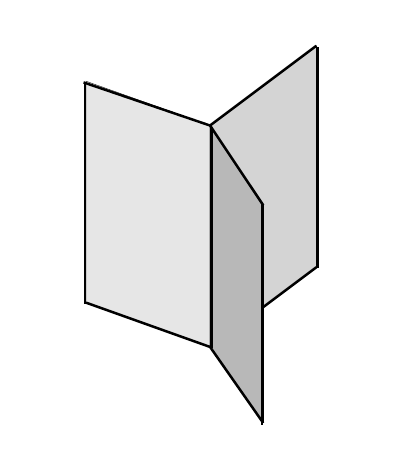} \hskip 2cm\includegraphics[width=0.25\textwidth]{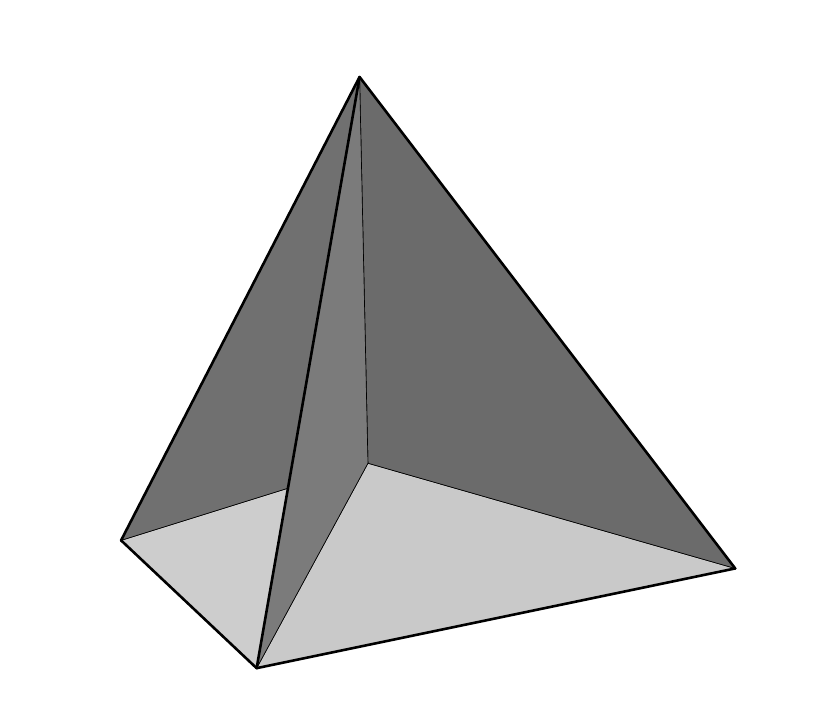}}

In higher dimensions, even in dimension 4, the list of minimal cones is still very far from clear. Except for those three minimal cones that already exists in $\R^3$, the only 2-dimensional minimal cone that was known before this paper is the union of two orthogonal planes, whose minimality could be proved by a simple projection argument.

%%On the one hand, every competitor of $E$ in $U$ is a competitor of $E$ in $U\times\R^m$; on the other hand, for any competitor $F$ of $E$ in $U\times\R^m$, its projection $\pi(F)$ into $U$ is also a competitor of $E$ in both $U$ and $U\times \R^m$, which has a smaller area, so if $E$ minimizes area among competitors in $U$, it minimizes area among competitors in $U\times \R^n$ automatically. So $E$ is minimal in $U$ if and only if $E$ is minimal in $U\times\R^m$. 
%Then, a simple projection argument gives the minimality of the union of two orthogonal planes. It is the first minimal cone which does not exist in $\R^3$. 

It was not even known whether this is the only minimal cone with this shape, i.e., as a union of two planes. 

Note that in $\R^3$, no small perturbation of any minimal cones ever preserves the minimality, and even, each minimal cone admits a different topology from the others. But in this article we will show the following theorem, which exhibits a phenomenon that does not exist in $\R^3$:

\begin{thm}There exists a continuous one-parameter family of 2-dimensional minimal cones $\R^n$ for $n\ge 4$.
\end{thm}

In fact, we are going to prove the following more precise theorem, of which Theorem 1.1 is a direct but noticeable corollary.

\begin{thm}[minimality of the union of two almost orthogonal planes]\label{main}There exists $0<\theta<\frac\pi2$, such that if $P^1$ and $P^2$ are two planes in $\R^4 $ whose characteristic angles $(\alpha_1,\alpha_2)$ satisfy $\a_2\ge\a_1\ge \theta$, then their union $P^1\cup P^2$ is a minimal cone in $\R^4$.
\end{thm}

The characteristic angles of two planes $P^1,P^2$ describe their relative position, the definition will be given in Section 2. But an equivalent statement of Theorem \ref{main} that might be easier to understand is that any almost orthogonal union of two planes in $\R^4$ is minimal, i.e. there exists $a>0$, such that if for any $v_1\in P^1,v_2\in P^2$,
\be |<v_1,v_2>|\le a||v_1||||v_2||\ee 
then $P^1\cup P^2$ is minimal. Or intuitively, there is an "open set" of unions of two planes, which contains the orthogonal union of two planes, such that each element in this set (which is a union of two planes) is a minimal cone.

When $a$ is small, (1.3) is almost equivalent to saying that $\alpha_1,\alpha_2>\frac\pi2-Ca$. Hence we really have a one parameter family of minimal cones (which are not $C^1$ equivalent), whose intersections with the unit ball have the same Hausdorff measure.

Moreover, it is not hard to check that each cone in this one parameter family verifies the "full-length property", hence by Thm 1.15 of \cite{DEpi}, if a blow-up limit of a minimal set $E$ at a point $x$ is one of these unions of planes $P_\a$, then locally $E$ is $C^1$ equivalent to $P_\a$.

This is also a general question for unions of higher dimensional planes. Unions of orthogonal planes are always minimal, but for non-orthogonal cases we know very little. Some new information can be given, using the arguments and results of this paper and the notion of topological minimal sets, see \cite{topo}.

\medskip

For the proof of the theorem, the argument is not simply a perturbation analysis. Due to the higher codimension and the lack of results for minimal sets in $\R^4$, the classical methods (such as calibration, slicing) for proving minimality do not apply, hence we are obliged to use  tools in harmonic analysis, which are rarely used to prove minimality results. In particular, we use a stopping time argument to divide minimal sets into good and bad parts, and we treat them differently (projection for the bad part, harmonic extension for the good part), to get the estimates for their Hausdorff measure. It would also be much simpler if we had the existence results for minimal sets, because this would allow us to do directly the stopping time argument. However we can only use a partial existence result of \cite{Fv}. We use it to get some sort of minimal set, and we have to prove some necessary regularity for it, so that the stopping time argument can be carried on well. Some of the arguments, especially the geometric constructions, are somehow painful. But we have not found any other easier way to do it. 

\begin{rem}

1) To the author's knowledge, the only proof of minimality that does not use a calibration argument (probably not available here) is the proof of minimality for the cone $T$, suggested by J.Taylor in \cite{Ta}, using a topological argument, and the list of all the other minimal cones. But in our case the topology is more complicated, existence results are weaker, and moreover the list of the other minimal cones are far from known.

2) Frank Morgan gave a conjecture in \cite{Mo93} on the angle condition (which is imposed on the sum of the two characteristic angles) under which a union of two planes is minimal. It was proved by Gary Lawlor \cite{La94} that if the sum of their two characteristic angles is less than $\frac{2\pi}{3}$, then the union of the two planes will not be minimal. The remaining part is still open. 

Our theorem also solves partially the remaining part of this conjecture, because if $\a_1+\a_2>\frac\pi 2+\theta$ ($\theta$ is as in Theorem \ref{main}), then we'll have $\a_2\ge \a_1>\theta$ (since both of them are no more than $\pi/2$), and Theorem \ref{main} implies that the union of two planes with this pair of angles is minimal.
\end{rem}

%Traditional tools for proving minimality are mostly calibrations and slicing arguments. For example, the same question about minimality of unions of planes has also been asked for mass minimizing current, and has has been well solved by Gary Lawlor \cite{La89} and Dana Nance \cite{Na}, where they proved the results by constructing some beautiful calibrations and surfaces. In particular a union of two planes in $\R^4$ is minimal if and only if the two planes are orthogonal, therefore in this different setting, the union of two orthogonal planes is the only minimizing one, thus the minimality is not stable under small perturbations. 
%
%However in our case, 

%Here we just give an idea of what happens. Those definitions ask that
%\be H^d(E\cap W_1)\le H^d(\varphi(E\cap W_1)+h(\d)\d^d,\ee
%for all deformation satisfying $diam(\widehat W)<\d$, where $\varphi_t,W_t,\widehat W$ are as in Definition 1.1, and we ask also that the gauge function $h$ tends to 0 as $\d$ tends to 0.
\bigskip

The strategy of the proof of Theorem \ref{main} is the following.

Since our objects are cones centered at the origin, to prove that their measure could not be decreased by any deformation in any compact set in $\R^n$ (which is the definition for minimal sets, see Subsection 1.1 for the precise definition), we can just look at deformations in the unit ball $B=B(0,1)$. For $\theta=(\theta_1,\theta_2)$, denote by $P_\theta$ the union of two planes with characteristic angles $\theta$. We begin by verifying that $P_{(\frac\pi2,\frac\pi2)}$ is minimal, and moreover is the only minimal set in the unit ball among all sets with the same boundary and surjective orthogonal projections on the two planes. Here the boundary of a set is just the intersection of its closure with the unit sphere.

Next we begin to prove Theorem \ref{main}. Suppose the conclusion of the theorem is false. Then there exists a sequence $\{\theta_k\}$ which converges to $\theta_0:=(\frac\pi2,\frac\pi2)$, such that $P_{\theta_k}$ is not minimal. Set $P_k=P_{\theta_k}$, and $P_0=P_0^1\cup_{\theta_0}P_0^2=P_0^1\cup_\perp P_0^2$. So we have
\be\inf\{H^2(F):\ F\subset\overline B(0,1)\ \mbox{is a deformation of }P_k\mbox{ in }B\}<H^2(P_k\cap B)=2\pi.\ee

The next step would be easier if we could find, for each $k$, a deformation $E_k$ of $P_k$ in $B$ which minimizes $H^2(E\cap B)$ among all deformations $E$ of $P_k$ in $B$. Unfortunately, no such existence theorem is known.

However, Vincent Feuvrier showed in his thesis \cite{Fv} a partial existence result, which says that given an original set $E$ in $\overline B$, there exists a certain set $F$ that is minimal in $B$, which is the limit of a minimizing sequence of deformations of $E$ in $B$, and such that
\be H^2(F)\le\inf\{H^2(E')\ :\ E' \subset\overline{B(0,1)}\mbox{ is a deformation of }E\mbox{ in }B(0,1)\}.\ee

We cannot say that this set $F$ obtained by limit of deformations is still a deformation of $E$, but luckily we shall not need to know that. We show first that for $i=1,2$, the projection of each set $E_k$ on the plane $P_k^i$ contains $P_k^i\cap B$; then we manage to show that $E_k$ has the same boundary as $P_k$. We can also get
\be H^2(E_k)<2\pi=H^2(P_k\cap B).\ee

Hence we get our favorite sequence of minimal sets $E_k$, whose boundaries converge to the boundary $P_0\cap\partial B$ of $P_0$ (recall that $P_0=P_0^1\cup_\perp P_0^2$). We take a convergent subsequence, that we still denote by $\{E_k\}$. By \cite{GD03} Thm 4.1, the limit of a sequence of minimal sets is still minimal, hence the limit $E_\infty$ of $\{E_k\}$ is minimal. Now $E_\infty\cap\partial B=P_0\cap\partial B$, hence the uniqueness theorem for $P_0$ tells us that $E_\infty$ cannot be anyone other than $P_0$. Moreover, the Hausdorff distance between $E_k$ and $P_k$ also converges to 0, because $P_k$ converges also to $P_0$. See Section 4 for detail.

Now for every fixed small $\e$, we will use a stopping time argument to decompose every $E_k$. More precisely, we use a stopping time argument to find, for each $k$, a critical radius $r_k=r_k(\e)$, such that $E_k$ is $\e$-near $P_k$ or a translation outside a ball $B(o_k,r_k)$, with $o_k$ fairly near the origin. Moreover in $B(o_k,r_k)$, $E_k$ begins to get away a little from $P_k$, that is, $r_k$ is the scale in which $E_k$ begins to go away from $P_k$, even though we do not understand how. In spite of that, we will manage to show that the the projections of the part of $E_k$ inside $B(o_k,r_k)$ on the two planes are surjective. (The proof of this projection property is essentially derived from the construction in the proof of the partial existence result of [9], which is somehow complicated).

We know that $r_k$ is the very first radius such that $E_k$ is not $\e r_k$ near $P_k$, so when $\e$ is sufficiently small, by the regularity of minimal sets near flat points, we manage to show that outside the small ball $B(o_k,\frac14 r_k)$, $E_k$ is the union of two $C^1$ graphs on $P_k^1,P_k^2$ respectively. Then we decompose $E_k$ into two parts, which are outside and inside $B(o_k,\frac14 r_k)$ respectively.

For the part outside, we already know that $E_k$ is composed of two graphs $G_i$ of $C^1$ functions $f_i$ from the annuli $P_k^i\cap B(0,1)\bs B(o_k,\frac14 r_k)$ to ${P_k^i}^\perp$ respectively. But the union of two parts is not $\e r_k$-near any translation of $P_k$ in $B(o_k,r_k)\bs B(o_k,\frac14 r_k)$ (in fact, we have only shown that $E_k$ is not $\e r_k$-near any translation of $P_k$ in $B(o_k,r_k)$ rather than $B(o_k,r_k)\bs B(o_k,\frac14 r_k)$, but then we can prove by a compactness argument that $E_k$ is not too close to any translation of $P_k$ in $B(o_k,r_k)\bs B(o_k,\frac14 r_k)$ either, maybe at the price of making $\e$ smaller. See Section 8 for more detail). Thus, at least one of the two graphs is $\frac 12\e r_k$ far from any translation of its domain. Suppose this is the case for $i=1$. This means that $f_1$ oscillates of order $\e r_k$ in the region $B(o_k,r_k)\bs B(o_k,\frac14 r_k)$. Then some argument using harmonic functions gives that the Dirichlet energy $\int_{P_k^1\cap B(0,1)\bs B(o_k,\frac14 r_k)}|\nabla f_1|^2 $ of $f_1$ on  $P_k^1\cap B(0,1)\bs B(o_k,\frac14 r_k)$ is larger than $C(\e)r_k^2$. But this is almost equivalent to the difference between the measure of the graph $G_1$ and the measure of the annulus $P_k^1\cap B(0,1)\bs B(o_k,\frac14 r_k)$, so we have
\be H^2(E_k\cap B(0,1)\bs B(o_k\frac 14 r_k))\ge H^2(P_k^1\cap B(0,1)\bs B(o_k,\frac14 r_k))+C(\e)r_k^2.\ee
This means we gain some measure $C(\e) r_k^2$, where the constant $C(\e)$ does not depend on $k$.

For the inside part, by a projection argument we can show that 
\be \begin{split}H^2(E_k\cap B(o_k,\frac 14r_k))&\ge (1- C\times (\frac\pi2-\theta_k))H^2(P_k\cap B(o_k,\frac14 r_k))\\
&=H^2(P_k\cap B(o_k,\frac14 r_k))-C\times(\frac\pi2-\theta_k)r_k^2,\end{split}\ee
where (recall that) $\frac\pi 2-\theta_k\to 0$ as $k\to\infty$.

Altogether we get
\be H^2(E_k)\ge H^2(P_k\cap B(0,1))+[C(\e)-C\times (\frac\pi2-\theta_k)]r_k^2.\ee
When $k\to\infty$, the term $[C(\e)-C\times (\frac\pi2-\theta_k)]r_k^2$ is strictly positive, so the measure of $E_k\cap B$ is strictly larger than $2\pi$. But this cannot happen, because of (1.7). 

Thus we finish the proof of Theorem \ref{main}.

Here we have to point out that if we could control how fast are the $E_k$ converging to $P_0$, then we would be able to give an estimate on the term $C(\e)-C\times (\frac\pi2-\theta_k)$, and thus to give an estimate on the angle $\theta$ in Theorem \ref{main}. But such a control might need a deeper understanding of the uniqueness theorem of $P_0$ (Thm \ref{unicite}).

The same kind of argument for proving Theorem \ref{main} can be generalized (but not trivially) to prove the minimality of the union of $n$ almost orthogonal planes of dimension $m$ in $\R^{mn}$ (which are also related to Morgan's conjecture), and also the union of a plane and a $\Y$ in $\R^5$. See \cite{XY10} for detail.

\bigskip

The plan for the rest of this article is the following. 

Subsection 1.1 will give some notation and conventions that we will use frequently afterwards.

Section 2 is devoted to the estimation of the sum of projections of a simple unit 2-vector to two planes, estimation that depends on the characteristic angles of the two planes. Based on this we give a comparison between the measure of a rectifiable set and the sum of the measures of its projections on the two planes.

In Section 3 we prove that $P_0$ is the only minimal set in $B(0,1)$ with the given boundary and surjective projections.

In Section 4 we show the existence of the minimal sets $E_k$, as well as some of their properties.

Section 5 is devoted to finding the critical radius $r_k$, and giving some properties of $E_k$ outside the``ball" $D(o_k,\frac14 r_k)$.

In Section 6 we prove that the projections of $E_k$ on $P_k^1$ and $P_k^1$ are surjective, in $B(0,1)$, as well as in $B(o_k,t)$ for all $t\in[\frac14r_k,r_k]$. We also give the $C^1$ regularity of $E_k$ outside the ball $B(o_k,\frac14 r_k)$.

Section 7 contains an argument of harmonic extension, which gives a lower bound for the measure of the graph of a function, depending on the size of the order of its oscillations near the boundary.

Finally we arrive at our conclusion in Section 8, by combining all the information we gathered before.

%Sections 9-11 are devoted to generalizations of Theorem \ref{main}. In Section 9 we prove the minimality of the union of two almost orthogonal planes of dimension $m$; in Section 10 we prove the minimality of the union of $n$ almost orthogonal $m-$planes; in Section 11 we prove the minimality of the almost orthogonal union of a plane and a $\Y$ set in $\R^5$. Amusingly, we do not know how to treat the almost orthogonal union of two $\Y$ sets in $\R^6$, although we know that their orthogonal union is minimal.

\textbf{Acknowledgement:} I would like to thank Guy David for many helpful discussions and for his continual encouragement. I also wish to thank Thierry De Pauw for having given many useful suggestions. The results in this paper were part of the author's doctoral dissertation at University of Paris-Sud 11 (Orsay). This work was supported by grants from R\'egion Ile-de-France.

\subsection{Preliminaries}

$\bullet$ In all that follows, minimal set means Almgren minimal set;

\noindent$\bullet$ $[a,b]$ is the line segment with end points $a$ and $b$;

\noindent$\bullet$ $[a,b)$ is the half line issued from the point $a$ and passing through $b$;

\noindent$\bullet$ $B(x,r)$ is the open ball with radius $r$ and centered on $x$;

\noindent$\bullet$ $\overline B(x,r)$ is the closed ball with radius $r$ and center $x$;

\noindent$\bullet$ For any set $E\subset\R^n$, and $r>0$, $B(E,r):=\{x\in \R^n, d(x,E)<r\}=\cup_{x\in E}B(x,r)$;

\noindent$\bullet$ $\overrightarrow{ab}$ is the vector $b-a$;

\noindent$\bullet$ $H^d$ is the Hausdorff measure of dimension $d$ ;

\noindent$\bullet$ $d_H(E,F)=\max\{\sup\{d(y,F):y\in E,\sup\{d(y,E):y\in F\}\}$ is the Hausdorff distance between two sets $E$ and $F$.

%\noindent$\bullet$ $d_{x,r}$ is the relative distance with respect to the ball $B(x,r)$, defined as
%$$ d_{x,r}(E,F)=\frac 1r\max\{\sup\{d(y,F):y\in E\cap B(x,r)\},\sup\{d(y,E):y\in F\cap B(x,r)\}\}.$$

\medskip

In the next definitions, fix integers $0<d<n$. We are mostly interested in $d=2$ and $n=4$ here.

\begin{defn}[Almgren competitor (Al competitor for short)] Let $E$ be a closed set in an open subset $U$ of $\R^n$. An Almgren competitor for $E$ is a closed set $F\subset \overline U$ that can be written as $F=\varphi_1(E)$, where $\varphi_t:U\to U$ is a family of continuous mappings such that 
\be \varphi_0(x)=x\mbox{ for }x\in U;\ee
\be\mbox{ the mapping }(t,x)\to\varphi_t(x)\mbox{ of }[0,1]\times U\mbox{ to }U\mbox{ is continuous;}\ee
\be\varphi_1\mbox{ is Lipschitz,}\ee
  and if we set $W_t=\{x\in U\ ;\ \varphi_t(x)\ne x\}$ and $\widehat W=\bigcup_{t\in[0.1]}[W_t\cup\varphi_t(W_t)]$,
then
\be \widehat W\mbox{ is relatively compact in }U.\ee
 
Such a $\varphi_1$ is called a deformation in $U$, and $F$ is also called a deformation of $E$ in $U$.
\end{defn}

\begin{defn}[(Almgren) minimal sets]
Let $0<d<n$ be integers, $U$ an open set of $\R^n$. A closed set $E$ in $U$ is said to be minimal of dimension $d$ in $U$ if 
\be H^d(E\cap B)<\infty\mbox{ for every compact ball }B\subset U,\ee
and
\be H^d(E\bs F)\le H^d(F\bs E)\ee
for all Al competitors $F$ for $E$.
\end{defn}

\begin{rem}When the ambient set $U$ is $\R^n$, or a ball, we can also take the class of local Almgren competitors to define the same notion of minimal set. Keep the $E$, $U$, $n$ and $d$ as before; a local Almgren competitor of $E$ in $U$ is a set $F=f(E)$, with
\be f=id\mbox{ outside some compact ball }B\subset U,\ee
\be f(B)\subset B,\ee 
and $f$ is Lipschitz. 

A such $f$ is called a local deformation in $U$, or a deformation in $B$, and $F=f(E)$ is also called a local deformation of $E$ in $U$, or a deformation of $E$ in $B$.

Note that in this case, the condition (1.18) becomes 
\be H^d(E\cap B)\le H^d(F\cap B).\ee

We say that a set $E$ closed in an open set $U$ is locally minimal if (1.17) holds, and for any compact ball $B\subset U$, and any local Almgren competitor $F$ for $E$ in $B$, (1.22) holds.

One can easily verify that when $U$ is $\R^n$ or a ball, the class of Al competitors coincides with the class of local Al competitors, so the two classes define the same kind of minimal sets. However, if the ambient set $U$ has a more complicated geometry, then the class of local Al competitors is strictly smaller, so a set minimizing the Hausdorff measure among local Al competitors might fail to be Al-minimal. 
\end{rem}

\begin{rem}The notion of minimal sets does not depend much on the ambient dimension. One can easily check that $E\subset U$ is $d-$dimensional minimal in $U\subset \R^n$ if and only if $E$ is minimal in $U\times\R^m\subset\R^{m+n}$, for any integer $m$.\end{rem}

\begin{rem}
For a description of a little more general kind of sets, such as soap bubbles (where the pressures are different in different connected components of their complementary), or taking gravity into consideration, we can use a notion of almost minimal sets, where we add some error term in the condition (1.18). To these sets, almost all the argument for minimal sets can be applied, and we have essentially the same regularity results. See \cite{DJT} for definitions and more detail.
\end{rem}

\begin{defn}[reduced set] Let $U\subset \R^n$ be an open set. For every closed subset $E$ of $U$, denote by
\be E^*=\{x\in E\ ;\ H^d(E\cap B(x,r))>0\mbox{ for all }r>0\}\ee
 the closed support (in $U$) of the restriction of $H^d$ to $E$. We say that $E$ is reduced if $E=E^*$.
\end{defn}

It is easy to see that
\be H^d(E\bs E^*)=0.\ee
In fact we can cover $E\bs E^*$ by countably many balls $B_j$ such that $H^d(E\cap B_j)=0$.

\begin{rem}If $E$ is minimal, then $E^*$ is also minimal, because if $\varphi_1$ is a deformation in $U$ as defined in Definition 1.11, then it is Lipschitz, hence $H^d(\varphi(E\bs E^*))=H^d(E\bs E^*)=0$. So the condition (1.18) is the same for $E^*$ as for $E$. As a result it is enough to study reduced minimal sets.
\end{rem}

\begin{defn}[blow-up limit] Let $U\subset\R^n$ be an open set, let $E$ be a relatively closed set in $U$, and let $x\in E$. Denote by $E(r,x)=r^{-1}(E-x)$. A set $C$ is said to be a blow-up limit of $E$ at $x$ if there exists a sequence of numbers $r_n$, with $\lim_{n\to \infty}r_n=0$, such that the sequence of sets $E(r_n,x)$ converges to $C$ for the Hausdorff distance in any compact set of $\R^n$.
\end{defn}

\begin{rem}A set $E$ might have more than one blow-up limit at a point $x$.
\end{rem}

\medskip

We now state some regularity results that will be used throughout this paper.

\begin{defn}[bi-H\"older ball for closed sets] Let $E$ be a closed set of Hausdorff dimension 2 in $\R^n$. We say that $B(0,1)$ is a bi-H\"older ball for $E$, with constant $\tau\in(0,1)$, if we can find a 2-dimensional minimal cone $Z$ in $\R^n$ centered at 0, and $f:B(0,2)\to\R^n$ with the following properties:

$1^\circ$ $f(0)=0$ and $|f(x)-x|\le\tau$ for $x\in B(0,2);$

$2^\circ$ $(1-\tau)|x-y|^{1+\tau}\le|f(x)-f(y)|\le(1+\tau)|x-y|^{1-\tau}$ for $x,y\in B(0,2)$;

$3^\circ$ $B(0,2-\tau)\subset f(B(0,2))$;

$4^\circ$ $E\cap B(0,2-\tau)\subset f(Z\cap B(0,2))\subset E.$  

We also say that B(0,1) is of type $Z$.

We say that $B(x,r)$ is a bi-H\"older ball for $E$ of type $Z$ (with the same parameters) when $B(0,1)$ is a bi-H\"older ball of type $Z$ for $r^{-1}(E-x)$.
\end{defn}

\begin{thm}[Bi-H\"older regularity for 2-dimensional minimal sets, c.f.\cite{DJT} Thm 16.1]\label{holder} Let $U$ be an open set in $\R^n$ and $E$ a reduced minimal set in $U$. Then for each $x_0\in E$ and every choice of $\tau\in(0,1)$, there is an $r_0>0$ and a minimal cone $Z$ such that $B(x_0,r_0)$ is a bi-H\"older ball of type $Z$ for $E$, with constant $\tau$. Moreover, $Z$ is a blow-up limit of $E$ at $x$.
\end{thm}

\begin{defn}[point of type $Z$] In the above theorem, we say that $x_0$ is a point of type $Z$ (or $Z$ point for short) of the minimal set $E$.
\end{defn}

\begin{rem} Again, since we might have more than one blow-up limit for a minimal set $E$ at a point $x_0\in E$, the point $x_0$ might of more than one type (but all blow-up limits at a point are bi-H\"older equivalent). However, if one of the blow-up limits of $E$ at $x_0$ admits the "full-length" property (see Remark \ref{ful}), then in fact $E$ admits a unique blow-up limit at the point $x_0$. Moreover, we have the following $C^{1,\a}$ regularity around the point $x_0$. In particular, the blow-up limit of $E$ at $x_0$ is in fact a tangent cone of $E$ at $x_0$.
\end{rem}

\begin{thm}[$C^{1,\a}-$regularity for 2-dimensional minimal sets, c.f. \cite{DEpi} Thm 1.15]\label{c1} Let $E$ be a 2-dimensional reduced minimal set in the open set $U\subset\R^n$. Let $x\in E$ be given. Suppose in addition that some blow-up limit of $E$ at $x$ is a full length minimal cone (see Remark \ref{ful}). Then there is a unique blow-up limit $X$ of $E$ at $x$, and $x+X$ is tangent to $E$ at $x$. In addition, there is a radius $r_0>0$ such that, for $0<r<r_0$, there is a $C^{1,\a}$ diffeomorphism (for some $\a>0$) $\Phi:B(0,2r)\to\Phi(B(0,2r))$, such that $\Phi(0)=x$ and $|\Phi(y)-x-y|\le 10^{-2}r$ for $y\in B(0,2r)$, and $E\cap B(x,r)=\Phi(X)\cap B(x,r).$ 

We can also ask that $D\Phi(x)=Id$. We call $B(x,r)$ a $C^1$ ball for $E$ of type $X$.
\end{thm}

\begin{rem}[full length, union of two full length cones $X_1\cup X_2$]\label{ful}We are not going to give the precise definition of the full length property. Instead, we just give some information here, which is enough for the proofs in this paper.

1) The three types of 2-dimensional minimal cones in $\R^3$, i.e. the planes, the $\Y$ sets, and the $\T$ sets, all verify the full-length property (c.f., \cite{DEpi} Lemmas 14.4, 14.6 and 14.27). Hence all 2-dimensional minimal sets $E$ in an open set $U\subset\R^3$ admits the local $C^{1,\a}$ regularity at every point $x\in E$. But this was known from \cite{Ta}.

2) (c.f., \cite{DEpi} Remark 14.40) Let $n>3$. Note that the planes, the $\Y$ sets and the $\T$ sets are also minimal cones in $\R^n$. Denote by $\mathfrak C$ the set of all planes, $\Y$ sets and $\T$ sets in $\R^n$. Let $X=\cup_{1\le i\le n}X_i\in \R^n$ be a minimal cone, where $X_i\in \mathfrak{C}, 1\le i\le n$, and for any $i\ne j$, $X_i\cap X_j=\{0\}$. Then $X$ also verifies the full-length property. 
\end{rem}

\bigskip

\section{Projections on two orthogonal or almost orthogonal planes}

In this section, we will give some estimates for the sum of the measures of projections of a rectifiable set on two transversal planes. These estimates are somewhat algebraic, and mainly use estimates for the sum of projections of simple unit 2-vectors in $\R^4$. Here unit simple 2-vectors are used to represent planes and their relative positions, or are treated as surface elements or derivatives of functions between two rectifiable sets when we try to do some integration. As a result the orientation of a simple 2-vector will be ignored, in other words, we will essentially not need to distinguish between $x\wedge y$ and $y\wedge x$.

Denote by $\wedge_2(\R^4)$ the space of all 2-vectors in $\R^4$. Let $x,y$ be two vectors in $\R^4$, we denote by $x\wedge y\in\wedge_2(\R^4)$ their exterior product. If $\{e_i\}_{1\le i\le 4}$ is an orthonormal basis, then $\{e_i\wedge e_j\}_{1\le i<j\le 4}$ forms a basis of $\wedge_2 (\R^4)$. We say that an element $v\in\wedge_2 (\R^4)$ is simple if it can be expressed as the exterior product of two vectors. 

The norm on $\wedge_2 (\R^4)$, denoted by $|\cdot |$, is defined by 
\be|\sum_{i<j}\lambda_{ij}e_i\wedge e_j|=\sum_{i<j}|\lambda_{ij}|^2.\ee
Under this norm $\wedge_2(\R^4)$, is a Hilbert space, and $\{e_i\wedge e_j\}_{1\le i<j\le 4}$ is an orthonormal basis. For a simple 2-vector $x\wedge y$, its norm is 
\be |x\wedge y|=||x||||y||\sin\angle_{x,y},\ee
where $\angle_{x,y}\in[0,\pi]$ is the angle between the vectors $x$ and $y$, and $||\cdot ||$ denotes the Euclidean norm on $\R^4$. A unit simple 2-vector is a simple 2-vector of norm 1. Notice that $|\cdot|$ is generated by the scalar product $<,>$ defined on $\wedge_2(\R^4)$ as follows: for $\xi=\sum_{1\le i<j\le 4} a_{ij}e_i\wedge e_j, \zeta=\sum_{1\le i<j\le 4} b_{ij}e_i\wedge e_j$, 
\be <\xi,\zeta>=\sum_{1\le i<j\le 4} a_{ij}b_{ij}.\ee

One can easily verify that if two pairs of vectors $x,y$ and $x',y'$ generate the same 2-dimensional subspace of $\R^4$, then there exists $r\in\R\bs\{0\}$ such that $x\wedge y=rx'\wedge y'$.

Now given a unit simple 2-vector $\xi$, we can associate it to a 2-dimensional subspace $P(\xi)\in G(4,2)$, where $G(4,2)$ denotes the set of all 2-dimensional subspaces of $\R^4$:
\be P(\xi)=\{v\in\R^4,v\wedge\xi=0\}.\ee 
In other words, $P(x\wedge y)$ is the subspace generated by $x$ and $y$.

From time to time, when there is no ambiguity, we write also $P=x\wedge y$, where $P\in G(4,2)$ and the two unit vectors $x,y\in\R^4$ are such that $P=P(x\wedge y)$. In this case $x\wedge y$ represents a plane.

For the side of linear maps, if $f$ is a linear map from $\R^4$ to $\R^4$, then we denote by $\wedge_2f$ (and sometimes $f$ if there is no ambiguity) the linear map from $\wedge_2(\R^4)$ to $\wedge_2(\R^4)$ such that
\be \wedge_2f(x\wedge y)=f(x)\wedge f(y).\ee

And for the side of $G(4,2)$ (the set of all planes, without considering orientations), for a unit simple 2-vector $\xi\in\wedge_2\R^4$, we have always $P(\xi)=P(-\xi)$, so that we can define $|f(\cdot)|:G(4,2)\to\R^+\cup\{0\}$ by
\be |f(P(\xi))|=|\wedge_2f(\xi)|.\ee
One can easily verify that the value of  $|f(P(\xi))|$ does not depend on the choice of $\xi$ that generates $P$. So $|f(\cdot)|$ is well defined.

\subsection{Some estimate for the sum of projections of simple 2-vectors on two transversal planes}

Let us recall the definition of characteristic angles between two planes. Let $P^1,P^2$ be two 2-dimensional planes in $\R^4$. Among all pairs of unit vectors $(v,w)$ with $v\in P^1,\ w\in P^2$, we choose $(v_1,w_1)$ which minimizes the angle between them. We denote by $\a_1$ this angle. Next we look at all the pairs of unit vectors $\{(v',w'):v'\in P^1,w'\in P^2,v'\perp v_1,w'\perp w_1\}$ (here $P^1$ and $P^2$ are 2-dimensional, so once $w_1,v_1$ are fixed, generally we only have four such pairs, $(\pm v_2,\pm w_2)$), and we choose $(v_2,w_2)$ which minimize the angle among all such pairs. Denote by $\a_2$ this angle. Then  $(\a_1,\a_2)$ (with $\a_1\le\a_2$) is the pair of characteristic angles between $P^1$ and $P^2$.

Characteristic angles characterize absolutely the relative position between two planes, in the sense that we can find an orthonormal basis $\{e_i\}_{1\le i\le 4}$ of $\R^4$, such that
\be P^1=e_1\wedge e_2\mbox{ and }P^2=(\cos\a_1 e_1+\sin\a_1 e_3)\wedge (\cos\a_2 e_2+\sin\a_2 e_4).\ee

Notice that two planes are orthogonal if their pair of characteristic angles is $(\frac\pi2,\frac\pi2)$.

Now we want to estimate the sum of projections of a unit simple 2-vector on them. 
Denote by $p^i$ the orthogonal projection on $P^i,i=1,2$. Then $p^i$ is a linear map. We also denote by $p^i$ the linear map $\wedge_2p^i$ from $\wedge_2(\R^4)$ to itself, as defined in (2.5).

Now let $\xi$ be a unit simple 2-vector, then there exists two unit vectors $x,y$ such that $\xi=x\wedge y$ and $x\perp y$. We write $x,y$ in the basis : 
\be x=ae_1+be_2+ce_3+de_4,y=a'e_1+b'e_2+c'e_3+d'e_4\ee
with
\be a^2+b^2+c^2+d^2=a'^2+b'^2+c'^2+d'^2=1\ee
and
\be aa'+bb'+cc'+dd'=0.\ee

We calculate the projections $|p^j(\xi)|$.
\be |p^1(\xi)|=|<e_1\wedge e_2,\xi>|=|ab'-a'b|,\ee
and
\be \begin{split}|p^2(\xi)|=&\ |< (\cos\a_1 e_1+\sin\a_1 e_3)\wedge (\cos\a_2 e_2+\sin\a_2 e_4),\xi>|\\
=&\ |(ab'-a'b)\cos\a_1\cos\a_2+(ad'-a'd)\cos\a_1\sin\a_2\\
&+(cb'-c'b)\sin\a_1\cos\a_2+(cd'-c'd)\sin\a_1\sin\a_2|.
\end{split}\ee
Then when $\a_1=\a_2=\frac\pi 2$, we have
\be \begin{split}
|p^1(\xi)|+|p^2(\xi)|&=|ab'-a'b|+|cd'-c'd|\le |ab'|+|a'b|+|cd'|+|c'd|\\
&\le\frac12(a^2+b'^2+a'^2+b^2)+\frac12(c^2+d'^2+c'^2+d^2)=1\end{split}\ee
because of (2.9). So we get the following lemma.

\begin{lem}Let $P^1,P^2$ be two orthogonal planes, then for every unit simple 2-vector $\xi\in\wedge_2\R^4$ we have 
\be|p^1(\xi)|+|p^2(\xi)|\le 1.\ee
\end{lem}

More precisely, the next lemma gives exactly which are those unit simple 2-vectors satisfying equality in (2.15). Denote by $\Xi$ the set of all unit simple 2-vectors $\xi\in\bigwedge_2\R^4$ such that
\be |p^1(\xi)|+|p^2(\xi)|=1.\ee Then $P(\Xi)=\{P(\xi),\xi\in\Xi\}$ is a compact subset of $G(4,2)$ (c.f. \cite{Fe},1.6.2).
\begin{lem}If $x\wedge y\in\Xi$, then there exists $\a\in[0,\frac\pi2]$, $v_i,u_i,i=1,2$ four unit vectors such that $v_i\in P^1,\ u_i\in P^2$, $v_1\perp v_2,u_1\perp u_2$, so that
\be x=\cos\a v_1+\sin\a u_1\mbox{ and }y=\cos\a v_2+\sin\a u_2.\ee
\end{lem}

\nd This is just Wirtinger's inequality stated in 1.8.2 of \cite{Fe}, with $\nu=2$, $\R^4=\C_1\oplus\C_2$, $P^1=\C_1$, ${P^1}^\perp=\C_2$, $\mu=1$; the argument consists in checking the equality cases in (2.13). 
\qed

Now we look at unions of two almost orthogonal planes. 

\begin{pro} Let $0\le\alpha_1\le \alpha_2\le\frac\pi2$, and let $P^1,P^2\subset\R^4$ be two planes with characteristic angles $\alpha_1\le \alpha_2$. Denote by $p^i$ the orthogonal projection on $P^i$, $i=1,2$. Then for any unit simple 2-vector $\zeta\in \bigwedge_2\R^4$, its projections on these two planes satisfy:
\be |p^1\zeta|+|p^2\zeta|\le 1+2\cos\alpha_1.\ee
\end{pro}

\begin{rem}
Notice that when $\alpha_1$ tends to $\frac\pi2$, $\cos\alpha_1$ tends to $0$. So Proposition 2.19 implies that for $\e$ small, there exists $\theta=\theta(\e)\in]0,\frac\pi2[$ such that if
$\alpha_2\ge \alpha_1\ge\theta$, then for all unit simple vectors $\zeta\in\bigwedge_2\R^4$,
\be |p^1\zeta|+|p^2\zeta|\le 1+\e.\ee 
\end{rem}

\noindent Proof of Proposition 2.19.

Let $0\le\alpha_1\le\alpha_2\le\frac\pi2$ be two arbitrary angles, and let $P^1$ and $P^2$ be a pair of planes with characteristic angles $\alpha_1,\alpha_2$. Then there exists an orthonormal basis $\{e_i\}_{1\le i\le 4}$ of $\R^4$ such that $P^1=e_1\wedge e_2,\ P^2=(\cos\a_1 e_1+\sin\a_1 e_3)\wedge (\cos\a_2 e_2+\sin\a_2 e_4)$ .

Denote by $p$ the projection on the plane $P(e_3\wedge e_4)$. For each unit simple $\zeta\in\bigwedge^2(\R^4)$, we have 
\be|p^1\zeta|+|p^2\zeta|=|p^1\zeta|+|p\zeta|+(|p^2\zeta|-|p\zeta|)\le |p^1\zeta|+|p\zeta|+|p^2\zeta-p\zeta|.\ee
By Lemma 2.14, 
\be|p^1\zeta|+|p\zeta|\le 1.\ee
Hence we just need to estimate $|(p^2-p)\zeta|$. By definition,
\be\begin{split}
(p^2-p)(\zeta)=&<(\cos\a_1e_1+\sin\a_1e_3)\wedge(\cos\a_2e_2+\sin\a_2e_4)-e_3\wedge e_4,\zeta>\\
=&<\cos\a_1\cos\a_2e_1\wedge e_2+\cos\a_1\sin\a_2 e_1\wedge e_4+\sin\a_1\cos\a_2 e_3\wedge e_2\\
&+(\sin\a_1\sin\a_2-1)e_3\wedge e_4),\zeta>.
\end{split}\ee

Notice that $\zeta$ is a unit 2-vector, so we have $\zeta=\sum_{1\le i<j\le 4}a_{ij}e_i\wedge e_j,$ where $\sum_{1\le i<j\le 4}a_{ij}^2=1.$ Hence
\be \begin{split}&|(p^2-p)(\zeta)|\\
=&|a_{12}\cos\a_1\cos\a_2+a_{14}\cos\a_1\sin\a_2-a_{23}\sin\a_1\cos\a_2+a_{34}(\sin\a_1\sin\a_2-1)|\\
\le &[\cos^2\a_1\cos^2\a_2+\cos^2\a_1\sin^2\a_2+\sin^2\a_2\cos^2\a_2+(\sin\a_1\sin\a_2-1)^2]^\frac12\\
&\times[a_{12}^2+a_{14}^2+a_{23}^2+a_{34}^2]^\frac12\\
%\le & [\cos^2\a_1\cos^2\a_2+\cos^2\a_1\sin^2\a_2+\sin^2\a_2\cos^2\a_2+(\sin\a_1\sin\a_2-1)^2]^\frac12\\
\le & [\cos^2\a_1+\cos^2\a_1+\cos^2\a_2+(1-\sin^2\a_1)^2]^\frac12\\
\le & [3\cos^2\a_1+\cos^4\a_1]^\frac12\le[4\cos^2\a_1]^\frac12=2\cos\a_1\end{split}\ee
by Cauchy-Schwarz and $\a_1\le\a_2$.

Combining (2.26), (2.23) and (2.24) we obtain the conclusion.\qed

\subsection{Comparison of the measure of a rectifiable set with the sum of measures of its projections on two planes}

We will apply the estimates on simple 2-vectors obtained in the last subsection to rectifiable sets. Let $F$ be a 2-rectifiable set. Then for almost all $x\in F$ the approximate tangent plane of $F$ at $x$ exists (c.f.\cite{Ma} Thm 15.11), and we denote it by $T_xF$. Then $T_xF\in G(4,2)$. For each linear map $f:\R^4\to\R^4$, $|f(T_xF)|$ is defined as in (2.6).

\begin{lem}\label{projection}Let $P^1, P^2$ be two planes in $\R^4$, let $F\subset\R^4$ be a 2-rectifiable set. Denote by $p^i$ the projection on $P^i$. If $\lambda$ is such that for almost all
$x\in F$, the approximate tangent plane of $F$ at $x$ $T_xF\in G(4,2)$ is such that
\be |p^1(T_xF)|+|p^2(T_xF)|\le \lambda,\ee 
then
\be H^2(p^1(F))+H^2(p^2(F))\le \lambda H^2(F).\ee
\end{lem}

\nd

Denote by $f$ the restriction of $p^1$ on $F$, then $f$ is a Lipschitz function from a 2-rectifiable set to  a subset of a plane $P^1$. Since $F$ is 2-rectifiable, for $H^2-$ almost all $x\in F$, $f$ has an approximate differential
\be ap Df(x): T_xF\to P^1 \ee (c.f.\cite{Fe}, Thm 3.2.19). Moreover this differential is such that $||\bigwedge_2 ap Df(x)||\le 1$ almost everywhere, because $f$ is $1-$Lipschitz.

Now we can apply the area formula to $f$, (c.f. \cite{Fe} Cor 3.2.20). For all $H^2|_F$-integrable functions $g\ : \ F\to\bar\R$, we have
\be \int_F (g\circ f)\cdot ||\wedge_2 apDf(x)||dH^2=\int_{P^1} g(z)N(f,z)dH^2z,\ee where $ N(f,z)=\s\{f^{-1}(z)\},$
and for $z\in p^1(F)$ we have $N(f,z)\ge 1$.
Take $g\equiv1$, we get
\be\int_F ||\wedge_2 apDf(x)||dH^2=\int_{P^1} N(f,z)dH^2z\ge \int_{p^1(F)}dH^2=H^2(p^1(F)).\ee

Recall that $p^1$ is linear, hence its differential is itself. As a result $apDf(x)$ is the restriction of $p^1$ on the 2-subspace $T_xF$, which implies that if $\{u,v\}$ is an orthonormal basis of $T_xF$, then
\be ||\wedge_2apDf(x)||=|\wedge_2p^1(u\wedge v)|=|p^1(T_xF)|\ee
by (2.6).
 Hence by (2.32)
\be\int_F |p^1(T_xF)| dH^2(x)\ge H^2(p^1(F)).\ee
A similar argument gives also:
\be\int_F |p^2(T_xF)| dH^2(x)\ge H^2(p^2(F)).\ee
Summing over $i=1,2$ we get
\be\begin{array}{ll} H^2(p^1 F)+H^2(p^2 F)&\le\int_F |p^1 T_xF|+|p^2 T_xF|dH^2(x)\\
&\le\int_F\ \lambda\ dH^2(x)=\lambda H^2(F)\end{array}\ee
since $|p^1 T_xF|+|p^2 T_xF|\le \lambda$. \qed

Here are two corollaries of Lemma \ref{projection}.

\begin{cor}The set $P_0=P_0^1\cup_\perp P_0^2 \subset \R^4$ is minimal. 
\end{cor}

\nd

Let $F$ be an Al competitor of $P_0$ in $B(0,1)$. By Remark 1.19, this means that there exists a Lipschitz deformation $\varphi$ in $\R^4$, with 
\be\varphi|_{B(0,1)^C}=Id,\ \varphi(B(0,1))\subset B(0,1), \mbox{ and }\varphi(P_0)\cap B(0,1)=F.\ee 
We will compare the measure of $F$ with that of $P_0\cap B(0,1)$.

Denote by $p_0^i,i=1,2$ the projections of $\R^4$ on $P_0^i$. Since $F$ is 2-rectifiable, the approximate tangent plane $T_xF$ of $F$ at $x$ exists for almost all $x\in F$. By representing planes by unit simple 2-vectors, and by Lemma 2.14,  
\be |p_0^1 (T_xF)|+|p_0^2 (T_xF)|\le 1.\ee

As a result, Lemma \ref{projection} gives
\be H^2(p_0^1(F))+H^2(p_0^2(F))\le H^2(F).\ee

Denote by $F^i=p_0^i\circ \varphi (P^i_0\cap B(0,1)),i=1,2$, then $F^i$ is a deformation of $P^i_0\cap B(0,1)$, and by the minimality of planes (in arbitrary ambient dimension) we have
\be H^2(F^i)\ge H^2(P^i_0\cap B(0,1)).\ee

Now since
\be p_0^i(F)=p_0^i\circ \varphi(P_0\cap B(0,1))\supset p_0^i\circ \varphi (P^i_0\cap B(0,1))= F^i\ee
we have
\be H^2(p_0^i(F))\ge H^2(F^i)\ge H^2(P^i_0\cap B(0,1)).\ee
We sum over $i$ and get
\be\begin{array}{ll}
H^2(F)&\ge H^2(p_0^1 F)+H^2(p_0^2 F)\\
&\ge H^2(P^1_0\cap B(0,1))+H^2(P_0^2\cap B(0,1))=H^2(P_0\cap B(0,1)).\end{array}
\ee
This implies that $P_0$ is minimal.$\hfill\Box$

\begin{cor}\label{bigproj}Suppose $\e>0$ is such that $\arccos(\e/2)\le\a_1\le\a_2$, and $P^1,\ P^2$ are two planes with characteristic angles $(\alpha_1,\alpha_2)$. Denote by $p^i$ the orthogonal projection on $P^i,i=1,2$. Then if $E$ is a closed 2- rectifiable set satisfying $p^i(E)\supset B(0,1)\cap P^i$, we have
\be H^2(E)\ge \frac{2\pi}{1+\epsilon}.\ee
\end{cor}

\nd by Proposition 2.19, we know that for almost all point $x\in E$, the approximate tangent plane of $E$ at $x$ $T_xE$ satisfies
\be|p^1(T_xE)|+|p^2(T_xE)|\le 1+\e.\ee
We apply Lemma \ref{projection} to $E$ and obtain the conclusion. $\hfill\Box$

\bigskip

\section{Uniqueness of $P_0$}

In this section, we will add more information on Corollary 2.37. More precisely, we will prove the following Theorem \ref{unicite}. It looks quite natural, but its proof is more complicated than the author imagined.

\begin{thm}\label{unicite} [uniqueness of $P_0$] Let $P_0=P_0^1\cup_\perp P_0^2$, as in the previous section, and denote by $p_0^i$ the orthogonal projection on $P_0^i,i=1,2$. Let $E\subset \overline B(0,1)$ be a 2-dimensional closed reduced set such that $E\cap B(0,1)$ is minimal in $B(0,1)\subset\ \R^4$, and
\be p_0^i(E\cap \overline B(0,1))\supset P_0^i\cap \overline B(0,1);\ee
\be E\cap \partial B(0,1)=P_0\cap \partial B(0,1);\ee
\be \begin{split}&H^2(E\cap B(0,1))\le 2\pi\\
\mbox{ (or equivalently }&H^2(E\cap B(0,1))=2\pi,\mbox{ see near (3.5))}.\end{split}\ee

Then $E=P_0\cap \overline B(0,1)$.

\end{thm}

In the rest of the section, we suppose that $E$ is a set that verifies all the hypotheses in Theorem \ref{unicite}.

First notice that $p_0^i(E)=P_0^i\cap \overline B(0,1)$, because $E\subset \overline B(0,1)$ and $p_0^i(E)\supset P_0^i\cap\overline B(0,1)$. As a result
\be\begin{array}{ll}H^2(E)&=H^2(E\cap B(0,1))=2\pi\\
&=H^2(P_0^1\cap B(0,1))+H^2(P_0^2\cap B(0,1))=H^2(p_0^1(E))+H^2(p_0^2(E)).\end{array}\ee
Compare with (2.36), and observe that here by Lemma 2.14 we can take $\lambda=1$, we have that all the inequalities in (2.36), and hence also in (2.32), are in fact equalities. This means that  
\be \mbox{ for almost all point }x\in E,\ |p_0^1(T_xE)|+|p_0^2(T_xE)|=1,\mbox{ or equivalently }T_xE\in P(\Xi),\ee 
and 
\be \mbox{for }i=1,2,\mbox{ for almost all }z\in p_0^i(E),\ N(p_0^i,z)=\sharp\{{p_0^i}^{-1}(z)\cap E\}=1.\ee

Now we are going to use these two conditions and Lemma 2.17 to get useful local properties of the set $E$. 

First of all, since $E\cap B(0,1)$ is a reduced minimal set in $B(0,1)$, we know that for all $x\in E$, there exists a bi-H\"older ball $B(x,r)\subset B(0,1)$, and in $B(x,r)$ the set $E$ is bi-H\"older equivalent to a minimal cone $C_x$ (c.f. Theorem \ref{holder}). The cone $C_x$ is a blow-up limit of $E$ at $x$. 

On the other hand, all that we know generally for a 2-dimensional minimal cone $C$ is that its intersection with the unit sphere $S=C\cap\partial B(0,1)$ is a finite collection of great circles and arcs of great circles. Every great circle is disjoint from the rest of $S$; at their ends, the arcs meet by set of three, with $120^\circ$ angles, and in particular no free ends exists (c.f.\cite{DJT}, Proposition 14.1). Hence each endpoint of any of these arcs is a one-dimensional $\Y$ point of the net $S$. Thus if a minimal cone is not the union of several transversal planes, its intersection with the unit sphere contains at least a $\Y$ point. As a result, there exist $\Y$ points in $C$ arbitrarily near the origin, since $C$ is a cone.

This discussion yields the following.

\begin{lem}
There is no point of type $\Y$ in $E\cap B(0,1)$.
\end{lem}

\nd If $x\in E$ is a point of type $\Y$, then it means that the (unique, by Theorem \ref{c1} and Remark \ref{ful}) tangent cone $C_x$ is composed of 3 closed half planes $\{P_i\}_{1\le i\le 3}$ which meet along a line $D$ passing through the origin. Denote by $Q_i,1\le i\le 3$ the plane that contains $P_i$. In this case we claim that
\be \mbox{at least one }Q_i\mbox{  doesn't belong to }P(\Xi).\ee

In fact, if we denote by $v$ the unit vector generating $D$, then there exist three unit vectors $w_i,1\le i\le 3$ such that $w_i\perp v$, $Q_i=P(v\wedge w_i)$, and the angle between any two of the $w_i$  is $120^\circ$. 

If $Q_1\not\in P(\Xi)$, then the claim (3.9) holds automatically. Suppose that $Q_1\in P(\Xi)$, and hence $v\wedge w_1\in \Xi$. Then by Lemma 2.17, there exist unit vectors $v_j,u_j,j=1,2$ with $v_j\in P_0^1, u_j\in P_0^2$, $v_1\perp v_2,u_1\perp u_2$, and $\a\in[0,\frac\pi2]$, such that
\be v= \cos\a\ v_1+\sin\a\ u_1,w_1=\cos\a\ v_2+\sin\a\ u_2.\ee

But then since $u_j,v_j,j=1,2$ generate $\R^4$, there exists $a,b,c,d\in\R$ with $a^2+b^2+c^2+d^2=1$ such that $w_2=av_1+bv_2+cu_1+du_2$. Therefore
\be\begin{split}|p_0^1(Q_2)|+|p_0^2(Q_2)|&=|\wedge_2p_0^1(v\wedge w_2)|+|\wedge_2p_0^2(v\wedge w_2)|\\
&=|b\cos\a|+|d\sin\a|\le (b^2+d^2)^\frac12 (\cos^2\a+\sin^2\a)^\frac 12\\
&=(b^2+d^2)^\frac12\le 1.
\end{split}\ee 

So if we want $Q_2$ belong to $P(\Xi)$, all inequalities in (3.11) should be equalities, therefore $b^2+d^2=1$, and $b:d=\pm\cos\a:\sin\a$. As a result, $(b,d)=(\pm\cos\a,\pm\sin\a)$, $a^2+c^2=0$, and hence $w_2=\pm w_1$ or $\pm (\cos\a v_2-\sin\a u_2)$.

Notice that angle between $w_1$ and $w_2$ is $120^\circ$, so there are only two possibilities for $\a$, which are $\frac\pi3$ and $\frac\pi 6$.

If $\a=\frac\pi3$, then $w_1=\frac12 v_2+\frac{\sqrt3}{2}u_2$ and $w_2=\frac 12 v_2-\frac{\sqrt3}{2}u_2$. But the argument above holds also for $w_3$, so we have $w_3=w_2$, which contradicts the fact that the angle between $w_2$ and $w_3$ is $120^\circ$. 

The situation is the same for $\a=\frac\pi6$.

Thus we have proved the claim (3.9).

Now suppose for example that $Q_1\not\in P(\Xi)$. Then since $P(\Xi)$ is closed in $G(4,2)$, there exists an open set $U\subset G(4,2)$ that contains $Q_1$ such that $U\cap P(\Xi)=\emptyset.$

As we said, $\Y$ is a cone with the full-length property. Hence by the theorem of $C^{1+\a}$-regularity for minimal sets (c.f. Theorem \ref{c1}), there exists $r>0$ such that $B(x,r)\subset B(0,1)$, and in $B(x,r)$, $E$ coincides with the image of $C_x+x$ of a $C^1$ homeomorphism $\varphi$ from $\R^4$ to $\R^4$ (recall that $C_x=\cup_{i=1}^3 P_i$). Denote by $R=\varphi(Q_1)$. Then the map $T:R\to G(4,2),T(x)=T_xR$ is continuous. As a result $T^{-1}(U)$ is open. Denote $R_+=\varphi(P_1)$. Then $T^{-1}(U)\cap \varphi^{-1}(R_+)\ne\emptyset$, since it contains $x$. Therefore $T^{-1}(U)$ is relatively open in $R_+$. But $R_+$ is a $C^1$ manifold with boundary, whose boundary $\varphi(D)$ is of zero measure, so $R_+\cap T^{-1}(U)\bs\varphi(D)$ is of positive measure. Now for every $y\in R_+\cap T^{-1}(U)\bs\varphi(D)$, $T_yR=T_yE\not\in P(\Xi)$, thus we have found a set of positive measure and all of its point do not satisfy (3.6); this gives a contradiction.\qed

After this lemma, we will obtain some useful description of the local structure around each point $x$ of $E$.
\begin{lem}
 For each $x\in E\cap B(0,1)$, every blow-up limit of $E$ at $x$ is either a plane belonging to $P(\Xi)$, or $P_0$.
\end{lem}

\nd Suppose that $X$ is a blow-up limit of $E$ at $x$, suppose also that $x=0$ for short. First we claim that
\be X\mbox{ doesn't contain any point of type }\Y.\ee

Suppose this is not true, then there exists $p\in X$ such that $p$ is of type $\Y$. Then $p$ is not the origin, because otherwise $X$ is of type $\Y$, and hence $0$ is of type $\Y$, which gives a contradiction with Lemma 3.8.

So $p$ is not the origin. Then since $X$ is a cone, for every $r>0$, $rp\in X$ is a point of type $\Y$. We can thus suppose that $||p||=1$. Then by our description of 2-dimensional minimal cones (above Lemma 3.8), there exists $0<r<\frac12$ such that in $B(p,r)$, $X$ coincides with a cone $Y$ of type $\Y$ centered at $p$.

Define $d_{x,r}(E,F)=\frac 1r\max\{\sup\{d(y,F):y\in E\cap B(x,r)\},\sup\{d(y,E):y\in F\cap B(x,r)\}\}$, which is the relative distance of two sets $E,F$ with respect to the ball $B(x,r)$. Now $X$ is a blow-up limit of $E$ at $0$, so that there exists $s>0$ (large) such that $d_{0,2}(X,sE)<\frac{r\e_2}{100}$, where $\e_2$ is the constant in Proposition 16.24 of \cite{DJT}. Equivalently, $d_{p,\frac r2}(sE,X)<\frac{\e_2}{50}$.

We want to show that \be d_{p,\frac r2}(sE,X)=d_{p,\frac r2}(sE,Y).\ee

Once we have this, we take a point $z\in sE$ such that $d(z,p)<\frac r2\times \frac{\e_2}{50}$, then $d_{z,\frac r4}(sE,Y+z-p)<\frac{\e_2}{10}$. Here $Y+z-p$ is a $\Y$ cone centered at $z$. But $sE$ is minimal (since $E$ is), therefore Proposition 16.24 of \cite{DJT} gives that $sE$ contains a $\Y$ point, and hence $E$, too. This contradicts Lemma 3.8. Thus we obtain our claim (3.13).

We should still prove (3.14). By definition,
\be d_{p,\frac r2}(sE,X)=\frac 2r\max\{\sup_{x\in sE\cap B(p,\frac r2)}d(x,X),\sup_{x\in X\cap B(p,\frac r2)}d(x,sE)\}.\ee
 
For the second term, we have $X\cap B(p,\frac r2)=Y\cap B(p,\frac r2)$, and hence
 \be\sup_{x\in X\cap B(p,\frac r2)}d(x,sE)=\sup_{x\in Y\cap B(p,\frac r2)}d(x,sE).\ee
 
 For the first term, we have
 \be d(x,X)=d(x,X\cap B(p, r))\mbox{ for all }x\in sE\cap B(p,\frac r2).\ee
  In fact, since $d_{p,\frac r2}(sE,X)<\frac{\e_2}{50}$, for each $x\in sE\cap B(p,\frac r2)$, $d(x,X)<\frac{\e_2}{50}\times\frac 2r$, so $d(x,X)=d(x,X\cap B(x,\frac{\e_2}{50}\times\frac 2r))\le d(x, X\cap B(p, r))$, because $B(x,\frac{\e_2}{50}\times\frac 2r)\subset B(p, r)$. On the other hand, $X\cap B(p,r)\subset X$, therefore $d(x,X)\ge d(x,X\cap B(p,r))$. Thus we have (3.17), and hence 
  \be \begin{split}\sup_{x\in sE\cap B(p,\frac r2)}d(x,X)&=\sup_{x\in sE\cap B(p,\frac r2)}d(x,X\cap B(p,r))\\
  &=\sup_{x\in sE\cap B(p,\frac r2)}d(x,Y\cap B(p,r))\ge \sup_{x\in sE\cap B(p,\frac r2)}d(x,Y).
  \end{split}\ee
  
Combining with (3.16) we have 
  \be d_{p,\frac r2}(sE,X)\le d_{p,\frac r2}(sE, Y).\ee
  A similar argument yields
  \be d_{p,\frac r2}(sE,Y)\le d_{p,\frac r2}(sE, X).\ee
  
So we have (3.14), and (3.13) follows.

Since $X$ is a minimal cone, as we have said before, $X\cap\partial B(0,1)$ is a finite collection of great circles and arcs of great circles that meet by 3 with angles of $120^\circ$. Then (3.13) implies that there is no such arcs, since $X$ does not have $\Y$ points. As a result, $X\cap\partial B(0,1)$ is a finite collection of great circles, and therefore $X$ is the union of a finite number of transversal planes.

By Remark \ref{ful}, $X$ is a full-length cone, so the $C^1$ regularity holds (Thm \ref{c1}). Then by the same argument as for $\Y$ in the proof of Lemma 3.8, every plane in $X$ belongs to $P(\Xi)$, since $P(\Xi)$ is closed.

As a result, if $X$ is not a plane, then $X=\cup_{i=1}^n Q_n$ with $Q_n\in \Xi$ with $n\ge 2$. These $Q_n$ are transversal, by Lemma 2.17.  Moreover in a small neighborhood $B(x,r)$ of the point $x$, the set $E$ is a union of $n$ transversal $C^1$ manifolds $S_i,1\le i\le n$, and the tangent plane to $S_i$ at $x$ is $Q_n$. 

By Lemma 2.17, for all $Q\in\Xi$ and $Q\ne P_0^2$, there exists $s=s_Q>0$ such that $p_0^1(Q\cap B(0,1))\supset P_0^1\cap B(0,s_Q)$. Hence if $X$ contains two planes $Q_i,Q_j, i\ne j$, that are different from $P_0^2$, there exists $s>0$ such that $p_0^1(Q_1\cap B(0,1))\cap p_0^1(Q_2\cap B(0,1))\supset B(0,s)\cap P_0^1$. Then since the tangent plane to $S_i$ at $x$ is $Q_i$ and that to $S_j$ at $x$ is $Q_j$,  there exists a neighborhood $U$ of $p_0^1(x)$ which is contained in both projections $p_0^1(S_i)$ and $p_0^1(S_j)$. Notice that $S_i$ and $S_j$ are transverse, hence $S_i\cap S_j=\{x\}$. This yields that each point $y\in U\bs p_0^1(x)$ admits at least two pre-images by $p_0^1$ in $E\cap B(x,r)$, one is in $S_i$ and the other is in $S_j$, i.e.,
\be (U\bs\{p_0^1(x)\})\cap P_0^1\subset\{z\in P_0^1,\sharp\{{p_0^1}^{-1}\{z\}\cap E\}\ge 2\}.\ee

But $U$ is open, so $(U\bs\{p_0^1(x)\})\cap P_0^1$ is of positive measure, which contradicts (3.7).

Hence in $X$, we have at most one plane different from $P_0^2$. The same argument gives also that we have at most one plane which is not $P_0^1$. But $X$ contains at least two planes, therefore $X=P_0^1\cup P_0^2=P_0$. \qed

Lemma 3.12 says that there exists only two types of blow-up limits in $E$, and both of them have the full-length property. As a result, by Thm \ref{c1}, around each point $x\in E$, $E$ is locally $C^1$ equivalent to a plane or to $P_0$.

The two next lemmas will describe more precisely what happens around each of the two types of singularities. We identify $\R^4$ with $\C\times\C$.

\begin{lem}
Let $\{e_1,e_2=ie_1,e_3,e_4=ie_3\}$ be an orthonormal basis of $\R^4$, with $P_0^1=e_1\wedge e_2$ and $P_0^2=e_3\wedge e_4$. Let $x\in E\cap B(0,1)$ be a point of type $\mathbb P$ such that $T_xE\ne P_0^2$. Then there exists $r=r(x)>0$ such that in $B(x,r),\ E$ is (under the given basis)  the graph of a complex analytic or anti-analytic function 
$\varphi=\varphi_x : P^1_0\to P_0^2 $. More precisely,
\be E\cap B(x,r)= \mbox{graph}(\varphi)\cap B(x,r).\ee
We have also a similar conclusion for $x$ such that $T_xE\ne P_0^1$, i.e., near $x$, $E$ is the graph of a complex analytic or anti-analytic function from $P_0^2$ to $P_0^1$. 
\end{lem}

\nd We will only prove it for $T_xE\ne P_0^2$. The other case is similar.

So let $x\in E$ be as in the lemma. Assume all the hypotheses in the lemma. Since $x$ is a $\P$ point, there exists $r_1>0$ such that in $B(x,r_1)$, $E$ is the graph of a $C^1$ function $\varphi_1:T_xE\to T_xE^\perp$ (c.f. Thm \ref{c1}). If we denote by $\pi$ the projection from $\R^4$ to $T_xE$,  and define $F:\R^4\to \R^2$ by $F(y)=y-\varphi_1(\pi(y))\in T_xE^\perp\cong\R^2$ for $y\in B(x,r_1)$, then $F$ is $C^1$ regular, and $E\cap B(x,r_1)=\{y\in B(x,r_1):F(y)=(0,0)\}$. As a result, $DF(x)$ is a linear map from $\R^4$ to $\R^2$, with $T_xE=Ker DF(x)$.

By the hypothesis, $T_xE\ne P_0^2$, hence by Lemmas 3.12 and 2.17, $T_xE$ can be represented by $T_xE=P((\cos\a v_1+\sin\a u_1)\wedge (\cos\a v_2+\sin\a u_2))$, with $\a\in [0,\frac\pi2)$, $v_i,u_i,i=1,2$ unit vectors, $v_i\in P_0^1,u_i\in P_0^2$ and $v_1\perp v_2,u_1\perp u_2$. In particular, $P_0^2\cap Ker DF(x)=P_0^2\cap T_xE=\{0\}$. This means that $DF(x)|_{P_0^2}$ is invertible. Then by the implicit function theorem, there exists $r_2>0$ and a $C^1$ map $\varphi:P_0^1\to P_0^2$, such that for every $y=(y_1,y_2)\in P_0^1\times P_0^2\cap B(x,r_2)$, $F(y)=0\Leftrightarrow y_2=\varphi(y_1)$. As a result, in $B(x,r_2)$, $E$ is the graph of $\varphi: P_0^1\to P_0^2$.

We still have to prove that $\varphi$ is complex analytic or anti-analytic. This will follow from a particular property of $P(\Xi)$.

More precisely, return to the orthonormal basis $\{e_1,e_2=ie_1,e_3,e_4=ie_3\}$ with $P_0^1=e_1\wedge e_2$ and $P_0^2=e_3\wedge e_4$. Use Lemma 2.17, and expand each unit vector under the basis; a simple calculation gives
\be\begin{split}\Xi=\{\pm& [ae_1+be_2+ce_3+de_4]\wedge[-be_1+ae_2\pm(-de_3+ce_4)],\\
&a,b,c,d\in\R,a^2+b^2+c^2+d^2=1\},\end{split}\ee
and hence
\be \begin{split}P(\Xi)=\{P&([ae_1+be_2+ce_3+de_4]\wedge[-be_1+ae_2\pm(-de_3+ce_4)]),\\
&a,b,c,d\in\R,a^2+b^2+c^2+d^2=1\}.\end{split}\ee

If $T_yE=[ae_1+be_2+ce_3+de_4]\wedge[-be_1+ae_2+(-de_3+ce_4)]$, this means that $d\varphi(y)(a+bi)=c+di,d\varphi(y)(-b+ai)=-d+ci$. But $d\varphi(y)$ is (real) linear, and $a^2+b^2\ne 0$, therefore $a+bi,-b+ai$ is a basis, and $d\varphi(y)$ is determined by its values at these two points. Notice that $\frac{c+di}{a+bi}=\frac{-d+ci}{-b+ai}$, and denote it by $A\in\C$, then we have $d\varphi(y)(z)=Az$, which is complex analytic. In other words,
\be\frac{d\varphi}{d\bar z}(y)=0.\ee

If $T_yE=[ae_1+be_2+ce_3+de_4]\wedge[-be_1+ae_2-(-de_3+ce_4)]$, a similar argument gives $d\varphi(y)(z)=B\bar z,$ and hence
\be\frac{d\varphi}{d z}(y)=0.\ee

So we have proved that $\varphi$ is a complex $C^1$ function such that at each point, it is either analytic, or anti-analytic.

Denote by $B=B(p^1(x),r)\cap P_0^1$, $B_1=\{y\in B, \frac{d\varphi}{d z}(y)\ne 0\}$. Then $B_1$ is open, since $\varphi$ is $C^1$. If $B_1=\emptyset$, then $\varphi$ is anti-analytic. Otherwise, $B_1$ is not empty, and set $g=\frac{\partial \varphi}{\partial z}$. Then $g$ is continuous on $B$, and $B_1=\{y\in B: g(y)\ne 0\}$. Moreover, since $\frac{d\varphi}{d z}(y)\ne 0$ on $B_1$, $\frac{d\varphi}{d \bar z}(y)=0$ on $B_1$ (since at each point $\varphi$ is either analytic or anti-analytic), and hence $\varphi$ is holomorphic on the open set $B_1$, so that its derivative $g$ is holomorphic, too. Then conclusion of Lemma 3.22 will follow from the following theorem (c.f.\cite{Ru} Thm 12.14) :

\begin{thm}[Rad\'o's theorem] Let $U\subset \C$ be an open domain, and $f$ be a continuous function on $\overline U$. Set $\Omega=\{z\in U:f(z)\ne 0\}$, and suppose that $f$ is holomorphic on $\Omega$. Then $f$ is holomorphic on $U$.
 \end{thm} 
 
In fact, we apply the theorem to $g$, and obtain that $g$ is holomorphic on $B$. But since $B_1\ne\emptyset$, $g\not\equiv 0$. As a result $B_1^C=\{y\in B: g(y)=0\}$ doesn't have any limit point. Notice that $B_1^C\supset B_2:=\{y\in B,\frac{d\varphi}{d\bar z}(y)\ne 0 \}$, and $B_2$ is open, so it is an open set without limit point. Therefore $B_2=\emptyset$, which means that $\varphi$ is complex analytic on $B$.

Hence $\varphi$ is complex analytic or anti-analytic on $B=B(p^1(x),r)\cap P_0^1$. \qed

\begin{lem}
If $x\in E\cap B(0,1)$ is of type $P_0$, then there exists $r=r(x)>0$ such that
\be E\cap B(x,r)=(P_0+x)\cap B(x,r).\ee 
\end{lem}

\nd Let $B(x,r')$ be a $C^1$ ball for the point $x$, in which $E$ coincides with the image of $P_0+x$ by a $C^1$ function $\varphi$ : $E\cap B(x,r')=\varphi ((P_0+x)\cap B(x,r'))$, where $\varphi$ is a diffeomorphism and $\varphi(x)=x,$ $D\varphi(x)=Id$ (c.f. Thm \ref{c1}). Set $A_i=\varphi(P_0^i+x)\cap B(x,r'),\ i=1,2$.  Then $A_i$ are transversal $C^1$ manifolds, $[A_1\cup A_2]\cap B(x,r')=E\cap B(x,r')$ and $A_1\cap A_2=\{x\}$. Moreover, the tangent plane to $A_i$ at the point $x$ is $P_0^i$, $i=1,2$. In particular, $P_0^1(A_1)$ contains a neighborhood $U$ of $p_0^1(x)$ in $P_0^1$.

By a similar argument as the one in Lemma 3.22, there exists $r_2<\frac{r_1}{2}$ such that in $B(x,r_2)$, $A_2$ is the graph of the analytic or anti-analytic complex function $\psi\ :\ P_0^2+x\to P_0^1$, where $\psi=p_0^1\circ\varphi$ on $P_0^2+x\cap B(x,r_2)$). Without loss of generality, suppose that it is complex analytic.

If $\psi$ is not constant, it will be an open map, because it is analytic. Therefore $p_0^1(A_2)$ also contains a neighborhood $U'$ of $p_0^1(x)$ in $P_0^1$. But $A_1\cap A_2=\{x\}$, hence for every point $y\in U\cap U'\bs \{p_0^1(x)\}$, ${p_0^1}^{-1}\{y\}$ admits at least two points in $E$, one in $A_1$, and the other in $A_2$. But $U\cap U'\bs\{p_0^1(x)\}$ is of non zero measure, this contradicts (3.7).
 
 Hence $\psi$ is constant. As a result, 
 \be\varphi((P_0^2+x)\cap B(x,r_2))=P_0^2\cap B(x,r_2).\ee
 
 We can do the same for $\varphi(P_0^1+x)$, and we obtain that there exists $r<r_2$ such that in $B(x,r)$, $\varphi(P_0^1+x)$ is $P_0^1+x$ itself. This completes the proof of Lemma 3.29. \qed
 
 Let us sum up a little. The minimal set $E$ has only two types of points : planar points, around which the set $E$ is a $C^1$ (locally analytic or anti analytic) manifold; and $P_0$ points, around which $E$ is locally $P_0$ itself. Hence in fact, 
 \be\begin{split}\mbox{The set }&E\mbox{ is composed of countably many transversal }C^1\mbox{ manifolds }S_1,S_2,\cdots, S_n,\cdots.\\ 
 \mbox{They}&\mbox{ are locally analytic or anti analytic, and at any of their intersections, }\\
 &E\mbox{ is locally a translation of }P_0.\end{split}\ee

\medskip

We know that (3.7) holds for almost all $z\in p_0^1(E)$. But now let us look at those points $x$ in $E$ where the projection $p_0^1$ is not injective.

\begin{lem}
Suppose that $y_1,y_2\in E\cap B(0,1)$ are such that $p_0^1(y_1)=p_0^1(y_2)$. Then at least one of the two $y_i$ is such that $E$ coincides with $P_0^2+y_i$ in a neighborhood of $y_i$.  
\end{lem}

\nd By the argument near (3.21), at least one of the tangent planes $T_{y_1}E,\ T_{y_2}E$ is $P_0^2$. Suppose $T_{y_1}E=P_0^2$ for example. Then by Lemma 3.22, there exists $s_1>0$ such that in $B(y_1,s_1)$, $E$ is the graph of a complex analytic or anti-analytic function (under the basis given before) $\varphi_1:P_0^2\to P_0^1$, with $\varphi_1(0)=p_0^1(y_1),D\varphi_1=0$.

First, if the tangent plane $T_{y_2}E$ of $E$ at $y_2$ is not $P_0^2$, then again by the argument near (3.21), there exists $s_2>0$ such that $B(y_1,s_1)\cap B(y_2,s_2)=\emptyset$, and $B(p_0^1(y_2),s_1)\cap P_0^1\subset p_0^1(B(y_2,s_2)\cap E)$. Therefore $p_0^1(B(y_1,s_1)\cap E)$ is of measure zero, because of (3.7) and $p_0^1(y_1)=p_0^1(y_2)$. Hence $\varphi_1$ has to be constant, otherwise it will be an open map, which means $p_0^1(B(y_1,s_1)\cap E)$ contains an open set in $P_0^1$, which is of positive measure, and thus leads to a contradiction. So in $B(y_1,s_1),\ E$ coincides with $P_0^2+y_1$. (See the picture 3-1 below, the two points on the left.)

\noindent\centerline{\includegraphics[width=0.5\textwidth]{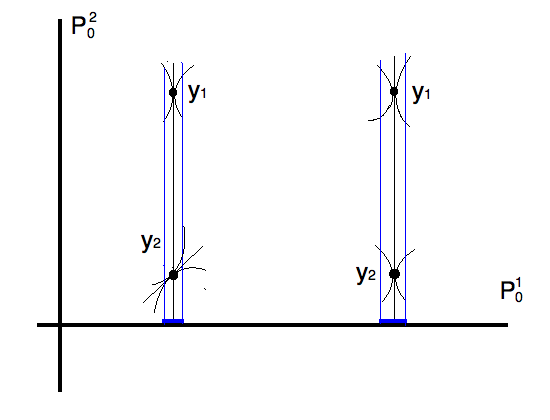}}
\nopagebreak[4]
\centerline{3-1}

If $T_{y_2}E$ is also $P_0^2$, then in $B(y_2,s_2)$, $E$ is the graph of a complex analytic or anti-analytic function $\varphi_2:P_0^2\to P_0^1$. Then at least one of the two $\varphi_i,i=1,2$ is constant, otherwise all of them are open maps, so that $p_0^1(B(y_1,s_1)\cap E))\cap p_0^1(B(y_2,s_2)\cap E))$ contains an open set in $P_0^1$, which contradicts (3.7).\qed

\begin{lem} If there exists a point in $P_0^1\cap B(0,1)$ whose pre-image by $p_0^1$ contains more than one point in $E\cap B(0,1)$, then $P_0^2\cap B(0,1)\subset E$. 
\end{lem}

\nd

Lemma 3.33 says that if there exists a point in $P_0^1\cap B(0,1)$ whose pre-image by $p_0^1$ contains more than one point in $E$, then there exists a piece of $P_0^2$ above it. Denote this piece by $(P_0^2+x)\cap B(x,r)$. On the other hand,  by (3.32), $E\cap B(0,1)$ is composed of at most countably many transversal $C^1$ manifolds $S_1,\cdots,S_l,\cdots$, so there exists $i$ such that $(P_0^2+x)\cap B(x,r)\subset S_i$.

We claim that $S_i=(P_0^2+x)\cap B(0,1)$. So set $A=\{y\in S_i:\mbox{ there exists }r=r_y>0\mbox{ such that }S_i\cap B(y,r)=(P_0^2+y)\cap B(y,r)\}$. By definition, $A$ is open in $S_i$. On the other hand, for any $y\in\bar A\cap S_i$, there exists a sequence $y_n\in A$ that converges to $y$. But $S_i$ is a $C^1$ manifold, hence the tangent plane $T_yS_i=\lim_{n\to\infty}T_{y_n}S_i=P_0^2$. Therefore there exists $r>0$ such that in $B(y,r)$, $S_i$ is the graph of a complex analytic (the anti-analytic case is exactly the same) $\psi:P_0^2\to P_0^1$. But $\psi'=0$ on $\{p_0^2(y_n)\}\cup\{p_0^2(y)\}$, which, if $y\not\in A$, is a infinite set with a limit point. So $\psi'\equiv 0$ in $p_0^2(B(y,r))$, so that $\psi$ is constant, which yields that $y\in A$. So $A$ is closed. But $S_i$ is connected, $A$ is open and closed and non-empty, hence $A=S_i$. But $S_i$ is a $C^1$ manifold, so the only possibility is that it is a piece of $(P_0^2+x)\cap B(0,1)$. However, $S_i$ is both closed and open in $(P_0^2+x)\cap B(0,1)$, because $A$ is both open and closed in $(P_0^2+x)\cap B(0,1)$ by the same argument as above, and $S_i=A$. As a result, $S_i=(P_0^2+x)\cap B(0,1)$.

By (3.3), $\overline S_i\cap \partial B(0,1)\subset E\cap \partial B(0,1)\subset P_0\cap \partial B(0,1)$, since $E$ is closed. This implies that $S_i=P_0^2\cap B(0,1)$. Hence $P_0^2\cap B(0,1)\subset E$.\qed

\begin{rem}By a similar discussion, if $p_0^2$ is not injective on $E\cap B(0,1)$, then $P_0^1\cap B(0,1)\subset E$. \end{rem}

\begin{lem} The set $E$ contains at least one point of type $P_0$.\end{lem}

\nd First of all, by Lemma 3.34 and Remark 3.35, if neither of the $p_0^i,i=1,2$ is injective, then $P_0\cap B(0,1)=(P_0^1\cup P_0^2)\cap B(0,1)\subset E$. In particular, there exists a $P_0$ point.

Thus if $E$ contains no $P_0$ point, then at least one of the $p_0^i,i=1,2$ is injective. Without loss of generality, we suppose $p_0^1$ is injective on $E\cap B(0,1)$, and $E$ contains no $P_0$ point. Notice that (3.2) is true for $E\cap\overline B(0,1)$, and we know that $p_0^1(E\cap\partial B(0,1))=\{0\}\cup (P_0^1\cap\partial B(0,1))$, hence $p_0^1(E\cap B(0,1))\supset P_0^1\cap B(0,1)\bs\{0\}$. Hence the injectivity gives
\be\mbox{ for all }z\in P_0^1\cap B(0,1)\bs\{0\}, \sharp\{{p_0^1}^{-1}(z)\cap E\cap B(0,1)\}=1,\ee
and 
\be \sharp\{{p_0^1}^{-1}(0)\cap E\cap B(0,1)\}\le1.\ee

So there exits $\psi:P_0^1\cap B(0,1)\bs\{0\}\to P_0^2$ such that $E\cap B(0,1)$ is its graph. Therefore above those points $y\in E$ whose tangent planes $T_yE\ne P_0^2$, $\psi$ is locally a complex analytic or anti-analytic function (by Lemma 3.22). 

We claim that there is no singular point $y\in E\cap B(0,1)$ with $T_yE=P_0^2$. In fact, if $y$ is a such point, then by Lemma 3.22 (recall that we have supposed that $E$ contains no $P_0$ point), in a neighborhood $B(y,r)$, $E$ is the graph of a complex analytic or anti-analytic function $g$ from $P_0^2$ to $P_0^1$. Suppose $p_0^2(y)=0$ for simplicity, hence $g'(0)=0$. Then $0$ is a zero of $g-g(0)$ of order at least 2. Thus in a  punctured neighborhood $O$ of $g(0)\in\C$, each point has at least 2 pre-images. So $O\subset \{z\in P_0^1:\sharp\{{p_0^1}^{-1}\{z\}\cap E\}\ge 2\}$.  But $O$ is of positive measure because it is open. This contradicts (3.7).

So there is no singular point for $\psi$, therefore $\psi$ is a $C^1$ function on $P_0^1\cap B(0,1)\bs\{0\}$, and $E\cap B(0,1)$ is its graph. For each point $x\in P_0^1\cap B(0,1)\bs\{0\}$, $\psi$ is analytic or anti analytic in a small neighborhood of $x$. Moreover $\psi$ is bounded.
By the same argument in Lemma 3.22, we know that $\psi$ is globally analytic or anti-analytic on any compact subdomain of $P_0^1\cap B(0,1)\bs\{0\}$, hence is globally analytic or anti-analytic on $P_0^1\cap B(0,1)\bs\{0\}$. Suppose it is analytic. Then since $\psi$ is bounded around $\{0\}$, we can extend $\psi $ to $P_0^1\cap B(0,1)$. But since $\psi$ is analytic with $\psi|_{\partial B(0,1)}=0$, the graph of $\psi$ has to be $P_0^1\cap B(0,1)$, hence 
\be E\cap B(0,1)\subset graph (\psi)=P_0^1\cap B(0,1).\ee
This contradicts the hypothesis (3.3).

Thus we obtain the existence of a point of type $P_0$ in $E$.\qed
%
%\begin{lem}Let $x\in E$ be a point of type $P_0$, then $E$ contains $(P_0+x)\cap B(0,1)$.
%\end{lem}
%
%\nd We claim that $(P_0^1+x)\cap E$ is relatively open in $(P_0^1+x)\cap B(0,1)$. In fact, for all $y\in (P_0^1+x)\cap E\backslash\{x\}$, $p_0^2(y)=p_0^2(x)$. By Remark 3.35,  hence at least one of $x$ and $y$ is such that $E$ is a disc parallel to $P_0^1$ in a neighborhood. But $x$ is already a point of type $P_0$, so there exists a ball $B_y$ centered at $y$ such that $E\cap B_y=(P_0^1+x)\cap B_y$. As for $x$, by Lemma 3.29, there exists a ball $B_x$ such that $B_x\cap E=B_x\cap (P_0+x)\supset B_x\cap (P_0^1+x)$. So $E$ is relatively open in $P_0^1+x$.
%
%On the other hand, $E$ is itself a closed set in $B(0,1)$, hence is relatively closed in $P_0^1+x$. Therefore the only choice is
%\be E\cap (P_0^1+x)\cap B(0,1)=(P_0^1+x)\cap B(0,1),\ee
%since $E\cap (P_0^1+x)$ is not empty.
%
%The same argument gives also
%\be E\cap (P_0^2+x)\cap B(0,1)=(P_0^2+x)\cap B(0,1).\ee
%Thus we get the conclusion. \qed

\noindent \textbf{Proof of Theorem \ref{unicite}.} We know by Lemma 3.36 that $E$ contains a point $x$ of type $P_0$. By Lemma 3.29, there exists $r>0$ such that $E\cap B(x,r)=(P_0+x)\cap B(x,r)$. But in this case, neither $p_0^1$ nor $p_0^2$ is injective on $E\cap B(0,1)$. Hence by Lemma 3.34 and Remark 3.35, $P_0\cap B(0,1)=(P_0^1\cup P_0^2)\cap B(0,1)\subset E$. But $H^2(E)=H^2(P_0\cap B(0,1))$, and $E\cap B(0,1)$ is reduced. As a result, $E\cap B(0,1)=P_0\cap B(0,1)$. This completes the proof of Theorem \ref{unicite}.\qed

\section{Existence of minimal sets}

Now we begin to prove Theorem \ref{main}.

Suppose that the conclusion of Theorem \ref{main} is not true. Then there exists a sequence of unions of 2 planes $P_k=P^1_k\cup_{\theta_k}P^2_k\subset\R^4$ which are not minimal, with $\theta_k\ge\frac\pi2-\frac 1k$. Moreover, we can suppose that all the $P_k^1$ are equal to $P_0^1$.

Our proof would be simpler if we could find a sequence of minimizers $E_k$ such that for each $k$, $E_k$ minimizes the measure among all deformations of $P_k\cap\overline B(0,1)$ in $B(0,1)$. (We shall call such a set a solution of Plateau's problem). It is a standard fact about the Hausdorff distance that we can find a subsequence (which we shall still denote by $\{E_k\}$ to save notation) such that $\{E_k\}$ converges to a limit $E_\infty$. Since $P_k\cap\partial B(0,1)=E_k\cap\partial B(0,1)$, and $P_k\cap \partial B(0,1)\to P_0\cap \partial B(0,1)$, we can hope to prove that the boundary of $E_\infty$ is $P_0\cap \partial B(0,1)$. Moreover, since each $E_k$ is minimal in $B(0,1)$, their limit $E_\infty$ has to be minimal in $B(0,1)$, too (c.f.\cite {GD03} Thm 4.1). Notice also that (3.2) is true since $E_\infty$ is the limit of a sequence of deformations, we would deduce that $E_\infty=P_0\cap \overline B(0,1)$, by the uniqueness theorem \ref{unicite}. We could begin our stopping time argument as soon as we got such a sequence.

Unfortunately we do not know any theorem which guarantees the existence of a solution of Plateau's problem. So we have to work more to get through this. We are disposed of a weaker existence theorem:

\begin{thm}[Existence of minimal sets\label{vincent}; c.f. \cite{Fv}, Thm 6.1.7]
Let $U\subset\R^n$ be an open domain, $0<d<n$, and let $\mathfrak F$ be a class of non-empty sets relatively closed in $U$ and satisfying  (1.17), which is stable by deformations in $U$. Suppose that 
 \be \inf_{F\in\mathfrak F} H^d(F)<\infty.\ee
 Then there exists $M>0$ (depends only on $d$ and $n$), a sequence $(F_k)$ of elements of $\mathfrak F$, and a set $E$ of dimension $d$ relatively closed in $U$ that verifies (1.17), such that:
 
 (1) There exists a sequence of compact sets $\{K_m\}_{m\in\N}$ in $U$ with $K_m\subset K_{m+1}$ for all $m$ and $\cup_{m\in N}K_m=U$, such that 
 \be\lim_{k\to\infty}d_H(F_k\cap K_m,E\cap K_m)=0\mbox{ for all }m\in \N;\ee
 
 (2) For all open sets $V$ such that $\overline V$ is relatively compact in $U$, from a certain rank,
 \be F_k\mbox{ is }(M,+\infty)\mbox{-quasiminimal in }V;\ee
 (see Definition 4.5 below) 
 
 (3) $H^d(E)\le\inf_{F\in\mathfrak F}H^d(F)$ ;
 
 (4) $E$ is minimal in $U$.
\end{thm}

\begin{defn}[Quasiminimal set]Let $0<d<n$ be two integers, $M>0,\delta>0$, let $U$ be open in $\R^n$. The set $E$ is said to be $(M,\delta)-$quasi minimal in $U$ ($E\subset QM(U,M,\delta)$ for short)  if $E$ is closed in $U$, (1.17) is true, and for every deformation $\varphi$ in $U$ as in Definition 1.11 such that diam$(\widehat W)<\d$, 
\be H^d(E\cap W_1)\le MH^d(\varphi_1(E\cap W_1)).\ee
\end{defn}

\begin{rem}In Theorem \ref{vincent}, we can also ask that $\{F_k\}$ is a minimizing sequence, i.e., $\lim_{k\to\infty} H^2(F_k)=\inf_{F\in\mathfrak F}H^d(F).$
\end{rem}

Theorem \ref{vincent} is a weaker result, which gives a certain type of minimizer, without preserving the boundary condition, and this minimizer may not be in the class $\mathfrak F$, either. But for our case, we do not need all these strong properties. Some weaker properties are sufficient for us to carry on the proof (as we will see soon).

Recall that $P_k=P_k^1\cup_{\theta_k}P_k^2$ is a sequence of unions of two planes which are not minimal, with $\theta_k>\frac\pi2-\frac 1k.$ Moreover we suppose that all the $P_k^1$'s are the same. Thus $P_k^1=P_0^1$. Denote by $P_0^2$ the plane orthogonal to $P_0^1$, $P_0=P_0^1\cup_\perp P_0^2$. Then $P_k\cap \overline B(0,1)$ converges to $P_0\cap \overline B(0,1)$ for the Hausdorff distance.

\begin{pro}For each $k$, there exists a closed set $E_k\subset \overline B(0,1)$ such that

(1) $ E_k\mbox{ is minimal in }\R^4\bs [P_k\backslash B(0,1)];$

(2) $\partial B(0,1)\cap E_k=\partial B(0,1)\cap P_k;$

(3) $p_k^i (E_k)\supset P_k^i\cap \overline B(0,1)$, where $p_k^i$ denotes the projection on $P_k^i,i=1,2$;

(4) $ H^2(E_k)<H^2(P_k\cap B(0,1))=2\pi;$

(5) $E_k$ is contained in the convex hull of $P_k\cap \overline B(0,1)$.
\end{pro}

\noindent\textbf{Proof of (1) and (4) of Proposition 4.8.}
Fix a $k$. Take $U=\R^4\bs[P_k\bs B(0,1)]$, and let $\mathfrak F$ be the class of all deformations of $P_k\cap\overline B(0,1)$ in $U$. Then by Theorem 4.1, and Remark 4.7, for $d=2$, there exists a minimizing sequence of sets $F_l\in \mathfrak F$ that are also uniformly quasiminimal in $U$, with a uniform constant $M$. Moreover the sequence converges under the Hausdorff distance. Denote by $E_k$ its limit. Then by the conclusion of Theorem 4.1, the terms (1) and (4) of Proposition 4.8 are automatically true. In (4) we have a strict inequality because we have supposed that $P_k$ is not minimal.\qed

For (3), we begin by proving the following lemma.

\begin{lem}Let $P$ be a plane in $\R^4$, and $p$ the projection on $P$. let $\varphi$ be a Lipschitz mapping from $\R^4$ to $\R^4$, such that $\varphi|_{P\cap B(0,1)^C}=id$. Then
\be p[\varphi(P\cap \overline B(0,1))]\supset P\cap\overline B(0,1).\ee
\end{lem}

\nd We prove it by contradiction.

Suppose that there exists  $x\in P\cap \overline B(0,1)$ such that $x\not\in p[\varphi(P\cap\overline B(0,1))]$. Then $x\in P\cap B(0,1)$, since $P\cap\partial B(0,1)=\varphi (P\cap\partial B(0,1))\subset p[\varphi(P\cap\overline B(0,1))]$ by hypothesis. Then there exists $r>0$ such that $B(x,r)\subset B(0,1)$. On the other hand, $P\cap\overline B(0,1)$ is compact, hence its image $p[\varphi(P\cap\overline B(0,1))]$ by the continuous map $p\circ\varphi$ is also compact, such that $\{p[\varphi(P\cap\overline B(0,1))]\}^C$ is open. As a result, there exists $r'<r$ such that $B(x,r')\cap p[\varphi(P\cap\overline B(0,1))]=\emptyset$. In other words, $\varphi(P\cap\overline B(0,1))\subset \R^4\bs p^{-1}[B(x,r')\cap P].$

Define $g:\R^4\bs p^{-1}[B(x,r')\cap P]\to p^{-1}[\partial B(0,1)\cap P]$ as follows. For $x\in \R^4\bs p^{-1}[B(0,1)\cap P]$, let $g(x)$ be the shortest distance projection of $\R^4$ onto $p^{-1}[B(0,1)\cap P]$, and for $y\in [p^{-1}[(B(0,1)\bs B(x,r'))\cap P]$, let $g(y)$ be the intersection of $p^{-1}[\partial B(0,1)\cap P]$ with the ray $[x,y)$ of end point $x$ and passing through $y$. Then $g$ is continuous.

Notice that $\varphi(P\cap\overline B(0,1))\subset \R^4\bs p^{-1}[B(x,r')\cap P]$, hence $p\circ g\circ\varphi$ is Lipschitz, and sends $P\cap\overline B(0,1)$ continuously to $P\cap\partial B(0,1)$, with all points of $P\cap\partial B(0,1)$ fixed. This is impossible. \qed

\noindent\textbf{Proof of Proposition 4.8 (3).} We know that $E_k$ is the limit of a sequence of deformations $F_l=\varphi_l(P_k)$ of $P_k$. Then for each $l$, we have, by Lemma 4.9,
\be p_k^i(F_l)\supset p_k^i[\varphi_l(P_k^i\cap \overline B(0,1))]\supset P_k^i\cap \overline B(0,1),i=1,2.\ee 
As a result,
\be p_k^i(E_k)\supset P_k^i\cap \overline B(0,1),i=1,2.\ee \qed

We still have to prove (2) and (5) of Proposition 4.8. In fact all we have to do is to prove (5), because it says that $E_k$ is contained in the convex hull $C$ of $P_k\cap\overline B(0,1)$, which gives (2).

Let us verify that $E_k\subset \overline B(0,1)$. First of all we claim that for each $\e>0$, there exists $N=N(\e)>0$ such that for all $l>N$,
\be F_l\subset  B(B(0,1)\cup P_k,\e),\ee
That is to say, $F_l$ is contained in an $\e-$neighborhood of the union of the unit ball and the boundary of $U$.

In fact, those $F_l$ are uniformly quasiminimal in $U$, hence are locally uniformly Ahlfors regular in $U$ (c.f. \cite{DS00} Proposition 4.1). This means that there exists a constant $C>1$ such that for all $l$, for all $x\in F_l$ and all $r>0$ with $B(x,2r)\subset U$,
\be C^{-1}r^2\le H^2(F_l\cap B(x,r))\le Cr^2.\ee

Then if there exists $x\in F_l$ such that $d(x, \overline B(0,1)\cup P_k)>\e$, we should have $B(x,\e)\subset U$, and by the Ahlfors regularity,
\be H^2(F_l\cap B(x,\frac 12\e))\ge (4C)^{-1}\e^2.\ee

Now we deform $F_l$ into $\overline B(0,1)$ by the radial projection $\pi$ on $\overline B(0,1)$. Here for each $l$, $F_l$ is a deformation of $P_k\cap\overline B(0,1)$ in $U$, it means that there exists a compact set $K\subset U$, such that $F_l\bs K=P_k\cap\overline B(0,1)\bs K$. The compactness of $K$ gives that $0<d=d(K,\partial U)=d(K,P_k\bs B(0,1))$, and there exists $R>0$ such that $K\subset B(0,R)$. Set $V=B(0,R)\bs\overline B(P_k\bs B(0,1),d)$. Then in $V$, $\pi$ is homotopic to the identity map. Now we set $\pi':U\to \overline B(0,1)\cap U,\pi'(x)=t(x)x+(1-t(x))\pi(x)$, where $t(x)=\min\{2d(x,\overline V)/d,1\}$. Then $\pi'$ is a deformation on $U$, and $\pi'|_V=\pi|_V$, which gives $\pi'(F_l)=\pi(F_l)$. Hence the set $\pi(F_l)$ is a deformation of $P_k$ in $U$. In other words $\pi(F_l)\in\mathfrak F.$

\begin{lem} Let $E$ be rectifiable such that $E\cap B(0,1+\e)=\emptyset$. Then
\be H^2(\pi(E))\le \frac{1}{(1+\e)^2} H^2(E).\ee
\end{lem}

\nd We are going to prove that $\pi$ is $\frac{1}{1+\e}$-Lipschitz on $E$, which gives (4.17).

Denote by $\pi_\e$ the shortest distance projection on the ball $\overline B(0,1+\e)$. Then $\pi_\e$ is 1-Lipschitz, and $\pi
_\e(\R^4\bs B(0,1+\e))\subset\partial B(0,1+\e)$. On the other hand. $\pi$ is $\frac{1}{1+\e}$-Lipschitz on $\partial B(0,1+\e)$. So if $E\subset \R^4\bs B(0,1+\e)$, then $\pi=\pi\circ\pi_\e$ is $\frac{1}{1+\e}$-Lipschitz on $E$. \qed

Let us return to our sets $F_l$. We have $B(x,\frac 12\e)\cap B(0,1+\frac12\e)=\emptyset$, since $x\not\in B(0,1+\e)$, and hence
\be\begin{split}
H^2(\pi(F_l))&\le H^2(\pi(F_l\bs B(x,\frac12\e)))+H^2(\pi(F_l\cap B(x,\frac12\e)))\\
&\le H^2 (F_l\bs B(x,\frac12\e))+\frac{1}{(1+\frac12\e)^2}H^2(F_l\cap B(x,\frac12\e))\\
&\le H^2(F_l)-\frac{4\e+\e^2}{4+4\e+\e^2}H^2(F_l\cap B(x,\frac12\e)) \\
&\le H^2(F_l)-\frac{4\e+\e^2}{4+4\e+\e^2}\frac{\e^2}{4C}\\
&=H^2(F_l)-C(\e),
\end{split}\ee
where $C(\e)>0$ for all $\e>0$, and $C(\e)$ does not depend on $l$ for $l$ large.

We know that $\{F_l\}$ is a minimizing sequence, therefore for all $\e>0$, there exists $N>0$ such that for $l>N$, 
\be H^2(F_l)\le\inf_{E\in\mathfrak F}H^2(E)+\frac 12C(\e)< H^2(\pi(F_l))+C(\e).\ee 
Hence (4.18) is not true, which means that $F_l\subset B(B(0,1)\cup P_k,\e)\cap U$, thus we obtain (4.13).

As a result, since $E_k$ is the limit of $F_l$, 
\be E_k\subset\cap_\e B(B(0,1)\cup P_k,\e)\cap U\subset\overline B(0,1).\ee

\smallskip

\begin{lem}
Let $C\subset \R^n$ be a closed convex symmetric (with respect to the origin) set with non-empty interior. Then for all $\e>0$, there exists $\d>0$ and a 1-Lipschitz retraction $f$ of $\R^n$ onto $C$ such that $f$ is $1-\d$-Lipschitz on $\R^n\bs B(C,\e)$.
\end{lem}

\nd Denote by $||\cdot||_C$ the norm on $\R^n$ whose closed unit ball is $C$. Then there exists $A>1$ such that
\be A^{-1}||\cdot||_C\le ||\cdot||\le A||\cdot||_C,\ee
where $||\cdot||$ denotes the Euclidean norm.

For all $a\ge 0$, set $||\cdot||_a=||\cdot||_C+a||\cdot||$ and $C_{a,b}$ the closed ball of radius $b$ under the norm $||\cdot||_a$. Then $C_{0,1}=C$. Notice that for all $x\in \R^n\bs\{0\}$, $||x||_a$ is a strictly increasing continuous function of $a$, $||\cdot||_0=||\cdot||_C$, and that $C_{a,b}$ is continuous, decreasing with respect to $a$ and increasing with respect to $b$, that is, 
\be C_{a,b}\supset C_{a',b},\ C_{a,b}\subset C_{a,b'}\mbox{ for all }a<a',b<b',\ee
and
\be \bigcap_{a_n\to a-}C_{a_n,b}=\bigcap_{b_n\to b+}C_{a,b_n}=C_{a,b};\bigcup_{a_n\to a+}C_{a_n,b}=\bigcup_{b_n\to b-}C_{a,b_n}=C_{a,b}^\circ.\ee

Since the norm defined by $C_{a,b}$ contains a part of Euclidean norm, which is uniformly convex, it is easy to verify that 
\be\begin{split}\mbox{for all }&a,b>0, \mbox{ there exists a constant }M(a,b,A)>0,\mbox{ such that}\\
&\mbox{ for each }x,y\in \partial C_{a,b}\mbox{ with }\a_{x,y}<\frac\pi2,
 B(\frac{x+y}{2}, M(a,b,A)||x-y||^2)\subset C_{a,b},\end{split}\ee
%&||\frac{x+y}{2}||_a\le b-D(a,b,A)||x-y||^2,\end{split}\ee
where $\a_{x,y}<\pi$ denotes the angle between $\vec{Ox}$ and $\vec{Oy}$ for $x,y\ne 0$, and $B(\frac{x+y}{2}, M(a,b,A)||x-y||^2)$ denotes the euclidean ball centered at $\frac{x+y}{2}$ with radius $M(a,b,A)||x-y||^2$.

Now for all $\e>0$, let $w,v$ be two points in $\R^n\bs B(C_{a,b},\e)$ such that $\pi_{a,b}(w)=x,\pi_{a,b}(v)=y$, where $\pi_{a,b}$ denotes the shortest distance projection on the convex set $C_{a,b}$. We claim that the angle $\beta_1\in[0,\frac\pi 2]$ between $\vec{xw}$ and $\vec{yx}$ is smaller than $\arctan\frac{1}{2M(a,b,A)||x-y||}.$ (See the picture 4-1 below). In fact, denote by $P$ the plane containing $x,y$ and $w$, denote by $z\in P$ the point such that $[z,\frac{x+y}{2}]\perp [x,y]$ and $||z-\frac{x+y}{2}||=M(a,b,A)||x-y||^2$. Then $z\in C_{a,b}$, $[x,z]\in C_{a,b}$, and
\be\tan\angle zxy=2M(a,b,A)||x-y||.\ee

Then if $\beta_1>\arctan\frac{1}{2M(a,b,A)||x-y||}$, we have $\angle wxz<\frac\pi2$. Denote by $s$ the projection of $w$ on $L'$, the line passing through $x$ and $z$. Then $s$ is between $x$ and $z$, or $z$ is between $x$ and $s$. In both cases $(x,z)\cap (x,s)\ne \emptyset$. Take $x'\in  (x,z)\cap (x,s)\subset C_{a,b}$, then $x'\in C_{a,b}$, and $\angle wxz<\angle wx'z$. As a result
\be ||w-x||=||w-s||/\sin\angle wxz>||w-s||/\sin\angle wx'z=||w-x'||,\ee
which contradicts the fact that $x$ is the shortest distance projection of $w$ on $C_{a,b}.$

\noindent\centerline{\includegraphics[width=0.8\textwidth]{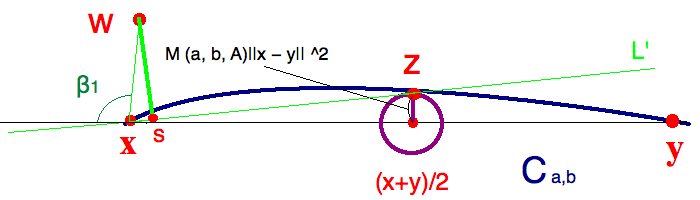}}
\nopagebreak[4]
\centerline{4-1}

Similarly we can prove that if $\beta_2$ denotes the angle between $\vec{yv}$ and $\vec{xy}$, then $\beta_2\le \arctan\frac{1}{2M(a,b,A)||x-y||}.$

Denote by $L$ the line passing through $x$ and $y$, and $p_L$ the orthogonal projections on them. Then
\be ||w-v||\ge ||p_L(w)-p_L(v)||=||w-x||\cos\beta_1+||x-y||+||v-y||\cos\beta_2.\ee
But $w,v\in\R^n\bs B(C_{a,b},\e)$, hence we have $||w-x||>\e,||v-y||>\e$, therefore
\be \begin{split}||w-v||&\ge ||x-y||+2\e\cos\arctan\frac{1}{2M(a,b,A)||x-y||}\\
%&=||x-y||+2\e(\frac{4M^2||x-y||^2}{1+4M(a,b,A)^2||x-y||^2})^\frac12\\
%&\ge ||x-y||+2\e\frac{2M(a,b,A)||x-y||}{(1+4M(a,b,A)^24R^2)}\\
&=(1+\e C(a,b,A))||x-y||,\end{split}\ee
with $C(a,b,A)>0$.

Notice that (4.29) is true for all pairs of $x,y$ such that $\a_{x,y}<\frac\pi2$. Hence $\pi_{a,b}$ is locally $\frac{1}{1+\e C(a,b,A))}$-Lipschitz on $\R^n\bs B(C_{a,b},\e)$.

Return to the proof of the lemma. Fix an arbitrary $\e>0$. Then by (4.23) and (4.24), there exists $a,b>0$ such that
\be C\subset C_{a,b}\subset B(C_{a,b},\frac\e2)\subset B(C,\e).\ee

Now denote by $\pi_C$ the shortest distance projection on $C$. Denote by $f=\pi_C\circ \pi_{a,b}$ for a pair of $a,b$ which satisfies (4.30). We want to prove $f$ is the desired map for Lemma 4.21. Since $\pi_C$ is 1-Lipschitz, for proving the lemma, it is sufficient to prove that $\pi_{a,b}$ is locally $1-\d$-Lipschitz on $\R^n\bs B(C_{a,b},\frac\e2)$. Then by (4.29), we take $\d$ such that $1-\d=\frac{1}{1+\frac12 \e C(a,b,A))}$, and we obtain the conclusion.\qed

\begin{cor}
Let $E\subset \R^n\bs B(C,\e)$ be a rectifiable set and $f$ be as in the lemma, then
\be H^d(f(E))\le (1-\d)^dH^d(E).\ee
\end{cor}\qed
 
Now let us return to the proof of Proposition 4.8. Recall that here $C$ is the convex hull of $P_k\cap \overline B(0,1)$.

\noindent\textbf{Proof of (2) and (5) of Proposition 4.8.}

We will prove that
\be \mbox{ for all }\e>0, E_k\subset B(C,2\e).\ee
If it is not true, i.e. if $E_k\bs B(C,2\e)\ne\emptyset$, then by Ahlfors regularity, $H^2(E_k\bs B(C,2\e))>0$. We apply Corollary 4.31 to $E_k\bs B(C,2\e)$ and the convex set $B(C,\e)$, and obtain that there exists a Lipschitz map $f_\e$ of $\R^n$ in $B(C,\e)$, such that
\be H^2(f_\e(E_k))<H^2(E_k),\ee
where $f_\e$ is as in Lemma 4.21. Then we will get (4.33) if we know that $f_\e(E_k)$ is a deformation of $E_k$ in $U$. 

But we already know that $E_k\subset\overline B(0,1)$. So $W_\e:=\{x\in E_k,\pi_\e(x)\ne x\}$ is contained in $\overline B(0,1)\bs B(C,\e)$, which is compact and far from the boundary of $U$. Denote by $d(W_\e,\partial U)=d$, and set
\be g:\R^n\to B(C,\e),f(x)=t(x)x+(1-t(x))f_\e(x),t(x)=\min\{2d(x,W_\e)/d,1\},\ee
then $g$ is a deformation in $U$, and
\be g(E_k)=f_\e(E_k).\ee

Therefore, if $E_k\bs B(C,2\e)\ne\emptyset$, then $g$ decreases strictly the measure of $E_k$, which contradicts the fact that $E_k$ is minimal. Hence we have (4.33). But (4.33) is true for all $\e>0$, hence we have
\be E_k\subset C,\ee
from which (5) follows. (2) is a direct corollary of (5). This completes the proof of Proposition 4.8.\qed
 
 \begin{cor}For each subsequence $\{n_k\}$ such that $E_{n_k}$ converges in $\overline B(0,1)$ for the Hausdorff distance, the limit is $P_0\cap \overline B(0,1)$.
 \end{cor}
 
 \nd Take such a sequence $\{n_k\}$. Denote by $E_\infty$ the limit of $E_{n_k}$. We want to apply Theorem \ref{unicite}, so we should check its hypotheses for $E_\infty$.
 
First, by \cite{GD03} Thm 4.1, which says that the limit under the Hausdorff distance of a sequence of minimal sets is minimal, we know that $E_\infty$ is also minimal in $B(0,1)$, because each $E_{n_k}$ is. 

Next, (3.2) follows by Proposition 4.8(3) and the fact that $E_\infty$ is the limit of $E_{n_k}$.
 
To verify (3.3), if we denote by $C_k$ the convex hull of $P_k\cap \overline B(0,1)$ (and $C_0$ the convex hull of $P_0\cap\overline B(0,1)$), then by Proposition 4.8, $E_{n_k}\subset C_{n_k}$. Since $E_\infty$ is the limit of $E_{n_k}$, for each $m>0$, there exists $K(m)>0$ such that for all $k>K(m)$, $E_\infty\subset B(E_{n_k},\frac1m)\subset B(C_{n_k},\frac 1m)\subset\cup_{k=K(m)}^\infty B(C_{n_k},\frac 1m)$. We can also ask that $K(m)>k(l)$ when $m>l$. As a result we have
  \be E_\infty\subset\cap_{m=1}^\infty\cup_{k=K(m)}^\infty B(C_{n_k},\frac 1m)=C_0,\ee
  and therefore 
  \be E_\infty\cap \partial B(0,1)\subset C_0\cap\partial B(0,1)=P_0\cap\partial B(0,1).\ee
  
On the other hand, since $E_\infty$ is the limit of the $E_{n_k}$, it contains $\lim_{k\to\infty}E_{n_k}\cap\partial B(0,1)=P_0\cap \partial B(0,1)$.
  
Hence
 \be E_\infty\cap \partial B(0,1)=P_0\cap\partial B(0,1).\ee which gives (3.3).
   
For (3.4), Proposition 4.8 (4) gives that that $H^2(E_{n_k})<2\pi$. But $E_{n_k}$ and $E_\infty$ are minimal sets, by the lower semi continuity of Hausdorff measures for minimal sets (\cite{GD03} Thm 3.4),
\be H^2(E_\infty\cap B(0,1))\le\liminf_k H^2(E_{n_k}\cap B(0,1))\le 2\pi.\ee

The equality follows from (3.2), and Lemmas \ref{projection} and 2.14.

Hence by Theorem \ref{unicite}, $E_\infty=P_0\cap\overline B(0,1)$.
\qed
 
Notice that $\overline B(0,1)$ is compact, hence there exists a subsequence of $\{E_k\}$ that converges. We still denote this subsequence by $\{E_k\}$ for short. Then by Corollary 4.38, the limit is $P_0\cap \overline B(0,1)$. So from now on, we will concentrate on this sequence $\{E_k\}_{k>0}$ that converges to $P_0\cap\overline B(0,1)$.

\bigskip

\section{Critical radius}

For our sequence $\{E_k\}$, once again, if for each $k$, $E_k$ is a deformation of $P_k$, then there should be some sort of pinch at the center, because otherwise, the deformation is injective, and $E_k$ is therefore an essentially disjoint union of images of $P_k^1$ and $P_k^2$. But since $P_k^i$ is minimal, $i=1,2$, its image always has a larger measure than $P_k^i$, therefore $E_k$ is not a better competitor.

Now since we have to pinch, we know that $E_k$ has to get far away from $P_k$ somewhere, and intuitively we have to pay a price for this. In fact, the main issue of the proof is to understand why a pinch costs more than all that it could save for us.

However $E_k$ is not necessarily a deformation of $P_k$, so what we said just now is just for giving an idea on what we are going to do, that is, to find the place where $E_k$ begins to get away from $P_k$.

So for $\e$ sufficiently small, we want to find a center $o$, which is not far from the origin, and a scale $r$, such that $E_k$ is $\e r'$ near some translation of $P_k$ in $B(o,r')$ for all $r'\ge 2r$, but is not $\e r$ near any translation of $P_k$ in $B(o,r)$. Then this is the place that we are looking for, and $r$ is the critical radius. 

As a result, in the small ball $B(o,\frac14r)$ where the pinch takes place, probably we cannot see very clearly what happens, so we only control the measure of $E_k$ by an argument of projection, using Corollary 2.45. However, for the part outside the small ball, $E_k$ is near the two planes, and hence is regular, by the regularity property of minimal sets. So we can treat this part more precisely. And finally we will be able to prove that pinching will make us lose more measure on these regular parts than that we can gain in the small ball.

So let's begin by looking for the critical radius in this section, by a stopping time argument.

\bigskip

For each, $k\in \N$ and $i=1,2$, set
\be C_k^i(x,r)=(p^i_k)^{-1}(B(0,r)\cap P^i_k)+x,\ee
where $p_k^i$ is the orthogonal projection on $P_k^i$, and
\be D_k(x,r)=C_k^1(x,r)\cap C_k^2(x,r).\ee
So $C_k^i(x,r)$ is a cylinder and $D_k(x,r)$ is the intersection of two cylinders. It is not hard to see that $D_k(x,r)\supset B(x,r)$ and $D_k(0,1)\cap P_k=B(0,1)\cap P_k$. 

We say that two sets $E,F$ are $\e r$ near each other in an open set $U$ if 
\be d_{r,U}(E,F)<\e,\ee
where 
\be d_{r,U}(E,F)=\frac 1r\max\{\sup\{d(y,F):y\in E\cap U\},\sup\{d(y,E):y\in F\cap U\}\}.\ee

%Et quand $U$ is fix\'e avant, on \'ecrit parfois $d_{x,r}$ au lieu de $d_{x,r,U}$ quand il n'y a pas d'ambigu\'it\'e.

We set also
\be \begin{split}&d^k_{x,r}(E,F)=d_{r,D_k(x,r)}(E,F)\\
     &=\frac 1r\max\{\sup\{d(y,F):y\in E\cap D_k(x,r)\},\sup\{d(y,E):y\in F\cap D_k(x,r)\}\}.\end{split}\ee

\begin{rem}We should be clear about the fact that 
\be d_{r,U}(E,F)\ne\frac 1rd_H(E\cap U,F\cap  U).\ee
To see this, we can take $U=D_k(x,r)$, and set $E_n=\partial D_k(x,r+\frac 1n)$ and $F_n=\partial D_k(x,r-\frac 1n)$. Then we have 
\be d^k_{x,r}(E_n,F_n)\to 0\ee and \be \frac 1rd_H(E_n\cap D_k(x,r),F_n\cap D_k(x,r))=\frac 1rd_H(E_n\cap D_k(x,r),\emptyset)=\infty.\ee
So $d_{r,U}$ measures rather how the part of one set in the open set $U$ could be approximated by the other set, and vice versa. However we always have
\be d^k_{x,r}(E,F)\le\frac 1rd_H(E\cap D_k(x,r),F\cap D_k(x,r)).\ee
\end{rem}

\bigskip

Now let us recall that $\{E_k\}$ is a sequence of sets as in Proposition 4.8,  with $\theta_k>\frac\pi2-\frac1k$, which converges to $P_0\cap\overline B(0,1)$. 

\begin{pro}\label{rk}There exists $\e_0>0$, such that if $\e<\e_0$, then for all $k$ large, there exists $r_k\in]0,\frac12[$ and $o_k\in B(0,12\e)$ such that $E_k$ is $2\e r_k$ near $P_k+o_k$ in $D_k(o_k, 2r_k(1-12\e))$, but not $\e r_k$ near $P_k+q$ in $D_k(o_k,r_k)$ for any $q\in\R^4$.
\end{pro}

\begin{rem}We will also use the construction for information about intermediate scales in the proof.
\end{rem}

 \nd

We fix $\e$ and $k$, and set $s_i=2^{-i}$ for $i\ge 0$. Set $D(x,r)=D_k(x,r),d_{x,r}=d_{x,r}^k$ for short. 

We  proceed in the following way.

Step 1: Denote by $q_0=q_1=O$, then in $D(q_0,s_0)$, $E_k$ is $\e s_0$ near $P_k+q_1$ if $k$ is large, since $E_k\to P_0$ and $P_k\to P_0$ and hence $d_{0,1}(E_k,P_k)\to 0$.

Step 2: If in $D(q_1,s_1)$ $E_k$ is not $\e s_1$ near $P_k+q$ for any $q$, we stop; if not, there exists a $q_2$ such that $E_k$ is $\e s_1$ near $P_k+q_2$ in $D(q_1,s_1)$. Here we also ask that $\e$ be small enough (say, $\e<\frac{1}{100}$) so that $q_2\in D(q_1,\frac12s_1)$, thanks to the conclusion of step 1. Then in $D(q_1,s_1)$, we have simultaneously : 
\be d_{q_1,s_1}(E_k,P_k+q_1)\le s_1^{-1} d_{q_0,s_0}(E_k,P_k+q_1)\le 2\e\ ;\ d_{q_1,s_1}(E_k, P_k+q_2)\le\e.\ee

Let us verify that (5.13) implies that $d_{q_1,\frac12 s_1}(P_k+q_1,P_k+q_2)\le 12\e$ when $\e$ is small, say, $\e<\frac {1}{100}$. 
In fact, for each $z\in D(q_1,\frac12 s_1)\cap (P_k+q_1)$, we have $d(z,E_k)\le d_{q_0,s_0}(E_k,P_k+q_1)\le\e$, hence there exists $y\in E_k$ such that $d(z,y)\le\e$. But since $z\in D(q_1,\frac12 s_1)$, we have $y\in D(q_1,\frac12 s_1+\e)\subset D(q_1,s_1)$, and hence $d(y,P_k+q_2)\le s_1^{-1} d_{q_1,s_1}(E_k,P_k+q_2)\le 2\e$, therefore $d(z,P_k+q_2)\le d(z,y)+d(y,P_k+q_2)\le 3\e$.

On the other hand, suppose $z\in D(q_1,\frac12 s_1)\cap (P_k+q_2)$, we have $d(z,E_k)\le s_1^{-1}d_{q_1,s_1}(P_k+q_2,E_k)\le 2\e$, hence there exists $y\in E_k$ such that $d(z,y)\le 2\e$. But since  $z\in D(q_1,\frac12 s_1)$, we have $y\in D(q_1,\frac12 s_1+2\e)\subset D(q_0,s_0)$, and hence $d(y,P_k+q_1)\le d_{q_0,s_0}(E_k,P_k+q_1)\le \e$, which implies $d(z,P_k+q_1)\le d(z,y)+d(y,P_k+q_1)\le 3\e$.

As a result 
\be d_{q_1,\frac12 s_1}(P_k+q_1,P_k+q_2)\le (\frac12 s_1)^{-1}\times 3\e=12\e,\ee
hence $d_{q_1,\frac12 s_1}(q_1,q_2)\le 24\e$, and therefore $d(q_1,q_2)\le 6\e=12\e s_1$.

Now we define our iteration process (notice that it depends on $\e$, so we also call it a $\e$-process).

Suppose that all $\{q_i\}_{i\le n}$ have already been defined, with
\be d(q_i,q_{i+1})\le 12s_i\e=12\times 2^{-i}\e\ee
for $0\le i\le n-1$, and hence
\be d(q_i,q_j)\le 24\e s_{\min(i,j)}=2^{-\min(i,j)}\times 24\e\ee 
for $0\le i,j\le n$, moreover for all $i\le n-1$, $E_k$ is $\e s_i$ near $P_k+q_{i+1}$ in $D(q_i, s_i)$. We say that the process does not stop at step $n$. In this case

Step n+1 : We look at the situation in $D(q_n, s_n)$.

If $E_k$ is not $\e$ near any $P_k+q$ in this ball of radius $s_n$, we stop, since we have found the $o_k=q_n, r_k=s_n$ as desired. In fact, since $d(q_{n-1}, q_n)\le 12\e s_{n-1}$, we have
$D(q_n, 2s_n(1-12\e))=D(q_n, s_{n-1}(1-12\e))\subset D(q_{n-1}, s_{n-1})$, and hence
\be\begin{split}
d_{q_n,2s_n(1-12\e)}(P_k+q_n, E_k)&\le (1-12\e)^{-1}d_{q_{n-1},s_{n-1}}(P_k+q_n, E_k)\\
&\le \frac{\e}{1-12\e}.
\end{split}\ee
Moreover 
\be d(o_k, O)=d(q_n,q_1)\le 2^{-\min(1,n)}\times 24\e=12\e.\ee

Otherwise, we can find a $q_{n+1}\in \R^4$ such that $E_k$ is still $\e s_n$ near $P_k+q_{n+1}$ in $D(q_n, s_n)$. Then since $\e$ is small, $q_{n+1}\in D(q_n,\frac12 s_n)$. Moreover we have as before $d(q_{n+1},q_n)\le 12\e s_n$, and for $i\le n-1$,
\be d(q_i,q_{n+1})\le \sum_{j=i}^n d(q_j,q_{j+1})\le\sum_{j=i}^n 12\times 2^{-j}\e\le 2^{-j}\times 24\e=2^{-\min(i,n+1)}\times 24\e.\ee
Thus we have obtained our $q_{n+1}$.  

Now all we have to do is to prove that for every $\e$ small enough, this process has to stop at a finite step. For this purpose, we are going to estimate the measure of $E_k$. So we need the lemma below.

\begin{lem}\label{fn}There exists $\e_0>0$, such that for all $\e<\e_0$, $k$ large enough, and for every $n$ such that the $\e-$process does not stop before $n$ (which means in particular that there exists  $q_n\in B(q_{n-1},\frac12s_{n-1})$ such that $E_k$ is $\e s_{n-1}$ near $P_k+q_n$ in $D(q_{n-1},s_{n-1})$),  
\be E_k\cap (D(0,1)\backslash D(q_n,s_n))=F_n^1\cup F_n^2\ee
where $F_n^1, F_n^2$ do not meet each other. Moreover
\be P_k^i\cap (D(0,1)\backslash D(q_n,s_n))\subset p_k^i(F_n^i)\ee
where $p_k^i$ is the orthogonal projection on $P_k^i,i=1,2$.
\end{lem}

We will prove a more general version of this lemma in the next section (Proposition 6.1 (2)). So let us admit it for the moment.

Since $H^2(E_k)<2\pi$, there exists $n_k>0$ such that
\be\inf_{q\in\R^4}  H^2(P_k\backslash D(q, s_{n_k}))>H^2(E_k).\ee

Then our process has to stop before step $n_k$, because otherwise by the above lemma, we have the disjoint decomposition
\be E_k=[E_k\cap D(q_{n_k},s_{n_k})]\uplus F_{n_k}^1\uplus F_{n_k}^2,\ee
and hence
\be\begin{split}
H^2(E_k)&\ge H^2(F_{n_k}^1)+H^2(F_{n_k}^2)\ge H^2[p_k^1(F_{n_k}^1)]+H^2[p_k^2(F_{n_k}^2)]\\
&\ge H^2(P_k\backslash D(q_{n_k}, s_{n_k})>H^2(E_k).
\end{split}\ee
This is impossible.
\qed

\section{Projection properties and regularity of $E_k$}

As we have said in the previous section, the next step is to give some useful properties of $E_k$, including the regularity for the flat part of $E_k$ out of the small critical ball, and the surjectivity of the projections of $E_k$ inside the ball. This is the main aim of this section. It gives also the proof of Lemma \ref{fn}. 

We are sorry that the proof for (3) and (4) of Proposition 6.1 is surprisingly painful. We have to derive the property from the geometric construction of the proof for Theorem 4.1, which is already complicated. But we did not find any easier proof.

So let us first state the main proposition that we will prove in this section.

\begin{pro}There exists $\e_0>0$, such that for all $\e<\e_0$ and $k$ large, if the $\e-$process does not stop before the step $n$, then

(1) $E_k\cap (D_k(0,\frac {39}{40})\backslash D_k(q_n,\frac {1}{10}s_n))$ is composed of two disjoint pieces $G^i,i=1,2$, such that:
\be
G^i\mbox{ is the graph of a }C^1\mbox{ map }\ g^i:D_k(0,\frac{39}{40})\backslash D_k(q_n,\frac{1}{10}s_n)\cap P^i_k\to {P^i_k}^\perp\ee
with
\be||\nabla g^i||_\infty<1;\ee

(2)(A more general version of Lemma \ref{fn}) for every $\frac {1}{10}s_n\le t\le s_n$,
\be E_k\cap (D_k(0,1)\backslash D_k(q_n,t))=G_t^1\cup G_t^2\ee
where $G_t^1, G_t^2$ do not meet each other. Moreover
\be P_k^i\cap (D_k(0,1)\backslash C_k^i(q_n,t))\subset p_k^i(G_t^i)\ee
where $p_k^i$ is the orthogonal projection on $P_k^i,i=1,2$;

(3) for each $\frac{1}{10}s_n<t<s_n$, there exists a sequence $\{F_l^n(t)=f_l^n(t)((P_k+q_n)\cap \overline D_k(q_n,t+\frac 1l))\}_{l\ge 1}$ of deformations of  $(P_k+q_n)\cap \overline D_k(q_n,t+\frac 1l)$ in $\overline D_k(q_n,t+\frac 1l)$, with
\be f_l^n(t)((P_k+q_n)\cap\partial C_k^i(q_n,t+\frac 1l))\subset \partial C_k^i(q_n,t+\frac 1l),\ee 
that converge to $E_k\cap D_k(q_n,t)$ in $D_k(0,1)$;

(4) the orthogonal projections $p_k^i: E_k\cap D_k(q_n,t)\to P_k^i\cap C_k^i(q_n,t),i=1,2$ are surjective, for all $\frac {1}{10}s_n\le t\le s_n$.
\end{pro}

In order to prove (1) we will apply a regularity theorem on varifolds. First we give some useful notations below.

 $G(n,d)$ is the Grassmann manifold $G(\R^n,d)$; 
 
for every $T\in G(n,d)$, we denote by $\pi_T$ the orthogonal projection on the $d$-plane represented by $T$;
 
 for every measure $\nu$ on $\R^n$, $\theta^d(\nu,x)=\lim_{r\to 0}\frac{\nu B(a,r)}{\a(d)r^d}$ (if the limit exists) is the density of $\nu$ on $x$, where $\a(d)$ denote the volume of the $d$-dimensional unit ball; 
 
 $\mathbb V_d(\R^n)$ denote the set of all $d-$varifold in $\R^n$, i.e. all Radon measures on $G_d(\R^n)=\R^n\times G(n,d)$; 
 
 for each $V\in \mathbb V_d(\R^n),$ $||V||$ is the Radon measure on $\R^n$ such that for each $A\subset \R^n$, $||V||(A)=V(G_d(\R^n)\cap \{(x,S):x\in A\})$; 
 
 $\delta(V)$ denotes the first variation of $V$, that is, the linear map from $\mathfrak X(R^n)$ to $\R$, defined by
\be\delta V(g)=\int Dg(x)\cdot \pi_S dV(x,S)\ee for $g\in\mathfrak X(\R^n).$
Here $\mathfrak X (\R^n)$ is the vector space of all $C^\infty$ maps from $\R^n$ to $\R^n$ with compact support.

In our case, we are only interested in rectifiable varifolds. In fact, with each $d-$rectifiable set $E$ we associate a $d-$varifold, denoted by $V_E$, in the following sense : for each $B\subset \R^n\times G(n,d)$, we have 
\be V_E(B)=H^d\{x:(x,T_xE)\in B\}.\ee
Recall that $T_xE$ is the $d$-dimensional tangent plane of $E$ at $x$, it exists for almost all $x\in E$, because $E$ is $d-$rectifiable. Then $||V_E||=H^d|_E$. Moreover, the density $\theta^d(||V_E||,x)$ exists for almost all $x\in E$.

\begin{thm}[c.f.\cite{All72} Regularity theorem at the beginning of section 8]Suppose $2\le d<p<\infty,$ $q=\frac{p}{p-1}$. Corresponding to every $\e$ with $0<\e<1$ there is $\eta>0$ with the following property:

Suppose $0<R<\infty$, $0<\mu<\infty$, $V\in\mathbb{V}_d(\R^n)$, $a\in spt||V||$ and

1) $\theta^d(||V||,x)\ge \mu$ for $||V||$ almost all $x\in B(a,R)$;

2) $||V||B(a,R)\le(1+\eta)\mu\a(d)R^d ;$

3) $\delta V(g)\le\eta \mu^{\frac 1p}R^\frac{d}{p-1}(\int|g|^q\mu||V||)^\frac 1q$ whenever $g\in {\mathfrak X}(\R^n)$ and $spt\ g\subset B(a,R)$.

Then there are $T\in G(n,d)$ and a continuously differentiable function $F:T\to \R^n$, such that $\pi_T\circ F=1_T$,
\be||DF(y)-DF(z)||\le\e(|y-z|/R)^{1-\frac dp}\mbox{ whenever }y,z\in T,\ee and
\be B(a,(1-\e)R)\cap spt||V||=B(a,(1-\e)R)\cap\ image\ F.\ee
\end{thm}

\begin{rem}

1) In the theorem, since $\pi_T\circ F=1_T$,  we can see that $F$ is in fact the graph of a $C^1$ function $f$, defined by $f(t)=\pi_{T^\perp} F(t)$, with $t\in T$, $\pi_{T^\perp}$ the orthogonal projection on the orthogonal space $T^\perp$ of $T$. Moreover $||Df(t)||\le ||DF(t)||$ for all $t\in T$.

2) If $E$ is a locally minimal set, then $V_E$ is stationary, i.e. $\delta V_E=0$. Hence the condition 3) is automatically true. In fact if we set $g_t(x)=(1-t)x+tg(x)$, then 
\be\delta V_E(g)=\frac{d}{dt} H^d(g_t(E\cap spt g)),\ee 
which can be deduced from the area formula. Thus if $E$ is minimal, $\delta V_E=0$. 
\end{rem}

\begin{pro}For all $n>d>0$, there exists $\e_1=\e_1(n,d)>0$ such that the following is true.  Let $E$ be a locally minimal set of dimension $d$ in an open set $U\subset \R^n$, with $U\supset B(0,2)$ and $0\in E$. Then if $E$ is $\e_1$ near a $d-$plane $P$ in $B(0,1)$, then $E$ coincides with the graph of a $C^1$ map $f:P\to P^\perp$ in $B(0,\frac34)$. Moreover $||\nabla f||_\infty<1$.
\end{pro}

\nd We will prove it only for $d=2$. Proofs for other dimensions are similar. 

First let us verify the conditions in Theorem 6.9, with $a=0,R=1,\mu=1,$ $\eta<\frac{1}{10}$, and $\e$ small, to choose later.

1) Since $E$ is minimal, for each $x\in E$, the density of $E$ at $x$ is at least $1$, hence 1) is true.

2) We know that $E$ is $\e_1$ near a 2-dimensional affine plane $P$ in $B(0,1)$, and $H^2(P\cap B(0,1))\le\a(2)=\pi$. Then by Lemma 16.43 in \cite{DJT}, we can choose $\e_1$ (which depends on $\e$, since $\eta$ depends on $\e$) such that 2)  is true. In particular
\be H^2(B(0,1)\cap E)\le \frac{11}{10}\pi.\ee

3) comes from the minimality of $E$, by Remark 6.12 2), with any $p>2$.

Then when $p$ is sufficiently large, by Theorem 6.9, there exists a plane $T$ and a $C^1$ map $F$ from $T$ to $\R^4$ such that in $B(0,1-\e)$, $E$ coincides with the image of $F$, and
 \be ||DF(y)-DF(z)||\le \e(|y-z|/R)^{1-\frac 2p}\le (\frac 74)^{1-\frac 2p}\e\le 2\e\ee
for all $y,z\in F^{-1}(E\cap B(0,1-\e)).$

Notice that for each $y\in F^{-1}(E\cap B(0,1-\e))\subset T$, $DF(y)(T)$ is the tangent plane of $E$ at $F(y)$. Thus (6.16) means that the tangent planes of $E$ do not vary  a lot in $B(0,1-\e)$. 

Now denote by $Q$ the plane parallel to $P$ (the plane in the statement of Proposition 6.14) and passing by the origin, $\pi$ the projection on $Q$, and $\pi'$ the projection on $Q^\perp$. We claim that for each $y\in T$ such that $F(y)\in B(0,1-\e)$ and each $u\in T_yE=DF(y)(T)$, 
\be||\pi(u)||\ge\frac 34||u||.\ee

In fact, if (6.17) is not true, then there exists $y\in F^{-1}(E\cap B(0,1-\e))$ and $u\in T_yE$ such that $||\pi(u)||<\frac 34||u||.$ Denote by $t=DF(y)^{-1}(u)\in T$, then

\be ||\pi(DF(y)(t))||\le \frac 34||DF(y)(t)||.\ee 
By (6.16), for all $z\in F^{-1}(E\cap B(0,1-\e))$,
\be \begin{split}||\pi(DF(z)(t)||&\le ||\pi\circ(DF(z)-DF(y))(t)||+||\pi\circ DF(y)(t)||\\
&\le 2\e||t||+\frac 34||DF(y)(t)||\\
&\le 2\e||t||+\frac 34||DF(y)(t)-DF(z)(t)||+\frac 34||DF(z)(t)||\\
&\le 2\e||t||+\frac32\e||t||+\frac 34||DF(z)(t)||\\
&=\frac 72\e||t||+\frac 34||DF(z)(t)||.\end{split}\ee
But $\pi_T\circ F=1_T$ implies that for all $t\in T,$ \be||DF(z)(t)||\ge ||\pi_T\circ DF(z)(t)||=||t||,\ee 
therefore
\be ||\pi(DF(z)(t)||\le (\frac72\e+\frac 34)||DF(z)(t)||.\ee

So when $\e$ is sufficiently small, we have 
\be ||\pi'\circ DF(z)(t)||\ge \frac 12||DF(z)(t)||\ee
for all $z\in F^{-1}(E\cap B(0,1-\e))$.

Set $e_1=\pi'\circ DF(0)(t)/||\pi'\circ DF(0)(t)||$ a unit vector in $Q$. Then we have
\be <e_1,DF(0)(t)>>\frac 12||DF(0)(t)||.\ee

Then for all $z\in F^{-1}(E\cap B(0,1-\e))$, still by (6.16) and (6.20),
\be \begin{split}<e_1,DF(z)(t)>&=<e_1,(DF(z)-DF(0))(t)>+<e_1,DF(0)(t)>\\
&\ge <e_1,DF(0)(t)>-||<e_1,(DF(z)-DF(0))(t)>||\\
&\ge \frac 12||DF(0)(t)||-2\e||t||\\
&\ge \frac 12[||DF(z)(t)||-||(DF(0)-DF(z))(t)||]-2\e||t||\\
&\ge \frac 12[||DF(z)(t)||-2\e||t||]-2\e||t||\ge \frac 12||DF(z)(t)||-3\e||DF(z)(t)||\\
&\ge\frac 13||DF(z)(t)||.\end{split}\ee

As a result, if we take $z\in F^{-1}(E\cap \partial B(0,1-\e))$, such that $\vec z=\lambda t$ with $\lambda>0$, we have
\be \begin{split}<e_1,F(z)-F(0)>&=<e_1,\int_0^1 DF(sz)(t)ds>=\int_0^1<e_1,DF(sz)(t)>ds\\
&\ge\int_0^1\frac 13||DF(sz)(t)||ds=\frac 13\int_0^1||DF(sz)(t)||ds\\
&\ge \frac13||F(z)-F(0)||=\frac 13||F(z)||=\frac13(1-\e),\end{split}\ee
which implies that when $\e_1,\e$ are sufficiently small, there exists no translation $Q+x$ of $Q$ (including $P$) such that $B_{0,1}(Q+x,E)<\e_1$. Contradiction.

So we have (6.17). In other words, $D\pi$ is always injective in $E\cap B(0,1-\e)$. Then by the implicit function theorem, for all $x\in E\cap B(0,1-\e)$, there exists $r_x>0$ and $g_x:Q\to Q^\perp$ such that in $\pi^{-1}[B(\pi(x),r_x)\cap Q]\cap B(x,2r_x)$, $E$ coincides with the graph of $g_x$ on $B(\pi(x),r_x)$. Moreover by (6.17) 
\be ||\nabla g_x(x)||\le 1.\ee

Next let us verify that
\be \pi(E\cap B(0,1-\e))\supset Q\cap B(0,\frac 34).\ee 

Recall that $E$ is $\e_1$ near a plane $P$ parallel to $Q$ in $B(0,1)$, hence $E\cap B(0,1)\subset B(P,\e_1)$ and $d(0,P)\le\e_1$, therefore $d(Q,P)\le \e_1$. Hence $E\cap B(0,1-\e)\subset B(Q,2\e_1)$, such that $E\cap \partial B(0,1-\e))\subset B(Q,2\e)$, and therefore 
\be\mbox{for every }x\in  E\cap \partial B(0,1-\e)), ||\pi(x)||\ge \sqrt{(1-\e)^2-(2\e)^2}\ge\frac 34.\ee

But by Theorem $6.9$, $E\cap B(0,1-\e)$ is a topological disc, hence by a topological argument similar to that of Lemma 4.9, (6.28) gives $\pi(E\cap B(0,1-\e))\supset B(0,\frac 34)\cap Q$. Thus we have (6.27). 
 
 Now let $\Gamma$ be a connected component of $F:=E\cap B(0,1-\e)\cap \pi^{-1}(B(0,\frac34)\cap Q)$. Then it is both open and closed in $F$. But we know that $D\pi (x)$ is injective for all $x\in F$, so $\pi$ is an open map, such that $\pi(\Gamma)$ is open in $B(0,\frac34)\cap Q$. On the other hand, we claim that $\pi(\Gamma)$ is also closed in $B(0,\frac34)\cap Q$. In fact, suppose that $\{x_n\}\subset \pi(\Gamma)$ is a sequence of points that converge to a point $x_\infty\in B(0,\frac34)\cap Q$. For each $n$, take $y_n\in\Gamma$ such that $\pi(y_n)=x_n$. Then 
$\{y_n\}\subset \overline\Gamma$, which is compact, hence the sequence $\{y_n\}$ admits a limit point $y_\infty\subset \overline \Gamma$. Therefore we have $\pi(y_\infty)=x_\infty$. We want to prove that $y_\infty\in\Gamma$, hence we look at $\overline\Gamma\bs\Gamma$. Since $\Gamma$ is closed in $F=E\cap B(0,1-\e)\cap \pi^{-1}(B(0,\frac34)\cap Q)$,
\be \begin{split}&\overline\Gamma\bs\Gamma\subset E\cap \partial [B(0,1-\e)\cap \pi^{-1}(B(0,\frac34)\cap Q)]\\
 &=E\cap\{[\partial B(0,1-\e)\cap \pi^{-1}(\overline B(0,\frac34)\cap Q)]\cup [\partial (\pi^{-1}(B(0,\frac34)\cap Q))\cap \overline B(0,1-\e)]\}.
 \end{split}\ee

 We know that the distance $d(\partial B(0,1-\e)\cap \pi^{-1}(\overline B(0,\frac34)\cap Q),Q)>\sqrt{(1-\e)^2-(\frac34)^2}>\sqrt{\frac{7}{16}-2\e}$, and $d(P,Q)\le d(0,P)<\e_1$, since $0\in Q$ and $0\in E$. As a result,
 \be d(\partial B(0,1-\e)\cap \pi^{-1}(\overline B(0,\frac34)\cap Q),P)>\sqrt{\frac{7}{16}-2\e}-\e_1>\e_1\ee 
 when $\e$ and $\e_1$ are both small. Then the hypothesis says that for each $y\in E\cap B(0,1)$, $d(y,P)<\e_1$, hence $[E\cap B(0,1)]\cap \partial B(0,1-\e)\cap \pi^{-1}(\overline B(0,\frac34)\cap Q)=\emptyset$. As a result, $\overline\Gamma\subset E\cap B(0,1)$ does not meet $\partial B(0,1-\e)\cap \pi^{-1}(\overline B(0,\frac34)\cap Q)$. On combining with (6.29),
\be\overline\Gamma\bs\Gamma\subset \partial (\pi^{-1}(B(0,\frac34)\cap Q))\cap \overline B(0,1-\e).\ee
 
But after the hypothesis, $\pi(y_\infty)=x_\infty\in B(0,\frac34)\cap Q$, hence $y_\infty\not\in \partial (\pi^{-1}(B(0,\frac34)\cap Q))$. Therefore, $y_\infty\not\in\overline\Gamma\bs\Gamma$. Hence $y_\infty\in\Gamma$, and thus $x_\infty\in\pi(\Gamma)$. 

So $\pi(\Gamma)$ is also closed in $B(0,\frac34)\cap Q$. As a result, $\pi(\Gamma)=B(0,\frac34)\cap Q$.
 
Now we claim that
\be\pi:\Gamma\to B(0,\frac 34)\cap Q\mbox{ is a covering space of }Q\cap B(0,\frac 34).\ee
In fact, since $E$ is compact, the continuous map $\pi:E\to Q$ is proper, such that for each $x\in Q\cap B(0,\frac 34)$, $\pi^{-1}(x)\cap \Gamma$ is a finite set. Denote this set $\{y_1,\cdots, y_N\}$. Then by the conclusion before (6.26), for each $1\le j\le N$, there exists $r_j>0$ such that in $\pi^{-1}[B(x,r_j)\cap Q]\cap B(y_j,2r_j)$, $\Gamma$ coincides with the graph of a map $g_j:Q\to Q^\perp$ on $B(x,r_j)\cap Q$. Set $r=\min_jr_j$, then $\pi^{-1}[B(x,r)\cap Q]$ contains the finite disjoint union of these $g_j(B(x,r)\cap Q),1\le j\le N$. On the other hand, for each $y\in \pi^{-1}(B(x,r)\cap Q)$, take a connected component $\gamma$ of $\pi^{-1}(B(x,r)\cap Q)$ such that $y\in\gamma$, then by an argument similar to the one for $\Gamma$ on $B(0,\frac 34)\cap Q$ above, $\pi(\gamma)\supset B(x,r)\cap Q$ , in particular, there exists $1\le j\le N$ such that $y_j\in\gamma$. But in this case we have $\gamma=g_j(B(x,r)\cap Q)$, hence $y\in g_j(B(x,r)\cap Q)$. As a result, $\pi^{-1}[B(x,r)\cap Q]$ is just a finite disjoint union of $g_j(B(x,r)\cap Q), 1\le j\le N$, and on each of theses pieces, $\pi$ is an homeomorphism of $B(x,r)$, where (6.32) follows.

But $Q\cap B(0,\frac 34)$ is simply connected, $\Gamma$ is its connected covering space by $\pi$, hence $\pi$ has to be a homeomorphism. Then by the conclusion around (6.26), $\Gamma$ is the graph of a $C^1$ map from  $Q$ to $Q^\perp$ whose gradient is of $L^\infty$ norm less than 1. In particular, the measure of $\Gamma$ is larger than the measure of $B(0,\frac 34)\cap Q$, which is $\frac{9}{16}\pi$.

Now denote by $\Gamma_1,\cdots,\Gamma_n,\cdots$ all the connected components of $F=E\cap B(0,1-\e)\cap \pi^{-1}(B(0,\frac34)\cap Q)$, the measure of each $\Gamma_i$ is larger than $\frac{9}{16}\pi$. Then if there exists more than one $\Gamma_i$, we have 
\be H^2(E\cap B(0,1))>H^2(E\cap B(0,1-\e)\cap \pi^{-1}(B(0,\frac 34)\cap Q))\ge 2\times \frac{9}{16}\pi=\frac98\pi,\ee
which contradicts (6.15).

Thus $E\cap B(0,1-\e)\cap \pi^{-1}(B(0,\frac 34)\cap Q)$ is a trivial covering space of $Q\cap B(0,\frac 34)$. Combine with (6.27), $\pi^{-1}: B(0,\frac 34)\cap Q\to E\cap B(0,1-\e)\cap \pi^{-1}(B(0,\frac 34))$ is a $C^1$ map, which coincides with $g_x$ around $\pi(x)$ for each $x\in E\cap B(0,1-\e)\cap \pi^{-1}(B(0,\frac 34))$, and hence $||\nabla g||_\infty<1$. Now set $f=g\circ \pi:P\cap \pi^{-1}(B(0,\frac 34)\to E\cap B(0,1-\e)\cap \pi^{-1}(B(0,\frac 34))$. Thus since $P$ is parallel to $Q$, we get the desired conclusion.\qed

\begin{rem}For $d=2,n=4$, we can also get the same result by Theorem 1.15 of \cite{DEpi}, without all those complicated concepts such as varifolds, etc. \end{rem}

\begin{cor}\label{graphe}There exists $\e_1>0$ such that if $k$ is large enough, $E$ is a locally minimal sets in a domain $U\subset\R^4$, $D_k(0,1)\subset U$, and $E$ is $\e_1$ near a plane $P$ in $D_k(0,1)$, then in $D_k(0,\frac12)$, $E$ coincides with the graph of a $C^1$ map $f:P\to P^\perp$. Moreover $||\nabla f||_\infty<1$.\end{cor}

\nd When $k$ is large enough, we have $D_k(0,\frac12)\subset B(0,\frac34)\subset B(0,1)\subset D_k(0,1)$. Then if $E$ is $\e_1$ near a plane $P$ in $D_k(0,1)$, this implies that $E$ is $\e_1$ near $P$ in $B(0,1)$. Thus by Proposition 6.14, in $B(0,\frac34)$, $E$ is the graph of a $C^1$ map $C^1$ $f:P\to P^\perp$, with $||f||_\infty<1$. Therefore in $D_k(0,\frac12)$, too.

\qed

Now fix a large enough $k$, and denote by $D(x,r)=D_k(x,r),C^i(x,r)=C_k^i(x,r)$ for $i=1,2$, and $d_{x,r}=d^k_{x,r}$. 

\noindent \textbf{Proof of (1) of Proposition 6.1.}

Since $k$ is large, $P_k$ is very near $P_0$, there exists $0<\e_3<\frac{1}{100}$ (which does not depend on $k$ for $k$ large), such that for all $\e<\e_3$, if $x\in\R^4$ and $E$ is any set such that $d_{x,r}(E,P_k+q)<\e$, with $q\in\R^4$ and $d(x,q)<20\e r$, then in
$D(x,r)\backslash D(q,\frac {1}{100}r)$, $E$ is the disjoint union of two pieces $E^1,E^2$, such that in $D(x,r)$ minus a small hole, $E^1$ is near $P_k^1+q$, but far from $P_k^2+q$, and vice-versa for $E^2$. More precisely,
\be E^i\subset B((P_k^i+q)\cap D(x,r)\backslash D(q,\frac {1}{100}r),\e r)\ee
and
\be \begin{array}{ll}&d(B(P_k^1\cap D(x,r)\backslash D(q,\frac{1}{100}r),\e r),B(P_k^2\cap D(x,r),\e r))>\frac{1}{80}r;\\
&d(B(P_k^2\cap D(x,r)\backslash D(q,\frac{1}{100}r),\e r),B(P_k^1\cap D(x,r),\e r))>\frac{1}{80}r.\end{array}\ee
In particular,
\be d(E^1,E^2)\ge \frac{1}{80}r.\ee

Take $\e_0=\min\{\frac12\e_3,\frac{1}{40}\e_1,10^{-5}\}$. Then for every $\e<\e_0$, if the $\e$-process does not stop before the step $n$, $E_k$ is $\e s_{n-1}$ near $P_k+q_n$ in $D(q_{n-1},s_{n-1})$. Fix this $n$, and denote by $q=q_n,x=q_{n-1},r=s_{n-1}$ for short. Then since $\e<\e_3$, 
\be E_k\cap D(x,r)\backslash D(q,\frac {1}{100}r)\mbox{ is the disjoint union of }E^1,E^2\mbox{ such that (6.36)-(6.38) hold.}\ee

For each $y\in E^1\cap D(x,r-\frac{1}{40}r)\backslash D(q,\frac {1}{20}r)=E^1\cap D(q_{n-1},\frac{39}{40}s_{n-1})\backslash D(q_n,\frac{1}{10}s_n)$,  
\be d_{y,\frac{1}{40}r}(E_k,P_k+q)<40\e<\e_1,\ee
where by definition,
\be\begin{array}{lll}
d_{y,\frac{1}{40}r}(E_k,P_k+q)=\frac{1}{\frac{1}{40}r}\max\{\sup\{&d&(z,P_k+q):z\in E_k\cap D(y,\frac{1}{40}r)\},\\
&&\sup\{d(z,E_k):z\in P_k+q\cap D(y,\frac{1}{40}r)\}\}.
\end{array}
\ee
For the second term,
\be\sup\{d(z,E_k):z\in (P_k+q)\cap D(y,\frac{1}{40}r)\}\ge \sup\{d(z,E_k):z\in (P_k^1+q)\cap D(y,\frac{1}{40}r)\}\ee 
since $P_k^1\subset P_k$; for the first term, notice first that
% \be E_k\cap D(y,\frac{1}{40}r)=E^1\cap D(y,\frac{1}{40}r),\ee
% puisque $y\in E^1\cap D(x,r-\frac{1}{40}r)\bs D(q,\frac{1}{20}r)$ implique que $D(y,\frac{1}{40}r)\subset D(x,r)\bs D(q,\frac{1}{40}r)$, and \`a cause de (6.13). Ensuite, on sait d\'ej\`a que $d_{y,\frac{1}{40}r}(E_k,P_k+q)<40\e$, qui implique for chaque $z\in E_k\cap D(y,\frac{1}{40}r)$, il existe un $w\in P_k+q$ such that $d(z,w)<\e r$. Par cons\'equent 
% \be w\in B(E^1\cap D(y,\frac{1}{40}r),\e r)\subset B(D(y,\frac{1}{40}r),\e r))\subset D(y,\frac{1}{40}r)+\e r\ee
% 
\be\begin{array}{lll}
\sup\{d(z,P_k+q)&:&z\in E_k\cap D(y,\frac{1}{40}r)\}\\&=&\sup\{d(z,(P_k+q)\cap D(y,\frac{1}{40}r+\e r)):z\in E_k\cap D(y,\frac{1}{40}r)\}   
\end{array}\ee
since we already know that $d_{y,\frac{1}{40}r}(E_k,P_k+q)<40\e$, which implies that for each $z\in E_k\cap D(y,\frac{1}{40}r)$, there exists $w\in P_k+q$ such that $d(z,w)<\e r$. Set $W=\{w\in P_k+q,d(z,w)<\e r\}$. Then
\be d(z,P_k+q)=d(z,W).\ee

For all $w\in W$,
\be\begin{array}{lll} w &\in& (P_k+q)\cap D(y,\frac{1}{40}r+\e r)\subset (P_k+q)\cap D(y,\frac{1}{40}r+\frac{1}{100}r)\\
&\subset& (P_k+q)\cap D(x,r)\backslash D(q,\frac{1}{100}r)\\
&=&[(P_k^1+q) \cap D(x,r)\backslash D(q,\frac{1}{100}r)]\cup (P_k^2+q) \cap D(x,r)\backslash D(q,\frac{1}{100}r). 
   \end{array}\ee

Then $w$ has to belong to $(P_k^1+q) \cap D(x,r)\backslash D(q,\frac{1}{100}r)$, because otherwise 
\be z\in B(w,\e r)\cap E_k\subset B((P_k^2+q) \cap D(x,r)\backslash D(q,\frac{1}{100}r),\e r)\cap E_k=E^2\ee
by (6.36) and (6.37), which contradicts the fact that $z\in E^1$.
Hence
\be d(z,P_k+q)=d(z,W)\ge d(z,P_k^1+q)\mbox{ for }z\in E_k\cap D(y,\frac{1}{40}r),\ee
therefore
\be\sup\{d(z,P_k+q):z\in E_k\cap D(y,\frac{1}{40}r)\}\ge\sup\{d(z,P_k^1+q):z\in E_k\cap D(y,\frac{1}{40}r)\}.\ee
Add (6.42) and (6.48) together we obtain
\be d_{y,\frac{1}{40}r}(E^1,P_k^1+q)\le d_{y,\frac{1}{40}r}(E^1,P_k+q)<40\e<\e_1.\ee

Now $P_k^1+q$ is a plane, hence we can use Corollary \ref{graphe}, which gives 
\be\begin{array}{ll}
\mbox{for each }y\in E^1\cap D(x,\frac{39}{40}r)\backslash D(q,\frac{1}{20}r),\mbox{ in }D(y,\frac{1}{80}r),\\
E_k\mbox{ is the graph of a }C^1\mbox{ map } f_y:P_k^1\to {P_k^1}^\perp\mbox{ with }||\nabla f_y||<1.\end{array}\ee
But $E_k\cap D(y,\frac{1}{80}r)=E^1\cap D(y,\frac{1}{80}r)$, which implies that around every point $y\in E^1\cap D(x,\frac{39}{40}r)\backslash D(q,\frac{1}{20}r)$, $E^1$ is locally a $C^1$ graph on $P_k^1$.

Let us verify that in $D(x,\frac{39}{40}r)\backslash D(q,\frac{1}{20}r)$, $E^1$ coincides with the graph of a $C^1$ map on the whole $P_k^1$, whose gradient is of norm $L^\infty$ less than 1. However we have already our small local graph near every point, with small gradient, so we only have to show that the projection $p_k^1:E^1\to P_k^1\cap C^1(x,\frac{39}{40}r)\backslash C^1(q,\frac{1}{20}r)$ is bijective on $E^1\cap D(x,\frac{39}{40}r)\backslash D(q,\frac{1}{20}r)$.

Surjectivity: Set $A=p_k^1(E^1)\cap C^1(x,\frac{39}{40}r)\backslash C^1(q,\frac{1}{20}r)$. Then $A$ is non empty. We are going to show that $A=P_k^1\cap C^1(x,\frac{39}{40}r)\backslash C^1(q,\frac{1}{20}r)$.

First $A$ is closed in $P_k^1\cap C^1(x,\frac{39}{40}r)\backslash C^1(q,\frac{1}{20}r)$, since $E^1$ is compact in $D(x,\frac{39}{40}r)\backslash D(q,\frac{1}{20}r)$.

But $A$ is also open, because if $z\in A$, then there exists $y\in E^1\cap D(x,\frac{39}{40}r)\backslash D(q,\frac{1}{20}r)$ such that $p_k^1(y)=z$. Thus by (6.50), we know that $B(z,\frac{1}{80}r)\cap P_k^1\subset A$. Hence $A$ is open.

Notice that $P_k^1\cap C^1(x,\frac{39}{40}r)\backslash C^1(q,\frac{1}{20}r)$ is connected, hence $A=P_k^1\cap C^1(x,\frac{39}{40}r)\backslash C^1(q,\frac{1}{20}r)$, which gives the surjectivity.

Injectivity: Suppose $p_k^1$ is not injective. Then there exists $y_1,y_2\in E^1\cap D(x,\frac{39}{40}r)\backslash D(q,\frac{1}{20}r)$ such that $p_k^1(y_1)=p_k^1(y_2).$ In other words 
\be y_1-y_2\in{P_k^1}^\perp.\ee
We know that in $D(y_1,\frac{1}{80}r)$, $E^1$ is a graph, hence $y_2\not\in D(y_1,\frac{1}{80}r)$. In other words, $|y_1-y_2|>\frac{1}{80}r$. Hence there exists at least one point between $y_1,y_2$ whose distance to $P_k^1+q$ is larger than $\frac{1}{160}r>\e r$. This gives a contradiction with (6.48) and the fact that $d(z,P_k+q)\le rd_{x,r}(E_k,P_k+q)<\e r$.
 
Therefore $p_k^i(q)$ is injective. Denote by $f^1$ the map defined on $P_k^1\cap C^1(x,\frac{39r}{40}+\frac{r}{80})\bs C^1(q,\frac{r}{20}-\frac{r}{80})$ and which coincides with $f_y$ on every $B(p_k^1(y),\frac{1}{80}r)\cap P_k^1$; then in $E^1\cap D(x,\frac{39}{40}r)\backslash D(q,\frac{1}{20}r)$, $E^1$ is the graph of $f^1$ with $||\nabla f^1||_\infty<1$.
 
 By a similar argument we obtain also that $E^2\cap D(x,\frac{39}{40}r)\backslash D(q,\frac{1}{20}r)$ is the graph of a $C^1$ map $f^2$, which sends $P_k^2\cap C^2(x,\frac{39}{40}r)\backslash C^2(q,\frac{1}{20}r)$ in ${P_k^2}^\perp$. Recall that $D(x,r)=D(q_{n-1},s_{n-1})$, and by replacing $E^i$ by $E^i(n)$, $f^i$ by $f^i(n)$, we can obtain our graph $E^i(n)=f^i(n)(P_k^ii\cap C^i(q_{n-1},\frac{39}{40}s_{n-1})\backslash C^i(q_n,\frac{1}{10}s_n)$, on condition that the $\e$ process does not stop at the step $n$. But of course if it does not stop at step $n$, it does not stop at any step before $n$ neither.  Therefore for all $j\le n$, we have decompositions of $E_k\cap  D(q_{j-1},\frac{39}{40}s_{j-1})\backslash D(q_j,\frac{1}{10}s_j)$ as disjoint unions
\be E_k\cap  D(q_{j-1},\frac{39}{40}s_{j-1})\backslash D(q_j,\frac{1}{10}s_j)=E^1(j)\cup E^2(j),\ee and 
\be E^i(j)\mbox{ is the graph of }g^i(j)\mbox{ on }P_k^i\cap C^i(q_{j-1},\frac{39}{40}s_{j-1})\backslash C^i(q_j,\frac{1}{10}s_j).\ee

We can easily verify that if $j,l$ are such that $x\in P_k^i\cap[C^i(q_{j-1},\frac{39}{40}s_{j-1})\backslash C^i(q_j,\frac{1}{10}s_j)]\cap [C^i(q_{l-1},\frac{39}{40}s_{l-1})\backslash C^i(q_l,\frac{1}{10}s_l)]$, then $g^i(j)(x)=g^i(l)(x)\in E^i(j)\cap E^i(l)$. Hence set
\be \begin{split}g^i\ &:\ P_k^i\cap D(0,\frac{39}{40})\backslash D(p_k^i(q_n),\frac{1}{10}s_n)\to {P_k^i}^\perp\ ;\\
g^i(x)&=g^i(j)(x)\mbox{ on }P_k^i\cap D(p_k^i(q_{j-1}),\frac{39}{40}s_{j-1})\backslash D(p_k^i(q_j),\frac{1}{10}s_j),1\le j\le n\ ;\end{split}\ee
then $||\nabla g^i||_\infty<1$, and its graph is $G^i=[\cup_{j=0}^n E^i(j)]\cap D(0,\frac{39}{40})\backslash D(q_n,\frac{1}{10}s_n)$.

Thus all we have to do is to show that $G^1,G^2$ are disjoint. This is equivalent to saying that for $0\le j,l\le n$, $E^1(j)\cap E^2(l)=\emptyset$. This is true for $j=l$, so suppose that $j<l$. Then for all point $x\in E^1(j)$ there are two cases: either $x\in D(q_{l-1},\frac{39}{40}s_{l-1})\backslash D(q_l,\frac{1}{10}s_l)$, either not.  In the second case, $x\not\in E^2(l)$ automatically because $E^2(l)\subset D(q_{l-1},\frac{39}{40}s_{l-1})\backslash D(q_l,\frac{1}{10}s_l)$; in the first case, $x\in E_k\cap D(q_{l-1},\frac{39}{40}s_{l-1})\backslash D(q_l,\frac{1}{10}s_l)=E^1(l)\cup E^2(l)$. Then $D(q_{l-1},s_{l-1})\backslash D(q_{j-1},\frac{1}{10}s_{j-1})\ne\emptyset$ implies that $l-j\le 4$.

But $x\in E^1(j)$ implies that $d(x,P_k^2+q_j)>\frac{1}{80}s_{j-1}$ because of (6.37). Hence  
\be d(x,P_k^2+q_l)>\frac{1}{80}s_{j-1}-d(q_j,q_l).\ee
While by (5.16) we have $d(q_j,q_l)\le 24\e \times 2^{-\min(l,j)}=24\e_0\times 2^{-j}\le \frac{48}{10^5}s_{j-1}$, hence
\be d(x,P_k^2+q_l)>\frac{3}{400}s_{j-1}\ge 2^4 \frac{3}{400}s_{l-1}>\e s_{l-1},\ee
therefore $x\not\in E^2(l)$, because of (6.36). 

Hence for all $0\le j,l\le n$, $E^1(j)\cap E^2(l)=\emptyset$, therefore $G^1\cap G^2=\emptyset.$ Thus we complete the proof of (1).\qed

\medskip

\noindent\textbf{Proof of (2) of Proposition 6.1.} Since $E_k$ is $\e$ near $P_k$ in $D(0,1)$, we can see that $E_k\cap D(0,1)\bs D(0,\frac{1}{100})$ is the disjoint union of two pieces $E^1,E^2$ which satisfy (6.36)-(6.38). By (1), we know that in $D(0,\frac{39}{40})\backslash D(q_n,\frac{1}{10}s_n)$ $E_k$ is composed  of two disjoint graphs $G^1,G^2$ on $P_k^1\cap C^1(0,\frac{39}{40})\backslash C^1(q_n,\frac{1}{10}s_n)$ and $P_k^2\cap C^2(0,\frac{39}{40})\backslash C^2(q_n,\frac{1}{10}s_n)$ respectively. Then for $\frac{1}{10}s_n<t<s_n$, $G_t^i=E^i\cup G^i\backslash D(q_n,t)$, hence (6.4) is true. Moreover, since $G^i,i=1,2$ are graphs, we have 
\be p_k^i(G^i\backslash D(q_n,t))\supset P_k^i\cap D(0,\frac{39}{40})\backslash C_k^i(q_n,t),\ee
therefore (6.5) is also true if we replace $D(0,1)$ by $D(0,\frac{39}{40})$, because $G^i\backslash D(q_n,t)\subset G_t^i.$

So we just have to prove
\be P_k^i\cap D(0,1)\backslash D(0,\frac{39}{40})\subset p_k^i(G_t^i)\ee

We prove it for $i=1$ for example. We know that $E_k$ is $\e$ near $P_k$ in $D(0,1)$. Hence by (6.36), 
\be E^2\subset B(P_k^2\cap D(0,1)\backslash D(0,\frac{1}{100}),\e)\subset B(P_k^2,\frac{1}{100})\ee
However when $k$ is large, we have $p_k^1(P_k^2\cap D(0,1))\subset P_k^1\cap D(0,\frac 15)$, and therefore
\be p_k^1(E^2)\subset P_k^1\cap D(0,\frac 14).\ee
On the other hand, we know that
\be p_k^1(E_k\cap D(0,\frac{1}{100}))\subset P_k^1\cap D(0,\frac{1}{100})\ee
hence we have \be p_k^1(E_k\backslash E^1)=p_k^1[E^2\cup (E_k\cap D(0,\frac{1}{100}))]\subset D(0,\frac 14)\cap P_k^1.\ee
By Proposition 4.8, $p_k^1(E_k)\supset P_k^1\cap D(0,1)$, so we get that
\be p_k^1(E^1)\supset P_k^1\cap D(0,1)\backslash D(0,\frac 14).\ee
Therefore
\be\begin{array}{ll} &p_k^1(G_t^1)=p_k^1(G^1\backslash D(q_n,t))\cup p_k^1(E^1)\\
&\supset [P_k^1\cap D(0,\frac{39}{40})\backslash D(q_n,t)]\cup [P_k^1\cap D(0,1)\backslash D(0,\frac 14)]=p_k^1\cap D(0,1)\backslash D(q_n,t),
\end{array}\ee
where (6.5) follows for $i=1$. The proof for $i=2$ is similar.\qed

Now let us deal with (3) and (4). First of all we give a small remark. In fact all we want is 4), i.e., the surjectivity of the projections, to estimate the measure of $E_k$ for the part where we do not know much about its structure. Notice that by the proof of Theorem 4.1, we know that $E_k$ is the limit of a sequence $\{H_l\}$ of deformations of $P_k$ in $U=\R^4\bs[P_k\bs B(0,1)].$ Hence in $E$ the projection $p_k^i$ of $E_k$ is automatically surjective on $D(0,1)\cap P_k^i$. Next we look into $D(0,\frac 12)$, since $E_k$ is near $P_k$, the part of $E_k\bs D(0,\frac12)$ which is near $P_k^1$ has no projection in $D(0,\frac12)\cap P_k^1$, and the part of $E_k\bs D(0,\frac12)$ which is near $P_k^2$ has a very small projection on $P_k^1$, and this projection is very near the origin, hence we can say that outside a small ball, say $D(0,\frac{1}{100})$, the projection $p_k^1$ from $E_k\cap D(0,\frac12)$ to $P_k^1$ is surjective on $P_k^1\cap D(0,\frac12)\bs D(0,\frac{1}{100})$. However in $D(0,\frac{1}{100})$, we can not say directly that the projection of $E_k$ comes from $E_k\cap D(0,\frac 12)$, because part of it may come from the part near $P_k^2$. (Picture 6-1 below may give an idea).

\noindent\centerline{\includegraphics[width=0.4\textwidth]{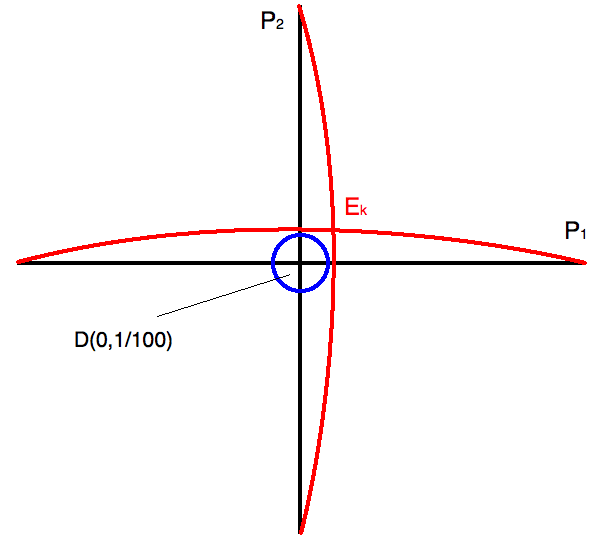}}
\nopagebreak[4]
\centerline{6-1}

Therefore the idea to prove 4) is still to prove that $E_k\cap \overline D(0,\frac 12)$ is the limit of a sequence of deformations of $P_k\cap \overline D(0,\frac 12)$. In other words, we can contract the part outside $\overline D(0,\frac12)$ into $\overline D(0,\frac12)$. This part is very flat and regular, hence intuitively such a contraction will not change essentially the structure of $E_k$ in $D(0,\frac 12)$.

Here we still have to say, if $E_k$ is itself a deformation of $P_k$ and is very near $P_k$, then we can easily contract the part $E_k\bs D(0,\frac 12)$ to $E_k\cap\partial D(0,\frac 12)$, just like we contract an annulus to the interior circle, because $E_k\bs D(0,\frac 12)$ is roughly composed of two pieces of $C^1$ graph on $P_k^i\bs D(0,\frac 12)$, thanks to (6.2). Then next we can carry on the same operation in $D(0,\frac12)$, if $E_k$ is still near some translation of $P_k$ in $D(0,\frac12)$. At last, we will arrive at the scale where the $\e-$process stops. Hence we can say that $E_k\cap D(q_n,t)$ is a deformation of $P_k$, too, because it is a deformation of $E_k$.

However $E_k$ may not be a deformation of $P_k$. Hence we will use the fact that $E_k$ is the limit of a sequence of deformations $\{H_l\}$, and we want to apply the argument above to prove that $H_l\cap \overline D(0,\frac 12)$ is a deformation of $P_k\cap \overline D(0,\frac 12)$. But this time, $H_l$ is not minimal, hence we cannot use (6.2) to say that $H_l\cap\partial D(0,\frac 12)$ is a very regular curve, and therefore it is not that easy to manage to contract $H_l$ directly on $H_l\cap\partial D(0,\frac12)$.

Then what we are going to do is, first use the shortest distance projection $\pi$ to project $H_l$ on $\overline D(0,\frac 12)$. Then the points of $H_l\bs D(0,\frac 12)$ are sent by $\pi$ to $\partial D(0,\frac 12)$ (but the image is not $H_l\cap\partial D(0,\frac12)$ anymore). To continue to project $\pi(H_l)\bs D(q_2,\frac14)$ in $\overline D(q_2,\frac 14)$, we use the fact that $\pi(H_l)\bs D(q_2,\frac14)$ is composed of two disjoint pieces near $P_k^1+q_2$ and $P_k^2+q_2$ respectively. To guarantee this, by the argument before (6.36), $\pi(H_l)$ should be $\e_3 s_2=\frac12\e_3$ near $P_k+q_2$. For the part inside $D(0,\frac12)$, there is no problem, because the $\e-$process does not stop here. However for the part on the boundary $\partial D(0,\frac12)$, things are complicated because this part contains the projection of the $H_l\bs D(0,\frac 12)$, which is $2\e_3 s_2$ near $p_k+q_2$, not $\e_3 s_2$ near it. But luckily we can manage first to do a contraction in $\partial D(0,\frac12)$, to contract everything near $P_k+q_2$, without moving points inside $D(0,\frac12)$. Thus we make things more complicated, but at last we can still carry on with the proof.

After all these, we can prove the surjectivity of projections of $E_k\cap D(q_n,t)$, since it is the limit of a sequence of deformations.

Now we begin to concretize the above idea.

\noindent\textbf{Proof of (3) and (4) of Proposition 6.1.} Fix a $\frac{1}{10}s_n\le t\le s_n$. 

We know that for all $j\le n$, in every $D(q_{j-1},s_{j-1})$, $E_k$ is $\e s_{j-1}$ near $P_k+q_j$, because the $\e-$process does not stop on step $j$. Then set $E_k$ is the limit of $H_l$, hence we can suppose that $l$ is large, so that in every $D(q_{j-1},s_{j-1})$, $H_l$ is $2\e s_{j-1}$ near $P_k+q_j$, for all $j\le n$. By definition of $\e_0$, we have $2\e<2\e_0<\e_3$, hence (6.36) and (6.37) give that in $D(q_{j-1},s_{j-1})\bs D(q_j,\frac{1}{100}s_{j-1}),\ H_l$ is the union of two disjoint pieces $H_l^1,H_l^2$ with
\be H_l^i\subset B(P_k^i+q_j\cap D(q_{j-1},s_{j-1})\backslash D(q_j,\frac{1}{100}s_{j-1}),2\e r).\ee

We will construct our deformation $F_l^n$ by recurrence on $j<n$.

For $j=1$, define $\pi_1 : H_l\to D(q_1,s_1)$ the shortest distance projection from $\R^4$ to $D(q_1,s_1)$. Notice that although $H_l\cap D(q_1,s_1)$ is $2\e s_1$ near $P_k^1+q_2$ in $D(q_1,s_1)$, $\pi_1(H_l\bs D(q_1,s_1))$ is not necessarily $2\e s_1$ near $P_k^1+q_2$ in $D(q_1,s_1)$. Therefore we will modify it a little, to be able to continue decompose it into two disjoint pieces that verify some conditions similar to (6.36)-(6.38).

By (6.36), in $D(0,1)\bs D(q_1,s_1)$, $H_l$ is the union of two disjoint pieces $H_l^1,H_l^2$, where $H_l^i$ is very near $P_k^i$, so we have
\be \pi_1(H_l^i)\cap D(0,1)\bs D(q_1,s_1))\subset\partial C^i(q_1,s_1)\cap B(P_k^i,2\e).\ee

As a result, the image by $\pi_1$ of each $H_l^i$ outside $D(q_1,s_1)$ is contained in the cylinder $C^i(q_1,s_1)$, and also contained in a very small neighborhood of the plane $P_k^i$. I.e., it is contained in a 3-dimensional thin band around $P_k^i\cap\partial C^i(q_1,s_1)$. In particular, it is far from the boundary $\partial C^j(q_1,s_1)$ for $i\ne j$. Therefore we can define $g_1$ on $\pi_1(H_l)$ by
\be g_1(x)=\left\{\begin{array}{ll}x\ ;\ x\in D(q_1,s_1)\ ;\\
                g_1^i(x)\ ;\ x\in \partial C^i(q_1,s_1).
               \end{array}\right.\ee
where $g_1^i$ is the orthogonal projection on $B(P_k^i+q_2,2\e s_1)$.

It is clear that for all points $x\in H_l\cap D(q_1,s_1)\subset \pi_1(H_l)$, neither $\pi_1$ nor $g_1^i$ can move it, since in $D(q_1,s_1)$ $H_l$ is $2\e s_1$ near $P_k^1+q_2$. Hence the action of $g_1^i$ is just to press all points on the boundary $\partial C_k^i(q_1,s_1)$ into a $2s_1\e$ neighborhood of $P_k^i$, without leaving the boundary. Then $g_1$ is 2-Lipschitz (locally 1-Lipschitz).
% In fact, for any two points $x,y\in D(q_1,s_1)$, $g_1$ is the identity, hence $|g_1(x)-g_1(y)|=|x-y|$ ; for $x,y\in \partial D(q_1,s_1)$, if $x,y$ belongs au m\^eme $\partial C^i(q_1,s_1)$, $g_1$ is juste une projection orthogonale, donc $|g_1(x)-g_1(y)|\le|x-y|$ ; si $x\in \partial C^1(q_1,s_1),y\in \partial C^2(q_1,s_1)$, alors $|x-y|>2\e$, and par d\'efinition de $g_1$, $|g_1(x)-x|$ and $|g_1(y)-y|$ sont inf\'erieurs \`a $\e s_1=\frac 12\e$, de sort que $|g_1(x)-g_1(y)|\le |g_1(x)-x|+|x-y|+|g_1(y)-y|\le \e+|x-y|\le\frac12|x-y|+|x-y|\le 2|x-y|$.
Set $h_1=g_1\circ \pi_1$, then $h_1$ is 2-Lipschitz, and moreover $h_1(H_l)\subset \overline D(q_1,s_1)$ is $2\e s_1$ near $P_k$. See Picture 6-2 and 6-3 below. 6-2 is the set, $H_l$, and in 6-3 we give the image of $H_l$ after $h_1$. 

\noindent\centerline{\includegraphics[width=0.8\textwidth]{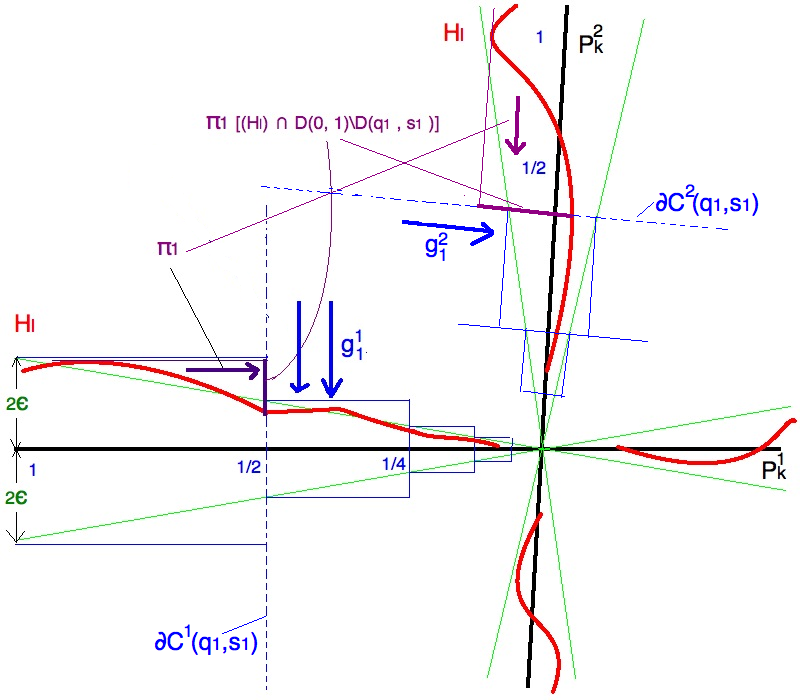}}
\nopagebreak[4]
\centerline{6-2}
\noindent\centerline{\includegraphics[width=0.8\textwidth]{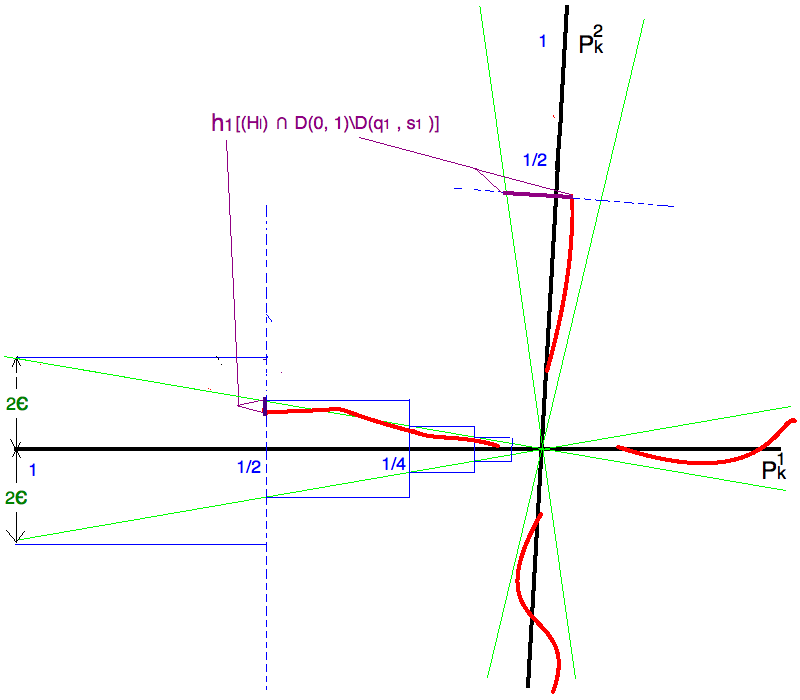}}
\nopagebreak[4]
\centerline{6-3}

Now we have defined $h_1$ in $D(0,1)=D(q_0,s_0)$, which deforms $H_l$ into $D(q_1,s_1)$ and keeps the image $2\e s_1$ near $P_k+q_2$. But $h_1$ is 2-Lipschitz and is defined on a compact set $H_l$, hence we can extend it on the whole $\R^4$. Still call the obtained function $h_1$. Then it is still 2-Lipschitz. 

Now suppose that for $j$, we have constructed a 2-Lipschitz deformation $h_j$ which deform $H_l$ into $D(q_j,s_j)$ and whose image is $2\e s_j$ near $P_k+q_{j+1}$ in $D(q_j,s_j)$. If $j<n-1$, then we can define $\pi_{j+1}$ the shortest distance projection to $D(q_{j+1},s_{j+1})$, and then similarly we project all point on the boundary of $C^i(q_{j+1},s_{j+1})$ into the $2\e s_{j+1}$ neighborhood of $P_k^i+q_{j+2}$. 
%We can defineOn peut le d\'efinir sans ambigu\"it\'e car $h_j(H_l)$ is $2\e s_j$ proche de $P_k+q_{j+1}$, de sort que $h_j(H_l)\bs D(q_{j+1},s_{j+1})$ se d\'ecompose en 2 parties ``plates'' ((6.36)-(6.38)).
Denote this projection by $g_{j+1}$, and set $h_{j+1}=g_{j+1}\circ \pi_{j+1}\circ h_j$. Then $h_{j+1}$ is a 2-Lipschitz deformation that maps $H_l$ in $D(q_{j+1},s_{j+1})$ and its image is $\e s_{j+1}$ near $P_k+q_{j+2}$ in $D(q_{j+1},s_{j+1})$.Thus we obtain our $h_{j+1}$. 

Therefore by recursion, we can define the $h_j$ until $h_{n-1}$. 

Now for each $l$, we will define a deformation $h_n(l)$, which will deform $H_l$ into $D(q_n,t+\frac 1l)$. Denote by $p_n(l)$ the shortest distance projection on $D(q_n,t+\frac 1l)$, and for $x\in \overline C^i(q_n,t+\frac1l)\bs C^i(q_n,t)$, set $g^i_n(l)(x)=(id,g^i)\circ p_k^i\circ p_n(l)$, where $g^i$ is as in (6.2). Then we define the deformation $h_n(l)$ of $P_k\cap D(0,1)$ as follows. For $x\in P_k\cap D(0,1)$:
\be h_n(l)(x)=\left\{\begin{array}{ll} h_{n-1}(x),&h_{n-1}(x)\in D(q_n,t)\ ;\\
                    sg^i_n(l)\circ h_{n-1}(x)+(1-s)h_{n-1}(x),&h_{n-1}(x)\in \partial C^i(q_n,t+\frac1ls).
                   \end{array}\right.
\ee

Then we can see that $h_n(l)(P_k)\cap \overline D(q_n,t)=H_l\cap \overline D(q_n,t)$, and $h_n(l)(P_k)\cap\partial D(q_n,t+\frac 1l)=E_k\cap\partial D(q_n,t+\frac 1l)$. Between $\overline D(q_n,t)$ and $\partial D(q_n,t+\frac 1l)$ the image of $h_n(l)$ is the image of a homotopy between $H_l$ and $E_k$. Then since $H_l$ converges to $E_k$, $F_l^n(t)=h_n(l)(P_k)$ converges to $E_k\cap (q_n,t)$ when $l$ goes to infinity.

Denote by $a_l(t)$ the affine deformation which sends $P_k\cap B(0,1)$ to $(P_k+q_n)\cap D(q_n,t+\frac1l)$, and set 
\be f_l^n(t)=h_n(l)\circ a_l(t).\ee

Then $F_l^n(t)=f_l^n(t)((P_k+q_n)\cap\overline D(q_n,t+\frac 1l))$ verifies all conditions in (3).

The condition (4) is an immediate corollary of (3). In fact we know that $F_l^n(t)$ is a deformation of $(P_k+q_n) \cap D(q_n,t+\frac1l)$ which sends $\partial C^i(q_n,t)$ in $\partial C^i(q_n,t)$, hence
\be p_k^i(F_l^n(t))\supset P_k^i\cap C_k^i(q_n,t+\frac 1l).\ee

Therefore since $E_k\cap D(q_n,t)$ is the limit of $F_l^n(t)$, each projection $p_k^i$ is surjective from $E_k\cap D(q_n,t)$ to $P_k^i\cap C_k^i(q_n,t)$.\qed

\section{Argument of harmonic extension}

In this section we will give some fundamental estimates on the measure of the graph of a $C^1$ function on an almost concentric annulus. One can easily skip this section and admit these estimates, and continue the proof of Theorem \ref{main} from the beginning of the next section.

\begin{pro}Suppose $0<r_0<\frac 12$ and $u_0\in C^1(\partial B(0,r_0)\cap \R^2,\R)$. Denote by $m(u_0)=\frac{1}{2\pi r_0}\int_{\partial B(0,r_0)}u_0$ its average. 

Then for all $u\in C^1((\overline{B(0,1)}\backslash B(0,r_0))\cap \R^2, \R)$ that satisfies
\be u|_{\partial B(0,r_0)}=u_0\ee
we have
\be\int_{B(0,1)\backslash B(0,r_0)}|\nabla u|^2\ge \frac14 r_0^{-1}\int_{\partial B(0,r_0)}|u_0-m(u_0)|^2.\ee
\end{pro}

 \nd

Let $u$ be a $C^1$ function as in the statement of Proposition 7.1. We define $\tilde u\in C((\overline{B(0,\frac{1}{r_0})}\backslash B(0,r_0))\cap \R^2, \R)$, which is also $C^1$ except for on $\partial B(0,1)$, by
\be\tilde u(x)=\left\{\begin{array}{ccc}&u(x),\  & x\in\overline{B(0,1)}\backslash B(0,r_0);\\
&u(\frac{x}{|x|^2}),\ & x\in \overline{B(0,\frac{1}{r_0})}\backslash B(0,1).
 \end{array}\right.\ee 
 
Then we have
\be \int_{\overline{B(0,\frac{1}{r_0})}\backslash B(0,r_0)}|\nabla\tilde u|^2=2\int_{\overline{B(0,1)}\backslash B(0,r_0)}|\nabla u|^2,\ee 
because $t\mapsto \frac{x}{|x|^2}$ is a conformal map, which does not change Dirichlet's energy.
 
 Now since $\tilde u$ satisfies the boundary condition
 \be\tilde u(x)|_{\partial B(0,r_0)}=u_0(x),\ \tilde u(x)|_{\partial B(0,\frac{1}{r_0})}=u_0(\frac{x}{|x|^2}),\ee
its Dirichlet's energy $\int_{\overline{B(0,\frac{1}{r_0})}\backslash B(0,r_0)}|\nabla\tilde u|^2$ is no less than that of the harmonic function $v$ satisfying the same boundary condition. So we are going to calculate $\int_{\overline{B(0,\frac{1}{r_0})}\backslash B(0,r_0)}|\nabla v|^2$.
 
We write
\be u(r_0,\theta)=u_0(\theta)=m(u_0)+\sum_{n=1}^\infty (A_n\cos n\theta+B_n\sin n\theta),\ee
and set
\be\begin{split}
 v& :\ \overline{B(0,\frac {1}{r_0})}\backslash B(0,r_0)\to\R,\\
  v(r,\theta)&=m(u_0)+\sum_{n=1}^\infty a_n(r^n+r^{-n})\cos n\theta+\sum_{n=1}^\infty b_n(r^n+r^{-n})\sin n\theta.\end{split}\ee 
Then $v$ is harmonic.

The boundary condition
\be v(x)|_{\partial B(0,r_0)}=u_0(x),\ v(x)|_{\partial B(0,\frac{1}{r_0})}=u_0(\frac{x}{|x|^2})\ee 
gives that 
\be a_n(r_0^n+r_0^{-n})=A_n,\ b_n(r_0^n+r_0^{-n})=B_n.\ee 

For estimate $\nabla v$, write
\be\begin{split}v(r,\theta)&=m(u_0)+\sum_1^\infty a_n(r^n+r^{-n})\cos n\theta+\sum_1^\infty b_n(r^n+r^{-n})\sin n\theta\\
&=m(u_0)+\sum_1^\infty\frac{r^n+r^{-n}}{r_0^n+r_0^{-n}}(A_n\cos n\theta+B_n\sin n\theta)\\
&=: m(u_0)+\sum_1^\infty v_n(r,\theta).
\end{split}\ee
It is not hard to verify that we can differentiate $v$ term by term. Hence we have

\be \frac{\partial v}{\partial \theta}=\sum_{n=1}^\infty n(r^n+r^{-n})(-a_n\sin n\theta+b_n\cos n\theta)\ee 
and
\be \frac{\partial v}{\partial r}=\sum_{n=1}^\infty n(r^{n-1}-r^{-n-1})(a_n\cos n\theta+b_n\sin n\theta).\ee 
Therefore
\be\begin{split}
|\nabla v|^2&=|\frac{\partial v}{\partial r}|^2+|\frac 1r\frac{\partial v}{\partial \theta}|^2\\
%&=[\sum_n n^2 (r^{n-1}-r^{-n-1})^2(a_n\cos n\theta+b_n\sin n\theta)^2\\
%&+2\sum_{n<m} nm (r^{n-1}-r^{-n-1})(r^{m-1}-r^{-m-1})(a_n\cos n\theta+b_n\sin n\theta)(a_m\cos m\theta+b_m\sin m\theta)]\\
%&+\frac{1}{r^2}[\sum_n n^2(r^n+r^{-n})^2(-a_n\sin n\theta+b_n\cos n\theta)^2\\
%&+2\sum_{n<m}nm(r^n+r^{-n})(r^m+r^{-m})(-a_n\sin n\theta+b_n\cos n\theta)(-a_m\sin m\theta+b_m\cos m\theta)]\\
%&=\sum_n n^2 \{(r^{2n-2}+r^{-2n-2})[(a_n\cos n\theta+b_n\sin n\theta)^2+(-a_n\sin n\theta+b_n\cos n\theta)^2]\\
%&-\frac{2}{r^2}[(a_n\cos n\theta+b_n\sin n\theta)^2-(-a_n\sin n\theta+b_n\cos n\theta)^2]\}\\
%&+2\sum_{n<m}nm \{(r^{n-1}r^{m-1}+r^{-n-1}r^{-m-1})[(a_n\cos n\theta+b_n\sin n\theta)(a_m\cos m\theta+b_m\sin m\theta)\\
%&+(-a_n\sin n\theta+b_n\cos n\theta)(-a_m\sin m\theta+b_m\cos m\theta)]\\
%&-(r^{n-1}r^{-m-1}+r^{-n-1}r^{m-1})[(a_n\cos n\theta+b_n\sin n\theta)(a_m\cos m\theta+b_m\sin m\theta)\\
%&-(-a_n\sin n\theta+b_n\cos n\theta)(-a_m\sin m\theta+b_m\cos m\theta)]\}\\
&=\sum_n n^2 (r^{2n-2}+r^{-2n-2})(a_n^2+b_n^2)-\frac{2}{r^2}\sum_n n^2(a_n^2\cos 2n\theta-b_n^2\cos 2n\theta+2a_nb_n\sin 2n\theta)\\
&+2\sum_{n<m}nm (r^{n-1}r^{m-1}+r^{-n-1}r^{-m-1})\\
&\{a_na_m\cos(n-m)\theta+a_nb_m\sin(m-n)\theta+a_mb_n\sin(n-m)\theta-b_nb_m\cos(n-m)\theta\}\\
&-2\sum_{n<m}nm (r^{n-1}r^{-m-1}+r^{-n-1}r^{m-1})\\
&\{a_na_m\cos(n+m)\theta-b_nb_m\cos(n+m)\theta+a_nb_m\sin(n+m)\theta+b_na_m\sin(n+m)\theta\}.
\end{split}\ee
Note that for $n\ne m$ and  $m,n\ge 1$,
\be\begin{split}
\int_0^{2\pi}&\cos(n-m)\theta d\theta=\int_0^{2\pi}\sin(n-m)\theta d\theta=\int_0^{2\pi}\cos(n+m)\theta d\theta\\
&=\int_0^{2\pi}\sin(n+m)\theta d\theta=\int_0^{2\pi}\cos 2n\theta d\theta=\int_0^{2\pi}\sin 2n\theta d\theta=0,
\end{split}\ee
therefore
\be\begin{split}
\int_0^{2\pi}|\nabla v|^2d\theta=2\pi\sum_nn^2(r^{2n-2}+r^{-2n-2})(a_n^2+b_n^2),
\end{split}\ee
hence
\be\begin{split}
\int_{B(0,\frac{1}{r_0})\backslash B(0,r_0)}|\nabla v|^2&=\int_{r_0}^\frac{1}{r_0} rdr\int_0^{2\pi}|\nabla v|^2d\theta\\
&=\int_{r_0}^\frac{1}{r_0} rdr \cdot 2\pi\sum_nn^2(r^{2n-2}+r^{-2n-2})(a_n^2+b_n^2)\\
%&=2\pi\sum_nn^2(a_n^2+b_n^2)\int_{r_0}^\frac{1}{r_0} (r^{2n-1}+r^{-2n-1})dr\\
%&=2\pi\sum_nn^2(a_n^2+b_n^2)(\frac{r^{2n}}{2n}+\frac{r^{-2n}}{-2n})|_{r_0}^\frac{1}{r_0}\\
%&=2\pi\sum_nn^2(a_n^2+b_n^2)\frac{1}{n}(r_0^{-2n}-r_0^{2n})\\
%&=2\pi\sum_nn(a_n^2+b_n^2)(r_0^{-2n}-r_0^{2n})\\
&=2\pi\sum_nn(a_n^2+b_n^2)(r_0^{-n}-r_0^n)(r_0^n+r_0^{-n}).
\end{split}\ee
But for $r_0<\frac12$ and $n\ge 1$, we have 
\be r_0^{-n}-r_0^n\ge \frac12(r_0^n+r_0^{-n})\ee 
and hence
\be\begin{split}
\int_{B(0,\frac{1}{r_0})\backslash B(0,r_0)}|\nabla v|^2
&\ge2\pi\sum_nn(a_n^2+b_n^2)\frac12(r_0^n+r_0^{-n})^2\\
&=\pi\sum_nn(A_n^2+B_n^2)\ge\pi\sum_n(A_n^2+B_n^2).
\end{split}\ee
But
\be \pi\sum_{n=1}^\infty (A_n^2+B_n^2)=\frac12\int_0^{2\pi}|u_0(\theta)-m(u_0)|^2 d\theta
= \frac12r_0^{-1}\int_{\partial B(0,r_0)} |u_0(s)-m(u_0)|^2ds.\ee 
Therefore
\be \int_{B(0,\frac{1}{r_0})\backslash B(0,r_0)}|\nabla v|^2\ge \frac12 r_0^{-1}\int_{\partial B(0,r_0)} |u_0(s)-m(u_0)|^2ds.\ee

Now return to the function $u$. We have 
\be\begin{split} &\int_{B(0,1)\backslash B(0,r_0)}|\nabla u|^2=\frac 12\int_{B(0,\frac{1}{r_0})\backslash B(0,r_0)}|\nabla\tilde u|^2\\
&\ge\frac12\int_{B(0,\frac{1}{r_0})\backslash B(0,r_0)}|\nabla v|^2\ge \frac14 r_0^{-1}\int_{\partial B(0,r_0)} |u_0(s)-m(u_0)|^2ds,\end{split}\ee
where the conclusion follows. 
$\hfill\Box$

\begin{cor}\label{cas1}Let $r_0>0$, $q\in\R^2$ be such that $r_0<\frac12d(q,\partial B(0,1))$, suppose $u_0\in C^1(\partial B(q,r_0)\cap \R^2,\R)$, and denote by $m(u_0)=\frac{1}{2\pi r_0}\int_{\partial B(q,r_0)}u_0$ its average. 

Then for all $u\in C^1((\overline{B(0,1)}\backslash B(q,r_0))\cap \R^2, \R)$ that satisfies
\be  u|_{\partial B(q,r_0)}=u_0\ee 
we have \be \int_{B(0,1)\backslash B(q,r_0)}|\nabla u|^2\ge \frac14r_0^{-1}\int_{\partial B(q,r_0)}|u_0-m(u_0)|^2.\ee 
\end{cor}

\nd

Set $R=d(q,\partial B(0,1))<1$, then $r_0<\frac12R$ and $B(q,R)\subset B(0,1)$. Thus we can apply Proposition 7.1 to $B(q,R)\backslash B(q,r_0)$ and by substituting $y=\frac{x-q}{R}$,
\be\begin{split}
\int_{B(q,R)\backslash B(q,r_0)}&|\nabla_x u|^2 dx=\int_{B(0,1)\backslash B(0,\frac{r_0}{R})}|\frac{1}{R}\nabla_y u|^2R^2dy\\
&=\int_{B(0,1)\backslash B(0,\frac{r_0}{R})}|\nabla_y u|^2dy\ge\frac14(\frac{R}{r_0})\int_{\partial B(0,\frac{r_0}{R})}|u-m(u_0)|^2dy
\end{split}\ee
since $\frac{r_0}{R}<\frac12$.
But
\be\int_{\partial B(0,\frac{r_0}{R})}|u-m(u_o)|^2dy=\int_{\partial B(q,r_0)}|u-m(u_0)|^2(\frac{1}{R})dx,\ee
hence
\be \int_{B(q,R)\backslash B(q,r_0)}|\nabla_x u|^2 dx\ge \frac14r_0^{-1}\int_{\partial B(q,r_0)}|u-m(u_0)|^2dx.\ee 
Now since $B(q,R)\subset B(0,1)$, we have 
\be \int_{B(0,1)\backslash B(q,r_0)}|\nabla u|^2\ge \frac14r_0^{-1}\int_{\partial B(q,r_0)}|u-m(u_0)|^2dx.\ee 
$\hfill\Box$

\begin{lem}
Let $0<r_0<1,$ and $u\in C^1(\ B(0,1)\backslash B(0,r_0),\R)$ be such that $u|_{\partial B(0,r_0)}=\delta r_0$ and $u|_{\partial B(0,1)}=0$; then we have 
\be \int_{B(0,1)\backslash B(0,r_0)} |\nabla u|^2\ge \frac{2\pi\delta^2r_0^2}{|\log r_0|}.\ee 
\end{lem}

\nd

Set $f(r,\theta)=A\log r$ with $A=\frac{\delta r_0}{\log r_0}$. Then $f$ is the harmonic extension with the given boundary value, and
\be \frac{\partial f}{\partial r}=\frac Ar,\ \frac{\partial f}{\partial\theta}=0.\ee 
Hence
\be |\nabla f|^2=|\frac{\partial f}{\partial r}|^2+|\frac1r\frac{\partial f}{\partial\theta}|^2=\frac{A^2}{r^2}.\ee
As a result
\be\begin{split}
\int_{B(0,1)\backslash B(0,r_0)} |\nabla f|^2&=\int_0^{2\pi}d\theta\int_{r_0}^1rdr|\nabla f|^2=2\pi \int_{r_0}^1rdr\frac{A^2}{r^2}\\
&=2\pi A^2|\log r_0|=\frac{2\pi\delta^2r_0^2}{|\log r_0|}
\end{split}\ee
which gives
\be \int_{B(0,1)\backslash B(0,r_0)} |\nabla u|^2\ge \frac{2\pi\delta^2r_0^2}{|\log r_0|}\ee 
since $f$ is harmonic.$\hfill\Box$

\begin{cor}\label{level}
For all $0<\epsilon<1$, there exists $C=C(\epsilon)>100$ such that if $0<r_0<1,$ $u\in C^1(\ B(0,1)\backslash B(0,r_0),\R)$ and
\be u|_{\partial B(0,r_0)}>\delta r_0-\frac{\delta r_0}{C}\ and\ \ u|_{\partial B(0,1)}<\frac{\delta r_0}{C}\ee 
then \be \int_{B(0,1)\backslash B(0,r_0)} |\nabla u|^2\ge\epsilon \frac{2\pi\delta^2r_0^2}{|\log r_0|}.\ee 
\end{cor}

\nd We will use the following lemma:

\begin{lem} \label {1}Let $0<r<1$, let $f, g$ be two harmonic functions on $B(0,1))\backslash \overline{B(0,r)}$, with $g|_{\partial B(0,1)}=a<b=g|_{\partial B(0,r)}$, and $f\le g$ on $\partial B(0,1)$, $f\ge g$ on $\partial B(0,r)$. Then
\be \int_{B(0,1)\backslash B(0,r)}|\nabla f|^2\ge\int_{B(0,1)\backslash B(0,r)}|\nabla g|^2.\ee 
\end{lem}

Let us admit this lemma for the moment and use it to prove Corollary 7.36. For each $C$, set $r=r_0,\ f=u,$ and $g$ such that $g|_{\partial B(0,1)}=\frac\delta C r_0,\ g|_{\partial B(0,r)}=(1-\frac 1C)\delta r_0$. Then we get
\be \int_{B(0,1)\backslash B(0,r_0)} |\nabla u|^2\ge\int_{B(0,1)\backslash B(0,r_0)} |\nabla g|^2=(1-\frac 2C)^2\frac{2\pi \delta^2r_0^2}{|\log r_0|}\ee 
and for each $\epsilon<1$ we can find $C$ large enough such that $(1-\frac 2C)^2\ge\epsilon$, thus complete the proof of Corollary 7.36.

Now let us prove Lemma 7.39. Set $h=(f-g)\nabla(f+g)$; then \be {\rm div} h= |\nabla f|^2-|\nabla g|^2\ee  
since $\Delta f=\Delta g=0$. Denote by $U=B(0,1)\backslash B(0,r)$, then by Stokes formula:
\be\begin{split}
\int_U|\nabla f|^2-|\nabla g|^2&=\int_U div h=\int_{\partial U} h\cdot \vec n=\int_{\partial U} (f-g)\frac{\partial}{\partial\vec n}(f+g)\\
&=\int_{\partial U}(f-g)\frac{\partial}{\partial\vec n} (f-g)+2\int_{\partial U}(f-g)\frac{\partial}{\partial\vec n}g,
\end{split}\ee
where $\vec n$ is exterior unit normal vector.

For the first term, since $k=f-g$ is harmonic, by Green's formula 
\be \int_{\partial U} k\frac{\partial}{\partial\vec n}k=\int_U(\nabla k\cdot \nabla k)+\int_U (k\Delta k)=\int_U|\nabla k|^2\ge 0.\ee 

For the second, by the boundary condition, on $\partial B(0,1)$, $f-g\le 0$ and $\frac{\partial}{\partial\vec n} g<0$, hence $(f-g)\frac{\partial}{\partial\vec n}g\ge 0$; similar for $\partial B(0,r)$, thus we have
\be \int_{\partial U}(f-g)\frac{\partial}{\partial\vec n}g\ge 0.\ee 

As a result \be \int_U|\nabla f|^2-|\nabla g|^2\ge 0.\ee  $\hfill\Box$

\bigskip

\section{Conclusion}

After all the preparations in the previous sections, we are going to conclude for Theorem \ref{main} in this section.

So fix a $\e<\e_0$. We are going to estimate the Hausdorff measure of $E_k$ for $k$ large enough. For each $k$ fixed, we have chosen $o_k$ and $r_k$ as in Proposition 5.11. Then by Proposition 6.1(1), $E_k\cap D_k(0,\frac {39}{40})\bs D_k(o_k,\frac{1}{10}r_k)$ is composed of two disjoint pieces $G_k,i=1,2$ such that (6.2) and (6.3) hold, where we replace $q_n,s_n$ by $o_k,r_k$. Moreover we can also suppose that $r_k<2^{-5}$, since $k$ is large.

\begin{pro}\label{epsilondelta}For all $\epsilon>0$, there exists $0<\delta=\d(\e)<\e$ and $\theta_0=\theta_0(\e)<\frac\pi2$, that depend only on $\e$, and satisfy the following properties. If $\frac\pi2>\theta>\theta_0$ (which means $\theta=(\theta_1,\theta_2)$ with $\theta_0<\theta_1\le\theta_2<\frac\pi2,\ i=1,2$) and $E$ is minimal in $B(0,1)$
which is also $\delta$ near 
$P_\theta=P^1_\theta\cup_\theta P^2_\theta$ in $B(0,1)\backslash B(0,\frac 12)$, and if moreover
\be p_0^i(E)\supset P_0^i\cap B(0,\frac 34)\ee
where $p_0^i$ denotes the orthogonal projection on $P_0^i,i=1,2$,
then $E$ is $\epsilon$ near $P_\theta$ in $B(0,1)$.
\end{pro}

\nd

We prove it by contradiction. So suppose the proposition is not true. Then there exists $\epsilon>0$, two sequences $\delta_l\to 0$ and $\theta_l\to(\frac\pi2,\frac\pi2)$, and a sequence of minimal sets $E_l$ in $B(0,1)$ such that $E_l$ is $\delta_l$ near $P_l=P_{\theta_l}$ in $B(0,1)\backslash B(0,\frac 12)$, and
\be p_0^i(E_l)\supset P_0^i\cap B(0,\frac34),\ee
but $E_l$ is not $\epsilon$ near $P_l$ in $B(0,1)$.

Since $P_l$ converges to $P_0=P_0^1\cup_\perp P_0^2$, $\delta_l\to 0$ and $\e$ is fixed, there exists a sequence $\{a_l\}$ which converges to 0, such that $E_l$ is $a_l$ near $P_0$ in $B(0,1)\backslash B(0,\frac 12)$, but not $\frac\epsilon 2$ near $P_0$ in $B(0,1)$ (because $E_l$ is not $\e$ near $P_l$, and $P_l$ is $\frac\e2$ near $P_0$ when $l$ is large.) 

Modulo extracting a subsequence, we can suppose that $\{E_l\}$ converges to a limit $E_\infty$. Then $E_\infty\cap B(0,1)\bs B(0,\frac 12)=P_0\cap B(0,1)\bs B(0,\frac12)$. 

We want to prove that
\be H^2(E_l\cap D(0,\frac 34))<\frac98\pi+b_l,\ with\ b_l\to 0\ when\ l\to\infty,\ee 
where $D(x,r)$ denotes $D_0(x,r)$ for short.

In fact, since $E_l$ is very near $P_0$ in $B(0,1)\bs B(0,\frac12)$ when $l$ is large, by the $C^1$ regularity of minimal sets (c.f. Thm \ref{c1}), we know that $E_l\cap \partial D(0,\frac 34)=\Gamma^1_l\cup\Gamma^2_l$, where $\Gamma_l^1$ and $\Gamma_l^2$ are two disjoint curves, $\Gamma^i_l$ is the graph of a $C^1$ function $h^i_l\ :\ P_0^i\cap \partial D(0,\frac 34)\to {P_0^i}^\perp$, with 
\be ||h_l^i||_\infty\to 0, k\to\infty\ et\ ||\frac{d}{dx}h_l^i||_\infty\le 1,\ \forall l.\ee 

Now set $D_l^i=(P_0^i\cap D(0,\frac34))\cup A_l^i$, where $A_l^i=\{(x,y): x\in P_0^i\cap \partial D(0,\frac 34), y\in[0,h_l^i(x)]\}$, which is just a 2-dimensional thin surface between $P_0^i\cap \partial D(0,\frac 34)$ and $\Gamma_l^i$, the graph of $h^i_l$. We can also write 
\be A_l^i=\{(x, th_l^i(x))\ :\ x\in P_0^i\cap\partial D(0,\frac34),t\in[0,1]\}.\ee 

It is not hard to see that $D_l^i$ is a surface whose boundary is $\Gamma_l^i$, and $D_l=D_l^1\cup D_l^2$ contains a deformation of $E_l$ in $\overline D(0,\frac 34).$

Let us verify that
\be H^2(A_l^i)\le \frac{3\sqrt 2\pi}{2} ||h_l^i||_\infty .\ee  
Since $h_l^i$ is $\R^2$-valued, the simplest is to use a parameterization.
Set $g_l^i(x,t):\partial D(0,\frac 34)\times[0,1]\to \R^4\ ; g_l^i(x,t)=(x, th_l^i(x))$, then $A_l^i$ is its image. And 
\be \frac{\partial}{\partial x}g_l^i=(1,t\frac{d}{dx}h_l^i(x))\ ;\ \frac{\partial}{\partial t}g_l^i=(0, h_l^i(x))\ee 
Therefore \be |\overrightarrow{\frac{\partial}{\partial x}g_l^i}\times\overrightarrow{\frac{\partial}{\partial t}g_l^i}|
=\sqrt{(1+t^2|\frac{\partial}{\partial x}h_l^i|^2)(|h_l^i|^2)-(t\overrightarrow{\frac{\partial}{\partial x}h_l^i}\cdot \overrightarrow{h_l^i}})^2.\ee 
But we know that $|\frac{\partial}{\partial x}h_l^i|\le 1$, hence
\be |\overrightarrow{\frac{\partial}{\partial x}g_l^i}\times\overrightarrow{\frac{\partial}{\partial t}g_l^i}|\le \sqrt{(1+t^2)|h_l^i|^2}\le\sqrt 2 |h_l^i|,\ee 
and therefore
\be\begin{split}
H^2(A_l^i)&=\int_{\partial D(0,\frac 34)\times[0,1]}|\overrightarrow{\frac{\partial}{\partial x}g_l^i}\times\overrightarrow{\frac{\partial}{\partial t}g_l^i}|dtdx\\
&\le\sqrt 2\int_{\partial D(0,\frac 34)}|h_l^i|\le \frac{3\sqrt 2}{2}\pi ||h_l^i||_\infty.
\end{split}\ee

Thus we obtain
\be H^2(D_l^i)\le \frac {9}{16}\pi+\frac{3\sqrt 2}{2}\pi ||h_l^i||_\infty\ee 
and
\be H^2(D_l)\le \frac 98\pi+\frac{3\sqrt 2}{2}\pi(||h_l^1||_\infty+||h_l^2||_\infty).\ee 
But $E_l$ is locally minimal, hence 
\be H^2(E_l\cap D(0,\frac 12))\le H^2(D_l)=\frac 98\pi+\frac{3\sqrt 2}{2}\pi(||h_l^1||_\infty+||h_l^2||_\infty).\ee 

Take $b_l=\frac{3\sqrt 2}{2}\pi(||h_l^1||_\infty+||h_l^2||_\infty)$, thus we get (8.4), since $||h_l^1||_\infty+||h_l^2||_\infty$ converges to $0$ when $l\to\infty$.

On the other hand, $E_l$ converges to $P_0$ in $B(0,1)\bs B(0,\frac12),$ so we have
\be E_l\cap (B(0,1)\backslash B(0,\frac 12))\to P_0\cap (B(0,1)\backslash B(0,\frac 12))=E_\infty\cap (B(0,1)\backslash B(0,\frac 12)).\ee 
By the lower semi continuity of Hausdorff measures for minimal sets (\cite{GD03} Thm 3.4, and recall that $E_\infty$ is the limit of $E_l$)
\be H^2(E_\infty\cap D(0,\frac 34))\le \liminf_{k\to\infty} H^2(E_l\cap D(0,\frac 34)).\ee
So we have
\be\begin{split}H^2(E_\infty)&=H^2(E_\infty\cap (B(0,1)\backslash D(0,\frac 34))+H^2(E_\infty\cap D(0,\frac 34))\\
&\le H^2(P_0\cap (B(0,1)\backslash D(0,\frac 34)))+\liminf_{k\to\infty} H^2(E_l\cap D(0,\frac 34))\\
&=2\pi+\liminf_{k\to\infty} b_l=2\pi.
\end{split}\ee

Now by (8.3) and the fact that $E_\infty$ is the limit of $E_l$, we know that 
\be p_0^i(E_\infty)\supset P_0^i\cap B(0,\frac 34).\ee
By hypothesis, $E_l$ converges to $P_0$ in $B(0,1)\bs B(0,\frac12)$, hence $E_\infty\cap B(0,1)\bs B(0,\frac12)=P_0\cap B(0,1)\bs B(0,\frac 12)$, and therefore
\be p_0^i (E_\infty)\supset p_0^i(E_\infty\cap B(0,1)\bs B(0,\frac 12))=P_0^i\cap B(0,1)\bs B(0,\frac 12).\ee
Thus we have
\be p_0^i(E_\infty)\supset P_0^i\cap B(0,1),\ee
and
\be E_\infty\cap \partial B(0,1)=P_0\cap \partial B(0,1)\ee
because $E_l$ converges to $P_0$ in $B(0,1)\bs B(0,\frac 12)$. Then by Theorem \ref{unicite}, (8.17), (8.20) and (8.21) give that
\be E_\infty=P_0.\ee 

This is impossible, because $E_l$ is $\frac\epsilon 2$ far from $P_0$. Thus we complete the proof of Proposition 8.1.$\hfill\Box$

\medskip

Now for $0<\theta=(\theta_1,\theta_2)$ with $0<\theta_1\le\theta_2<\frac\pi2$, set, for each $x\in\R^4,r>0$, 
\be D_\theta(x,r)=x+\{{p_\theta^1}^{-1}[B(0,r)\cap P_\theta^1]\cap {p_\theta^2}^{-1}[B(0,r)\cap P_\theta^2]\},\ee
where $P_\theta=P_\theta^1\cup P_\theta^2$ is the union of two planes with characteristic angles $\theta_1\le\theta_2$, and denote by $p_\theta^i$ the orthogonal projection to $P_\theta^i,i=1,2$. Then we have

\begin{cor}
For all $\e>0$, there exists $0<\d<\e$ and $0<\theta_0<\frac\pi2$, which do not depend on $\e$, with the following properties. If $\theta_0<\theta<\frac\pi2$, and if $E$ is minimal in $D_\theta(0,1)$ and is $\delta$ near $P_\theta$ in $D_\theta(0,1)\bs D_\theta(0,\frac 14)$, and moreover
\be p_\theta^i(E)\supset P_\theta^i\cap B(0,\frac 34),\ee
then $E$ is $\e$ near $P_\theta$ in $D_\theta(0,1)$.
\end{cor}

\nd First observe that there exists $0<\phi<\frac\pi2$ such that for all $0<\theta<\phi$ we have 
\be B(x,r)\subset D_\theta(x,r)\subset B(x,2r).\ee

Then for $\e>0$, take $\d=\d(\e)$ and $\theta_0=\max\{\phi,\theta_0(\e)\}$, where $\theta_0(\e)$ and $\d(\e)$ are as in Proposition 8.1. Then if $E$ is $\d$ near $P_\theta$ in $D_\theta(0,1)\bs D_\theta(0.\frac 14)$, it is $\d$ near $P_\theta$ in $B(0,1)\bs B(0,\frac 12)$. By Proposition 8.1, $E$ is $\e$ near $P_\theta$ in $B(0,1)$. Hence $E$ is $\e=\max\{\d,\e\}$ near $P_\theta$ in
$B(0,1)\cup [D_\theta(0,1)\bs D_\theta(0,\frac 14)]=D_\theta(0,1)$. \qed

\bigskip

\noindent\textbf{Proof of Theorem \ref{main}.} Take all the notation at the beginning of this section.

 Fix a $k$ large, and denote by $D(x,r)=D_k(x,r),C^i(x,r)=C_k^i(x,r)$ for $i=1,2$, and $d_{x,r}=d_{x,r}^k$.

We know that in $D(o_k,r_k)$, $E_k$ is not $\e r_k$ near any translation of $P_k$, so by Corollary 8.24, $E_k$ is not $\d r_k$ near any translation of $P_k$ in $D(o_k,r_k)\bs D(o_k,\frac 14r_k)$. However by (6.2), $E_k\cap D(o_k,r_k)\bs D(o_k,\frac 14r_k)=[G^1\cup G^2]\cap [D(o_k,r_k)\bs D(o_k,\frac 14r_k)]$, where $G^i$ is a $C^1$ graph of $P_k^i\cap D(0,\frac{39}{40})\bs D(o_k,\frac{1}{10}r_k)$. Hence there exists $i\in\{1,2\}$ such that in $D(o_k,r_k)\bs D(o_k,\frac14 r_k)$, $G^i$ is not $\delta r_k$ near any translation of $P_k^i$. Without loss of generality, we can suppose this is the case for $i=1$. 

Then denote by $P=P_k^1$ for short, and let $g^1$ be as in (6.2); then $g^1$ is a map from $P$ to $P^\perp$, and is therefore from $\R^2$ to $\R^2$. Write $g^1=(\varphi_1,\varphi_2)$, where $\varphi_i:\R^2\to \R$. Then since the graph of $g^1$ is $\d r_k$ far from all translation of $P$, there exists $j\in\{1,2\}$ such that
\be \sup_{x,y\in P\cap D(o_k, r_k)\backslash D(o_k, \frac 14 r_k)}|\varphi_j(x)-\varphi_j(y)|\ge \frac 12 r_k\delta.\ee

Suppose this is true for $j=1$. Denote by 
\be K=\{(z, \varphi_1(z)): z\in (D(0,\frac 34)\backslash D(o_k,\frac14 r_k))\cap P\},\ee 
then
\be \begin{split}K\mbox{ is the orthogonal }&\mbox{projection of }G^1\cap D(0,\frac 34)\\
&\mbox{ on a 3-dimensional subspace of }\R^4.\end{split}\ee 

%%$\d-$ proche------$\d r$proche

For $\frac 14r_k\le s\le r_k$, define 
\be \Gamma_s=K\cap p^{-1}(\partial D(o_k, s)\cap P)=\{(x,\varphi_1(x))|x\in\partial D(o_k, s)\cap P\}\ee  the graph of $\varphi_1$ on $\partial D(o_k, s)\cap P$. 

We know that the graph of $\varphi_1$ is $\frac12\d_k$ far from $P$ in $D(o_k,r_k)\bs D(o_k,\frac 14r_k)$; then there are two cases:

1st case: there exists $t\in[\frac 14r_k, r_k]$ such that
\be \sup_{x,y\in \Gamma_t}\{|\varphi_1(x)-\varphi_1(y)|\}\ge \frac{\delta}{C}r_k,\ee 
where $C=4C(\frac 12)$ is the constant of Corollary \ref{level}.

Then there exists $a,b\in \Gamma_t$ such that $|\varphi_1(a)-\varphi_1(b)|>\frac\delta C r_k\ge\frac\delta Ct$. Since $||\nabla \varphi_1||_\infty\le ||\nabla \varphi||_\infty<1$, we have 
\be \int_{\Gamma_t}|\varphi_1-m(\varphi_1)|^2\ge \frac{t^3\delta^3}{4C^3}=(\frac43 t\delta)^3(\frac{27}{4^4C^3}).\ee 
Now in $D(0,\frac 34)$ we have $d(0,o_k)<6\e\le 10\e \cdot\frac34$, and $s<r_k<\frac 18<\frac 12\times\frac 34$, therefore we can apply Corollary \ref{cas1} and obtain
\be \int_{(D(0,\frac 34)\backslash D(o_k, t))\cap P}|\nabla \varphi_1|^2\ge C_1(\delta)t^2.\ee 

2nd case: for all $\frac 14r_k\le s\le r_k$, 
\be \sup_{x,y\in \Gamma_s}\{|\varphi_1(x)-\varphi_1(y)|\}\le \frac{\delta}{C}r_k.\ee 
However, since
\be\begin{split}
\frac12 r_k\delta&\le\sup\{|\varphi_1(x)-\varphi_2(y)|:x,y\in P\cap D(o_k,r_k)\backslash D(o_k,\frac14 r_k)\}\\
&=\sup\{|\varphi_1(x)-\varphi_2(y)|:s,s'\in[\frac14r_k,r_k],x\in\Gamma_s,y\in\Gamma_{s'}\},
\end{split}\ee
there existe $\frac 14r_k\le t<t'\le r_k$ such that
\be \sup_{x\in \Gamma_{t},y\in\Gamma_{t'}}\{|\varphi_1(x)-\varphi_1(y)|\}\ge \frac12r_k\delta. \ee 
Fix $t$ and $t'$, and without loss of generality, suppose that
\be \sup_{x\in \Gamma_{t},y\in\Gamma_{t'}}\{\varphi_1(x)-\varphi_1(y)\}\ge \frac12r_k\delta .\ee 
Then
\be \inf_{x\in\Gamma_{t}}\varphi_1(x)-\sup_{x\in\Gamma_{t'}}\varphi_1(x)\ge \frac12r_k\delta-2\frac\delta Cr_k=(1-\frac{2}{C(\frac12)})\frac{\delta}{2}r_k\ge (1-\frac{2}{C(\frac12)})\frac{\delta}{2}t'\ee 
because $C=4 C(\frac 12)$.

Now look at what happens in the ball $D(o_k,t')\cap P$. Apply Corollary \ref{level} to the scale $t'$, we get
\be \int_{(D(o_k,t')\backslash D(o_k, t))\cap P}|\nabla \varphi_1|^2\ge C(\delta,\frac 12)\frac{\pi(\frac\delta2)^2 t'^2}{\log\frac{t'}{t}}.\ee 
Then since $\frac{t'}{t}\le 4, t'>t$, we have 
\be \int_{((D(o_k,t')\backslash D(o_k, t))\cap P}|\nabla \varphi_1|^2\ge C_2(\delta)t^2.\ee 

In both cases we pose $t_k=t$. The discussion above yields that there exists a constant $C_0(\delta)=\min\{C_1(\delta), C_2(\delta)\}$, which depends only on $\delta$, such that
\be \int_{(D(0,\frac 34)\backslash D(o_k, t_k))\cap P}|\nabla \varphi_1|^2\ge C_0(\delta)t_k^2.\ee 

On the other hand, since $|\nabla \varphi_1|\le|\nabla g^1|<1$,
\be \sqrt{1+|\nabla\varphi_1|^2}>\sqrt{1+\frac 12|\nabla\varphi_1|^2+\frac{1}{16}|\nabla\varphi_1|^4} =1+\frac14|\nabla\varphi_1|^2.\ee 
Hence 
\be\begin{split}
H^2(K\bs C_k^1(o_k,t_k) )&=\int_{D(0,\frac 34)\backslash C_k^1(o_k,t_k)\cap P}\sqrt{1+|\nabla\varphi_1|^2}\ge\int_{D(0,\frac 34)\backslash C_k^1(o_k, t_k)\cap P}1+\frac14|\nabla\varphi_1|^2\\
&\ge H^2((D(0,\frac 34)\backslash C_k^1(o_k, t_k))\cap P))+\frac 14\int_{D(0,\frac 34)\backslash C_k^1(o_k, t_k)\cap P}|\nabla \varphi_1|^2\\
&= H^2((D(0,\frac 34)\backslash C_k^1(o_k, t_k))\cap P_k^1))+C_0(\d)t_k^2.
\end{split}\ee

Then by (8.29) we get
\be\begin{split} H^2&(G^1\cap D(0,\frac34)\bs D(o_k,t_k))\ge H^2(K\bs D(o_k,t_k))\\
&\ge H^2((P_k^1+o_k)\cap D(0,\frac 34)\backslash D(o_k, t_k))+C(\d)t_k^2\\
&=H^2(P_k^1\cap D(0,\frac 34)\backslash D(0, t_k))+C_0(\d)t_k^2.\end{split}\ee

Thus we obtain an estimate for the regular part of $E_k$. Next for all other parts of $E_k$ we will control their measures by projection.

So let us decompose $E_k$. Set $F_1=E_k\cap D(o_k,t_k)$, $F_2=G^2_{t_k}$, $F_3=G^1_{t_k}\bs D(0,\frac 34)$, and $F_4=G^1_{t_k}\cap D(0,\frac34)$, where $G^i_t$ is defined as in Proposition 6.1(2). Then $F_i$ are disjoint.

For $F_1$, by Proposition 2.19 and Lemma \ref{projection}
\be (1+2\cos\theta_k(1))H^2(F_1)\ge H^2(p^1_k(F_1))+H^2(p_k^2(F_1)),\ee
where $\theta_k=(\theta_k(1),\theta_k(2))$ with $\theta_k(1)\le\theta_k(2)$. (Recall that $E_k$ has the same boundary as $P_k=P^1\cup_{\theta_k}P^2$ with $\theta_k\ge\frac\pi2-\frac1k$).
However since $\theta_k\ge \frac\pi2-\frac1k$, we have
\be H^2(F_1)\ge (1-\frac3k)[H^2(p^1_k(F_1))+H^2(p_k^2(F_1))]\ee
when $k$ is large. On the other hand, by Proposition 6.1(4),  
\be p_k^i(F_1)\supset P_k^i\cap C_k^i(o_k,t_k).\ee

As a result
\be \begin{split}H^2(F_1)&\ge (1-\frac3k)H^2((P_k+o_k)\cap D(o_k,t_k))\\
&\ge H^2((P_k+o_k)\cap D(o_k,t_k))-\frac{C}{k}t_k^2=H^2(P_k\cap D(0,t_k))-\frac{C}{k}t_k^2.\end{split}\ee

For $F_2$, by Proposition 6.1(2), we have 
\be p_k^2(F_2)=p_k^2(G_{t_k}^2)\supset P_k^i\cap D(0,1)\bs C_k^2(o_k,t_k).\ee
Hence 
\be H^2(F_2)\ge H^2[P_k^2\cap D(0,1)\bs C_k^2(o_k,t_k)]=H^2[P_k^2\cap D(0,1)\bs D(0,t_k)].\ee

For $F_3$, still by Proposition 6.1(2), and by the definition of $F_3$, we have  
\be p_k^1(F_3)\supset p_k^1(G_{t_k}^1)\backslash p_k^1(D(0,\frac 34)\supset P_k^1\cap D(0,1)\bs D(0,\frac34).\ee
Hence
\be H^2(F_3)\ge H^2(P_k^1\cap D(0,1)\bs D(0,\frac 34).\ee

For the last part, the definition of $F_4$ gives
\be H^2(F_4)=H^2(G^1\bs D(o_k,t_k))\ge H^2(P_k^1\cap D(0,\frac 34)\backslash D(0, t_k))+C_0(\d)t_k^2.\ee

Now we add measures of these four pieces together and get 
\be\begin{split} H^2(E_k)&=H^2(F_1)+H^2(F_2)+H^2(F_3)+H^2(F_4)\\ 
&\ge H^2(P_k\cap D(0,1))+t_k^2(C_0(\d)-\frac Ck).\end{split}\ee

Then when $k$ is such that $C(\d)>\frac Ck$, we have
\be H^2(E_k)>H^2(P_k\cap D(0,1)),\ee
which contradicts Proposition 4.8(4). Thus the proof of Theorem \ref{main} is completed. $\hfill\Box$

\bigskip

As a final remark, we give two similar theorems below.
 \begin{thm}[minimality of the union of $n$ almost orthogonal $m-$dimensional planes]For each $m\ge 2$ and $n\ge 2$, there exists $0<\theta<\frac\pi2$, such that if $P^1,P^2,\cdots,P^n$ are $n$ planes of dimension $m$ in $\R^{nm}$ with characteristic angles $\alpha^{ij}=(\alpha^{ij}_1,\alpha^{ij}_2,\cdots,\a^{ij}_m)$ between $P^i$ and $P^j,1\le i<j\le n$, which verify $\theta<\a^{ij}_1\le\a^{ij}_2\le\cdots\le\a^{ij}_m\le\frac\pi2$ for all $1\le i<j\le n$, then their union $\cup_{i=1}^n P^i$ is a minimal cone.
\end{thm}
\begin{thm}[minimality of the union of a plane and a $\Y$ set which are almost orthogonal]There exists $0<\theta<\frac\pi2$, such that if $P$ and $Q$ are two subspaces of $\R^5$ of dimension 2 and 3 respectively, which verify
\be \mbox{for all unit simple vectors }u\in P,v\in Q\mbox{, then the angle between }u\mbox{ and }v\mbox{ is larger than }\theta,\ee
and if $Y$ is a $\Y$ set in $Q$ centered at the point of $P\cap Q$, 
then the union $P\cup Y$ is a minimal cone.
\end{thm}

The general idea for the proofs of the two theorems are somehow similar to the proof of Theorem \ref{main}. But there are also some non trivial modifications, especially for Theorem 8.57. See \cite{XY10} for more detail of the proofs. We could have continued to discuss one by one that if the almost orthogonal unions of other pairs or families of minimal cones are minimal. In fact Theorem 8.57 was the first try, where we replaced one plane by the simplest minimal cone $\Y$. But then the author noticed that the proof became already a bit too complicated, so she stopped at this step.

\renewcommand\refname{References}
\bibliographystyle{plain}
\bibliography{reference}
\end{document}